\documentclass[aop,citesort,MSNbibl,dvips]{arximspdf}
\usepackage{mathbh}
\usepackage{graphicx}

\doi{10.1214/10-AOP549}
\volume{39}
\issue{1}
\pubyear{2011}
\firstpage{1}
\lastpage{69}

\makeatletter

\renewcommand{\mid}{|}

\newproclaim{Remark}{Remark}
\newtheorem{theorem}{Theorem}
\newtheorem{lmm}{Lemma}
\newtheorem{prp}{Proposition}
\newproclaim{Assumption}{Assumption}

\newcommand{\dgr}{d_{\mathrm{gr}}}

\newcommand{\bq}{\mathbf{q}}

\newcommand{\bN}{\mathbf{N}}
\newcommand{\mm}{\mathbf{m}}

\newcommand{\g}{\mathcal{G}}
\newcommand{\R}{\mathbb{R}}
\newcommand{\K}{\mathbb{K}}
\newcommand{\N}{\mathbb{N}}
\newcommand{\F}{\mathbb{F}}
\newcommand{\Z}{\mathbb{Z}}
\newcommand{\Q}{\mathbb{Q}}
\newcommand{\E}{\mathbb{E}}
\newcommand{\D}{\mathbb{D}}
\newcommand{\FF}{\mathcal{F}}
\newcommand{\m}{\mathcal{M}}
\newcommand{\vf}{{\varphi}}
\newcommand{\tc}{\mathcal{T}^\circ}
\newcommand{\vep}{\varepsilon}
\newcommand{\ga}{\gamma}

\newcommand{\eps}{\varepsilon}
\newcommand{\ov}{\overline}
\newcommand{\wh}{\widehat}
\newcommand{\wt}{\widetilde}
\newcommand{\la}{\longrightarrow}
\newcommand{\da}{\downarrow}
\newcommand{\var}{\operatorname{Var}}
\newcommand{\ind}{\mathbh{1}}

\makeatother

\begin{document}
\begin{frontmatter}

\title{Scaling limits of random planar maps with large~faces}
\runtitle{Scaling limits of planar maps}

\begin{aug}
\author[A]{\fnms{Jean-Fran\c cois} \snm{Le Gall}\corref{}\ead[label=e1]{jean-francois.legall@math.u-psud.fr}} and
\author[A]{\fnms{Gr\'egory} \snm{Miermont}\ead[label=e2]{gregory.miermont@math.u-psud.fr}}
\runauthor{J.-F. Le Gall and G. Miermont}
\affiliation{Universit\'e Paris-Sud}
\address[A]{D\'epartement de Math\'ematiques\\
Universit\'e Paris-Sud\\
Centre d'Orsay\\
91405 ORSAY Cedex\\
France\\
\printead{e1}\\
\phantom{E-mail: }\printead*{e2}}
\end{aug}

\received{\smonth{7} \syear{2009}}
\revised{\smonth{3} \syear{2010}}

%
\begin{abstract}
We discuss asymptotics for large random planar maps under the
assumption that the distribution of the degree of a typical face is
in the domain of attraction of a stable distribution with index
$\alpha\in(1,2)$. When the number $n$ of vertices of the map tends
to infinity, the asymptotic behavior of distances from a
distinguished vertex is described by a random process called the
continuous distance process, which can be constructed from a
centered stable process with no negative jumps and index $\alpha$.
In particular, the profile of distances in the map, rescaled by the
factor $n^{-1/2\alpha}$, converges to a random measure defined in
terms of the distance process. With the same rescaling of distances,
the vertex set viewed as a metric space converges in distribution as
$n\to\infty$, at least along suitable subsequences, toward a
limiting random compact metric space whose Hausdorff dimension is
equal to $2\alpha$.
\end{abstract}

%

\end{frontmatter}

\section{Introduction}\label{sec1}

The goal of the present work is to discuss the continuous limits of
large random planar maps when the distribution of the
degree of a typical face has a heavy tail. Recall that a planar map
is a proper embedding of a finite connected graph in the
two-dimensional sphere. For technical reasons, it is convenient to
deal with rooted planar maps, meaning that there is a distinguished
oriented edge called the \textit{root edge}. One is interested in the
``shape'' of the graph and not in the particular embedding that is
considered. More rigorously, two rooted planar maps are identified if
they correspond via an orientation-preserving homeomorphism of the
sphere. The faces of the map are
the connected components of the complement of edges and the degree of
a face counts the number of edges that are incident to it. Large
random planar graphs are of particular interest in theoretical
physics, where they serve as models of random geometry \cite{ADJ}.

A simple way to generate a large random planar map is to choose it
uniformly at random from the set of all rooted $p$-angulations with
$n$ faces (a planar map is a $p$-angulation if all faces have degree
$p$). It is conjectured that the scaling limit of uniformly
distributed $p$-angulations with $n$ faces, when $n$ tends to infinity
(or, equivalently, when the number of vertices tends to infinity), does
not depend on the choice of $p$ and is given by the so-called Brownian
map. Since the pioneering work of Chassaing and Schaeffer
\cite{CSise}, there have been several results supporting this
conjecture. Marckert and Mokkadem \cite{MaMo} introduced the Brownian
map and proved a weak form of the convergence of rescaled uniform
quadrangulations toward the Brownian map. A stronger version,
involving convergence of the associated metric spaces in the sense of
the Gromov--Hausdorff distance, was derived in Le Gall \cite{legall06}
in the case of uniformly distributed $2p$-angulations. Because the
distribution of the
Brownian map has not been fully characterized, the convergence results
of \cite{legall06} require the extraction of suitable subsequences. Proving
the uniqueness of the distribution of the Brownian map is one of the
key open problems in this area.

A more general way of choosing a large planar map at random is to use
Boltzmann distributions. In this work, we restrict our attention to
bipartite maps, where all face degrees are even. Given a sequence
$q=(q_1,q_2,q_3,\ldots)$ of nonnegative real numbers and a bipartite
planar map $\mathbf{m}$, the associated Boltzmann weight is
%
\begin{equation}
\label{Bolt-weight}
W_q(\mm) = \prod_{f\in F(\mm)} q_{\operatorname{deg}(f)/2},
\end{equation}
where $F(\mm)$ denotes the set of all faces of $\mm$ and
$\operatorname{deg}(f)$ is the degree of the face~$f$. One can then generate a
large planar map by choosing it at random from the set of all planar
maps with $n$ vertices (or with $n$ faces) with probability weights
that are proportional to $W_q(\mm)$. Such distributions arise
naturally (possibly in slightly different forms) in problems involving
statistical physics models on random maps. This is discussed in
Section \ref{sec:some-motivation-from} below.

Assuming certain integrability conditions on the sequence of weights,
Marckert and Miermont \cite{jfmgm05} obtain a variety of limit
theorems for large random bipartite planar maps chosen according to
these Boltzmann distributions. These results are extended in Miermont
\cite{Mi1} and Miermont and Weill \cite{MW} to the nonbipartite case,
including large triangulations. In all of these papers, limiting
distributions are described in terms of the Brownian map. Therefore,
these results strongly suggest that the Brownian map should be the
universal limit of large random planar maps, under the condition that
the distribution of the degrees of faces satisfies some integrability
property. Note that, even though the distribution of the Brownian map
has not been characterized, many of its properties can be investigated
in detail and have interesting consequences for typical large planar
maps; see, in particular, the recent papers \cite{legall08} and
\cite{Mi2} (and Bettinelli \cite{betti}, for similar results, for
random maps
on surfaces of higher genus).

In the present work, we consider Boltzmann distributions such that,
even for large $n$, a random planar map with $n$ vertices will have
``macroscopic'' faces, which, in
some sense, will remain present in the scaling limit. This
leads to a (conjectured) scaling limit which is different from the
Brownian map. In fact, our limit theorems involve new random processes
that are closely related to the stable trees of \cite{duqleg02}, in
contrast to the construction of the Brownian map \cite{MaMo,legall06},
which is based on Aldous' continuum random tree (CRT).

Let us informally describe our main results, referring to the following
sections for more precise statements. For technical reasons, we
consider planar maps that are both rooted and pointed (in addition to
the root edge, there is a distinguished vertex, denoted by $v_*$).
Roughly speaking, we choose the Boltzmann weights $q_k$ in
(\ref{Bolt-weight}) in such a way that the distribution of the degree
of a (typical) face is in the domain of attraction of a stable
distribution with index $\alpha\in(1,2)$. This can be made more
precise by using the Bouttier--Di Francesco--Guitter bijection
\cite{BdFGmobiles} between bipartite planar maps and certain labeled
trees called \textit{mobiles}. A mobile is a (rooted) plane tree, where
vertices at even distance (resp., odd distance) from the root
are called white (resp., black) and white vertices are
assigned integer labels that satisfy certain simple rules; see
Section \ref{sec:mobiles-bouttier-di}. In the Bouttier--Di
Francesco--Guitter bijection, a (rooted and pointed) planar map $\mm$
corresponds to a mobile $\theta(\mm)$ in such a way that each face of
$\mm$ is associated with a black vertex of $\theta(\mm)$ and each
vertex of $\mm$ (with the exception of the distinguished vertex $v_*$)
is associated with a white vertex of $\theta(\mm)$. Moreover, the
degree of a face of $\mm$ is exactly twice the degree of the
associated black vertex in the mobile $\theta(\mm)$ (see Section
\ref{sec:mobiles-bouttier-di} for more details).

Under appropriate conditions on the sequence of weights $q$, formula
(\ref{Bolt-weight}) defines a finite measure $W_q$ on the set of all
rooted and pointed planar maps. Moreover, if $\mathbf{P}_q$ is the
probability measure obtained by normalizing $W_q$, then the mobile
$\theta(\mm)$ associated with a planar map $\mm$ distributed according
to $\mathbf{P}_q$ is a critical two-type Galton--Watson tree, with
different offspring distributions $\mu_0$ and $\mu_1$ for white and
black vertices, respectively, and labels chosen uniformly over all
possible assignments (see \cite{jfmgm05} and Proposition
\ref{mobile-map} below). The distribution $\mu_0$ is always geometric,
whereas $\mu_1$ has a simple expression in terms of the weights $q_k$.

We now come to our basic assumption. In the present work, we choose
the weights $q_k$ in such a way that $\mu_1(k)$ behaves like
$k^{-\alpha-1}$, when $k\to\infty$, for some
$\alpha\in(1,2)$. Recalling that the degree of a face of $\mm$ is
equal to twice the degree of the associated black vertex in the mobile
$\theta(\mm)$, we see that, in a certain sense, the face degrees of a
planar map distributed according to $\mathbf{P}_q$ are independent, with a
common distribution that belongs to the domain of attraction of a
stable law with index $\alpha$.

We equip the vertex set $V(\mm)$ of a planar map $\mm$ with the graph
distance $d_{\mathrm{gr}}$ and would like to investigate the properties
of this metric space when $\mm$ is distributed according to
$\mathbf{P}_q$ and conditioned to be large. For every integer $n\geq1$,
denote by $M_n$ a random planar map distributed according to
$\mathbf{P}_q(\cdot\mid\#V(\mm)=n)$. Our goal is to get information about
typical distances in the metric space $(V(M_n),d_{\mathrm{gr}})$ when $n$
is large and, if at all possible, to prove that these (suitably
rescaled) metric spaces converge in distribution as $n\to\infty$ in
the sense of the Gromov--Hausdorff distance. As a motivation for
studying the
particular conditioning $\{ \#V(\mm)=n\}$, we note that our results will
have immediate application to Boltzmann distributions on
\textit{nonpointed} rooted planar maps:
simply observe that a given rooted planar map with
$n$ vertices corresponds to exactly $n$ different rooted and pointed
planar maps.

To achieve the preceding goal, we use another nice feature of the Bouttier--Di
Francesco--Guitter bijection: up to an additive constant depending on
$\mm$, the distance between $v_*$ and an arbitrary vertex $v\in
V(\mm)\setminus\{v_*\}$ coincides with the label of the white vertex
of $\theta(\mm)$ associated with $v$. Thus, in order to understand the
asymptotic behavior of distances from $v_*$ in the map $M_n$, it
suffices to get information about labels in the mobile $\theta(M_n)$
when $n$ is large. To this end, we first consider the tree $\mathcal{T}(M_n)$
obtained by ignoring the labels in $\theta(M_n)$. Under our basic
assumption, the results of \cite{duqleg02} can be applied to prove
that the tree $\mathcal{T}(M_n)$ converges in distribution, modulo a rescaling
of distances by the factor $n^{-(1-1/\alpha)}$, toward the so-called
stable tree with index $\alpha$. The stable tree can be defined by a
suitable coding from the sample path of a centered stable L\'evy
process with no negative jumps and index $\alpha$, under an
appropriate excursion measure. The preceding convergence to the stable
tree is, however, not sufficient for our purposes since we are
primarily interested in labels. Note that, under the assumptions made in
\cite{jfmgm05} on the weight sequence $q$ (and, in particular, in the case
of uniformly distributed $2p$-angulations), the rescaled trees
$\mathcal{T}(M_n)$ converge toward the CRT and the scaling limit of
labels is
described in \cite{jfmgm05} as Brownian motion indexed by the
CRT or, equivalently, as the Brownian snake driven by a normalized
Brownian excursion. In our ``heavy tail'' setting, however, the scaling
limit of the
labels is \textit{not} Brownian motion indexed by the stable tree, but is
given by a new random process of independent interest, which we call
the \textit{continuous distance process}.

Let us give an informal presentation of the distance process---a
rigorous definition can be found in Section
\ref{sec:continuous-object} below. We view the stable tree as the
genealogical tree for a continuous population and the distance of a
vertex from the root is interpreted as its generation in the tree. Fix
a vertex $a$ in the stable tree. Among the ancestors of $a$, countably
many of them, denoted by $b_1,b_2,\ldots,$ correspond to a sudden
creation of mass in the population: each $b_k$ has a macroscopic
number $\delta_k>0$ of ``children'' and one can also consider the
quantity $r_k\in[0,\delta_k]$, which is the rank among these children
of the one that is an ancestor of $a$. The preceding description is
informal in our continuous setting (there are no children), but can be
made rigorous thanks to the ideas developed in \cite{duqleg02} and, in
particular, to the coding of the stable tree by a L\'evy process. We
then associate with each vertex $b_k$ a Brownian bridge
$(B_k(t))_{t\in[0,\delta_k]}$ (starting and ending at $0$) with
duration $\delta_k$, independently when $k$ varies, and we set
\[
D(a)=\sum_{k=1}^\infty B_k(r_k).
\]
The resulting process $D(a)$ when $a$ varies in the stable tree is the
continuous distance process. As a matter of fact, since vertices of
the stable tree are parametrized by the interval $[0,1]$ (using the
coding by a L\'evy process), it is more convenient to define the
continuous distance process as a process $(D_t)_{t\in[0,1]}$ indexed
by the interval $[0,1]$ (or even by $\R_+$ when we consider a forest
of trees).

Much of the technical work contained in this article is devoted to
proving that the rescaled labels in the mobile $\theta(M_n)$ converge
in distribution to the continuous distance process. The proper
rescaling of labels involves the multiplicative factor
$n^{-1/2\alpha}$ instead of $n^{-1/4}$, as in earlier work. This indicates
that the typical diameter of our random planar maps $M_n$ is of order
$n^{1/2\alpha}$, rather than $n^{1/4}$ in the case of maps with faces
of bounded degree. Because conditioning on the total number of
vertices makes the proof more difficult, we first establish a version
of the convergence of labels for a forest of independent mobiles
having the distribution of $\theta(\mm)$ under $\mathbf{P}_q$. The proof
of this result (Theorem \ref{functional}) is given in Section
\ref{forest-mobiles}. We then derive the desired convergence for the
conditioned objects in Section \ref{contour-conditioned}.

Finally, we obtain asymptotic results for the planar maps $M_n$ in
Section \ref{large-maps}. Theorem \ref{asymptodistance} gives precise
information about the profile of distances from the distinguished
vertex $v_*$ in $M_n$. Precisely, let $\rho_{M_n}^{(n)}$ be the
measure on $\R_+$ defined by
\[
\int\rho_{M_n}^{(n)}({d}x) \varphi(x)=\frac{1}{n}
\sum_{v\in V(M_n)} \varphi(n^{-1/2\alpha} d_{\mathrm{gr}}(v_*,v)).
\]
Then, the sequence of random measures $\rho_{M_n}^{(n)}$ converges in
distribution toward the measure $\rho^{(\infty)}$ defined by
\[
\int\rho^{(\infty)}({d}x) \varphi(x)
=\int_0^1 {d}t\, \varphi\bigl(c(D_t-\underline D)\bigr),
\]
where $c>0$ is a
constant depending on the sequence of weights and $\underline D
=\min_{t\in[0,1]} D_t$.

We also investigate the convergence of the suitably rescaled metric
spaces $V(M_n)$ in the Gromov--Hausdorff sense. Theorem \ref{GHconv}
shows that, at least along a subsequence, the random metric spaces
$(V(M_n),n^{-1/2\alpha}d_{\mathrm{gr}})$ converge in distribution toward a
limiting random compact metric space. Furthermore, the Hausdorff
dimension of this limiting space is a.s. equal to $2\alpha$,
which should be compared with the value $4$ for the
dimension of the Brownian map \cite{legall06}. The fact that the
Hausdorff dimension is bounded above by $2\alpha$ follows from
H\"older continuity properties of the distance process that are
established in Section \ref{sec:continuous-object}. The proof of the
corresponding lower bound is more involved and depends on some
properties of the stable tree and its coding by L\'evy processes,
which are investigated in \cite{duqleg02}. Similarly as in the case of
the convergence to the Brownian map, the extraction of a subsequence
in Theorem \ref{GHconv} is needed because the limiting distribution is
not characterized.

The paper is organized as follows. Section \ref{sec:discrete-picture} introduces Boltzmann
distributions on planar maps and formulates our basic assumption on
the sequence of weights. Section \ref{sec:coding-maps-with} recalls the Bouttier--Di
Francesco--Guitter bijection and the key result giving the distribution
of the random mobile associated with a planar map under the Boltzmann
distribution (Proposition \ref{mobile-map}). Section \ref{sec:coding-maps-with} also introduces
several discrete functions coding mobiles, in terms of which most of
the subsequent limit theorems are stated. Section \ref{sec:continuous-object} is devoted to the
definition of the continuous distance process and to its H\"older
continuity properties. In Section \ref{forest-mobiles}, we address the problem of the
convergence of the discrete label process of a forest of random
mobiles toward the continuous distance process of Section \ref{sec:continuous-object}. We then
deduce a similar convergence for labels in a single random mobile
conditioned to be large in Section \ref{contour-conditioned}. Section \ref{large-maps} deals with the
existence of scaling limits of large random planar maps and the
calculation of the Hausdorff dimension of limiting spaces. Finally,
Section \ref{sec:some-motivation-from} discusses some motivation coming from theoretical physics.

\subsection*{Notation} The symbols $K,K',K_1,K'_1,K_2,\ldots$
will stand for positive constants that may depend on the choice of the
weight sequence $q=(q_1,q_2,\ldots)$, but, unless otherwise indicated, do
not depend on other quantities. The value of these constants may vary
from one proof to another. The notation $C(\R)$ stands for the space
of all continuous functions from $\R_+$ into $\R$ and the notation
$\D(\R^d)$ stands for the Skorokhod space of all c\`adl\`ag functions
from $\R_+$ into $\R^d$. If $X=(X_t)_{t\geq0}$ is a process with c\`
adl\`ag
paths, $X_{s-}$ denotes the left limit of $X$ at $s$ for every $s>0$.
We denote the set of all finite
measures on $\R_+$ by $M_f(\R_+)$ and this set is equipped with the
usual weak topology. If $(a_k)$ and $(b_k)$ are two sequences of
positive numbers, the notation $a_k\sim b_k$ (as $k\to\infty$) means
that the ratio $a_k/b_k$ tends to $1$ as $k\to\infty$. Unless
otherwise indicated, all random variables and processes are defined on
a probability space $(\Omega,\mathcal{F},\mathbb{P})$.

\section{Critical Boltzmann laws on bipartite planar
maps}\label{sec:discrete-picture}

\subsection{Boltzmann distributions}

A rooted and pointed bipartite map is a pair $(\mm,v_*)$, where $\mm$
is a rooted bipartite planar map and $v_*$ is a
distinguished vertex of $\mm$. As in Section \ref{sec1}, the graph distance on
the vertex set $V(\mm)$ is
denoted by $\dgr$ and we let $e_-,e_+$ be, respectively, the origin and
the target of the root edge of $\mm$. By the bipartite nature of
$\mm$, the quantities $\dgr(e_+,v_*),\dgr(e_-,v_*)$ differ. Moreover,
this difference is at most $1$ in absolute value since $e_+$ and $e_-$
are linked by an edge. We say that $(\mm,v_*)$ is \textit{positive} if
\[
\dgr(e_+,v_*)=\dgr(e_-,v_*)+1 .
\]
It is called \textit{negative} otherwise,
that is, if $\dgr(e_+,v_*)=\dgr(e_-,v_*)-1$.

We let $\m_+^*$ denote the set of all rooted and pointed bipartite planar
maps that are positive. In the sequel, the mention of $v_*$ will
usually be
implicit, so we will simply denote the generic element of
$\m_+^*$ by $\mm$. For our purposes, it is useful to add an element
$\dagger$
to $\m_+^*$, which can be seen roughly as the \textit{vertex map} with no
edge and a single vertex $v_*$ ``bounding'' a single face of degree
$0$.

Let $q=(q_1,q_2,\ldots)$ be a sequence of nonnegative real numbers.
For every
$\mm\in\m_+^*\setminus\{\dagger\}$, set
\[
W_q(\mm) = \prod_{f\in F(\mm)} q_{\operatorname{deg}(f)/2},
\]
where $F(\mm)$ denotes the set of all faces of $\mm$. By
convention, we set $W_q(\dagger)=1$. This defines a
$\sigma$-finite measure on $\m_+^*$, whose total mass is
\[
Z_q=W_q(\m_+^*)\in[1,\infty] .
\]
We say that $q$ is \textit{admissible} if $Z_q<\infty$, in which case we
can define $\mathbf{P}_q=Z_q^{-1}W_q$ as the probability measure obtained by
normalizing $W_q$. The measure $\mathbf{P}_q$ is called the \textit{Boltzmann
distribution} on $\m_+^*$ with weight sequence $q$.

Following \cite{jfmgm05}, we have the following simple criterion for
the admissibility of~$q$. Introduce the function
%
\begin{equation}\label{eq:18}
f_q(x)=\sum_{k=1}^\infty N(k)q_k x^{k-1}
,\qquad x\geq0,
\end{equation}
where
\[
N(k)=\pmatrix{{2k-1}\cr{k-1}}.
\]
Let $R_q\geq0$ be the radius of convergence of this power
series. Note that by monotone convergence, the quantity
$f_q(R_q)=f_q(R_q-)\in[0,\infty]$ exists, as well as
$f'_q(R_q)=f'_q(R_q-)$.
\begin{prp}[\cite{jfmgm05}]
\label{sec:crit-boltzm-laws-2}
The sequence $q$ is admissible if and only if the equation
%
\begin{equation}
\label{eq:16}
f_q(x)=1-\frac{1}{x} ,\qquad x\geq1,
\end{equation}
has a solution. If this holds, then the smallest such solution
equals $Z_q$.
\end{prp}

On the interval $[0,R_q)$, the function $f_q$ is convex, so the
equation (\ref{eq:16}) has at most two solutions. Let us now pause for
a short
informal discussion, inspired by~\cite{jfmgm05}. For a ``typical''
admissible sequence $q$, the graphs of $f_q$ and of the function
$x\mapsto1-1/x$ will
cross at $x=Z_q$ without being tangent. In this case, the law of the
number of vertices of a $\mathbf{P}_q$-distributed random map will have an
exponential tail. An admissible sequence $q$ is called \textit{critical}
if the graphs are tangent at $Z_q$, that is, if
%
\begin{equation}
\label{eq:17}
Z_q^2f_q'(Z_q)=1 .
\end{equation}
For critical sequences, the law of the number of vertices of a
$\mathbf{P}_q$-distributed random map may have a tail heavier than
exponential. In the case where $R_q>Z_q$, \cite{jfmgm05}~shows that
this tail follows a power law with exponent $-1/2$. However, the law
of the degree of a typical face in such a random map will have an
exponential tail.

In the present paper, we will be interested in the ``extreme'' cases
where $q$ is a critical sequence such that $Z_q=R_q$. We will show
that in a number of these cases, the degree of a typical face in a
$\mathbf{P}_q$-distributed random map also has a heavy tail distribution.

\subsection{Choosing the Boltzmann
weights}\label{sec:choos-boltzm-weights}

We start from a sequence $q^\circ:=(q_k^\circ)_{k\in\N}$ of
nonnegative real numbers such that
%
\begin{equation}
\label{eq:13}
q_k^\circ\mathop{\sim}_{k\to\infty} k^{-a}
\end{equation}
for some real number $a>3/2$. In agreement with (\ref{eq:18}), we set
\[
f_\circ(x)=f_{q^\circ}(x)=\sum_{k=1}^\infty N(k)q^\circ_k x^{k-1}
\]
for every $x\geq0$. By Stirling's formula, we have
\[
N(k)\mathop{\sim}_{k\to\infty}
\frac{2^{2k-1}}{\sqrt{\pi k}}
\]
so that the radius of convergence of the series defining $f_\circ$ is $1/4$.
Furthermore, the condition $a>3/2$ guarantees that $f_\circ(1/4)$ and
$f'_\circ(1/4)$ are (well defined and) finite.
\begin{prp}
\label{sec:crit-boltzm-laws}
Set
\[
c=\frac{4}{4f_\circ(1/4) + f'_\circ(1/4)},\qquad
\beta=\frac{f'_\circ(1/4)}{4f_\circ(1/4) + f'_\circ(1/4)}
\]
and define a sequence $q=(q_k)_{k\in\N}$ by setting
%
\begin{equation}
\label{formula-q}
q_k= c(\beta/4)^{k-1} q^\circ_k.
\end{equation}
Then, the sequence $q$ is both admissible and critical, and
$Z_q=R_q=\beta^{-1}$.
\end{prp}
\begin{Remark*} As the proof will show, the choice given for the
constants $c$ and $\beta$
is the only one for which the conclusion of the proposition holds.
\end{Remark*}
\begin{pf*}{Proof of Proposition \ref{sec:crit-boltzm-laws}}
Consider a sequence $q=(q_k)_{k\in\N}$ defined as in the proposition,
with an arbitrary choice of the positive constants $c$ and $\beta$.
If $f_q$ is defined as in (\ref{eq:18}), it is immediate that
\[
f_q(x)=c f_\circ(\beta x/4).
\]
Hence, $R_q= \beta^{-1}$. Assume, for the moment, that the sequence
$q$ is admissible and $Z_q=R_q$. By Proposition \ref{sec:crit-boltzm-laws-2},
we have
$f_q(\beta^{-1})=1-\beta$ or, equivalently,
%
\begin{equation}
\label{choicebeta}
cf_\circ(1/4)=1-\beta.
\end{equation}
Furthermore, the criticality of $q$ will hold if and only if
$f'_q(\beta^{-1})=\beta^2$
or, equivalently,
%
\begin{equation}
\label{choice-c}
cf'_\circ(1/4)=4\beta.
\end{equation}
Conversely, if (\ref{choicebeta}) and (\ref{choice-c}) both hold,
then the sequence $q$ is admissible
by Proposition \ref{sec:crit-boltzm-laws-2}, the curves
$x\to f_q(x)$ and $x\to1-1/x$ are tangent at $x=\beta^{-1}$ and a simple
convexity
argument shows that $\beta^{-1}$ is the unique solution
of (\ref{eq:16}) so that $Z_q=\beta^{-1}=R_q$, again by Proposition
\ref{sec:crit-boltzm-laws-2}.

We conclude that the conditions (\ref{choicebeta}) and (\ref{choice-c})
are necessary and sufficient for the conclusion of the proposition to hold.
The desired result thus follows.
\end{pf*}

We now introduce our basic assumption, placing a further
restriction on the value of the parameter $a$.
\renewcommand{\theAssumption}{(A)}
\begin{Assumption}\label{assumpA} The sequence $q$ is
of the form given in Proposition \ref{sec:crit-boltzm-laws}, with
a sequence $q^\circ$ satisfying (\ref{eq:13})
for some $a\in(3/2,5/2)$. We set $\alpha:=a-1/2\in(1,2)$.
\end{Assumption}

This assumption will be in force throughout the remainder of this work,
with the
exception of the beginning of Section \ref{sec:boltzm-maps-branch}
(including Proposition \ref{mobile-map}), where we consider a general
admissible sequence $q$.

Many of the subsequent asymptotic results will be written in terms of the
constant $\beta$, which lies in the interval $(0,1)$, and the constant $c_0>0$
defined by
%
\begin{equation}
\label{formula-c-0}
c_0= \biggl(\frac{2c\Gamma(2-\alpha)}{\alpha(\alpha-1)
\beta\sqrt{\pi}} \biggr)^{1/\alpha} .
\end{equation}
The reason for introducing this other constant will become clearer
in Section \ref{sec:boltzm-maps-branch}.

\section{Coding maps with mobiles}\label{sec:coding-maps-with}

\subsection{The Bouttier--Di Francesco--Guitter
bijection}\label{sec:mobiles-bouttier-di}

Following \cite{BdFGmobiles}, we now recall how bipartite planar maps
can be coded by certain
labeled trees called \textit{mobiles}.

By definition, a plane tree $\mathcal{T}$ is a finite subset of the set
%
\begin{equation}
\label{namevertices}
\mathcal{U}=\bigcup_{n\geq0}^\infty\N^n
\end{equation}
of all finite sequences of
positive integers (including the empty sequence $\varnothing$) which
satisfies three obvious conditions. First, $\varnothing\in\mathcal{T}$.
Then, for
every $v=(u_1,\ldots,u_k)\in\mathcal{T}$ with $k\geq1$, the sequence
$(u_1,\ldots,u_{k-1})$ (the ``parent'' of $v$) also belongs to $\mathcal{T}$.
Finally, for every $v=(u_1,\ldots,u_k)\in\mathcal{T}$, there exists an
integer $k_v(\mathcal{T})\geq0$ (the ``number of children'' of $v$)
such that
$vj:=(u_1,\ldots,u_k,j)$ belongs to $\mathcal{T}$ if and only if $1\leq
j\leq
k_v(\mathcal{T})$. The elements of $\mathcal{T}$ are called
\textit{vertices}. The
generation of a vertex $v=(u_1,\ldots,u_k)$ is denoted by $|v|=k$. The notions
of an ancestor and a descendant in the tree $\mathcal{T}$ are defined
in an
obvious way.

For our purposes, vertices $v$ such that $|v|$ is
even will be called \textit{white} vertices and vertices $v$ such that
$|v|$ is
odd will be called \textit{black} vertices. We denote by
$\mathcal{T}^\circ$ (resp., $\mathcal{T}^\bullet$) the set of all white (resp.,
black) vertices of
$\mathcal{T}$.

A (rooted) mobile is a pair $\theta=(\mathcal{T},(\ell(v))_{v\in\mathcal
{T}^\circ})$
that consists of a plane tree and a collection of integer labels
assigned to
the white vertices of $\mathcal{T}$ such that the following properties hold:

\begin{enumerate}[(a)]
\item[(a)]
$\ell(\varnothing)=0$.
\item[(b)] Let $v\in\mathcal{T}^\bullet$, $v_{(0)}$ be the parent of $v$,
$p=k_v(\mathcal{T})+1$ and $v_{(j)}=vj$, $1\leq j\leq p-1$ be the
children of
$v$. Then, for every $j\in\{1,\ldots,p\}$, $\ell(v_{(j)})\geq
\ell(v_{(j-1)})-1$, where, by convention, $v_{(p)}=v_{(0)}$.
\end{enumerate}

\begin{figure}

\includegraphics{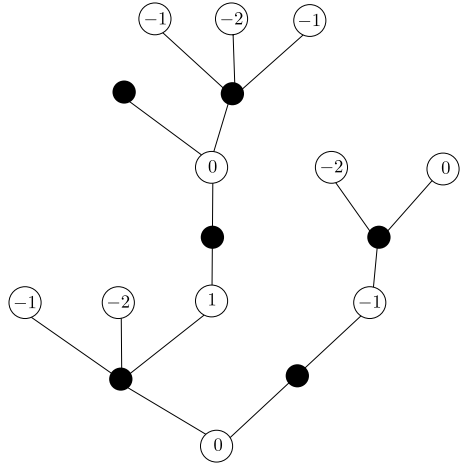}

\caption{A rooted mobile.}
\label{fig:mobile}
\end{figure}

Condition (b) means that if one lists the white vertices
adjacent to a given black
vertex in clockwise order, then the labels of these vertices can
decrease by at most~1
at each step. See Figure \ref{fig:mobile} for an example of a mobile.

%

We denote by $\Theta$ the (countable) set of all mobiles.
We will now describe the Bouttier--Di Francesco--Guitter (BDG)
bijection between $\Theta$
and
$\m_+^*$. This bijection can be found in Section 2 of
\cite{BdFGmobiles}, with the minor difference that \cite{BdFGmobiles}
deals with maps that are pointed, but not rooted.

Let $\theta=(\mathcal{T},(\ell(v))_{v\in\mathcal{T}^\circ})$ be a
mobile with $n+1$
vertices.
The contour sequence of $\theta$ is the
sequence $v_0,\ldots,v_{2n}$ of vertices of $\mathcal{T}$ which is
obtained by
induction as follows. First, $v_0=\varnothing$ and then, for every
$i\in\{0,\ldots, 2n-1\}$, $v_{i+1}$ is either the first child of
$v_i$ that has not yet appeared in the sequence $v_0,\ldots,v_i$ or
the parent of $v_i$ if all children of $v_i$ already appear in the
sequence $v_0,\ldots,v_i$. It is easy to verify that
$v_{2n}=\varnothing$ and that all vertices of $\mathcal{T}$ appear in the
sequence $v_0,v_1,\ldots,v_{2n}$. In fact, a given vertex $v$ appears
exactly $k_v(\mathcal{T})+1$ times in the contour sequence and each
appearance of
$v$ corresponds to one ``corner'' associated with this vertex.

The vertex $v_i$ is white when $i$ is even and black when $i$ is
odd. The contour sequence of $\mathcal{T}^\circ$, also called the white
contour sequence of $\theta$, is, by definition, the sequence
$v^\circ_0,\ldots,v^\circ_n$ defined by $v^\circ_i=v_{2i}$ for
every $i\in
\{0,1,\ldots,n\}$.

The image of $\theta$ under the BDG bijection is the element $(\mm
,v_*)$ of
$\m_+^*$ that is defined as follows. First, if $n=0$, meaning that
$\mathcal{T}=\{\varnothing\}$, we set \mbox{$(\mm,v_*)=\dagger$}. Suppose that
$n\geq1$ so that
$\mathcal{T}^\bullet$ has at least one element. We extend
the white
contour sequence of $\theta$ to a sequence $v^\circ_i$, $i\geq
0$, by periodicity, in such a way that $v^\circ_{i+n}=v^\circ_i$ for
every $i\geq0$. Then, suppose that the
tree $\mathcal{T}$ is embedded in the plane and add an extra vertex
$v_*$ not
belonging to the embedding. We construct a rooted planar map $\mm$
whose vertex set is
equal to
\[
V(\mm)=\mathcal{T}^\circ\cup\{v_*\}
\]
and whose edges are obtained by the
following device. For $i\in\{0,1,\ldots,n-1\}$, we let
\[
\phi(i)=\inf\{j> i\dvtx\ell(v^\circ_j)=\ell(v^\circ_i)-1\}\in
\{i+1,i+2,\ldots\}\cup\{\infty\} .
\]
We also set $v^\circ_\infty=v_*$, by convention. Then, for every
$i\in
\{0,1,\ldots,n-1\}$, we draw
an edge between $v^\circ_i$ and $v^\circ_{\phi(i)}$.
More precisely, the index $i$ corresponds to one specific ``corner'' of
$v^\circ_i$ and the associated edge starts from this corner.
The construction can then be made in such a way that edges
do not cross (and do not cross the edges of the tree) so that one
indeed gets a planar map.
This planar map $\mm$ is
rooted at the edge linking
$v^\circ_0=\varnothing$ to $v^\circ_{\phi(0)}$, which is oriented from
$v^\circ_{\phi(0)}$ to~$\varnothing$. Furthermore, $\mm$ is pointed
at the vertex
$v_*$, in agreement with our previous notation.

%
\begin{figure}[b]

\includegraphics{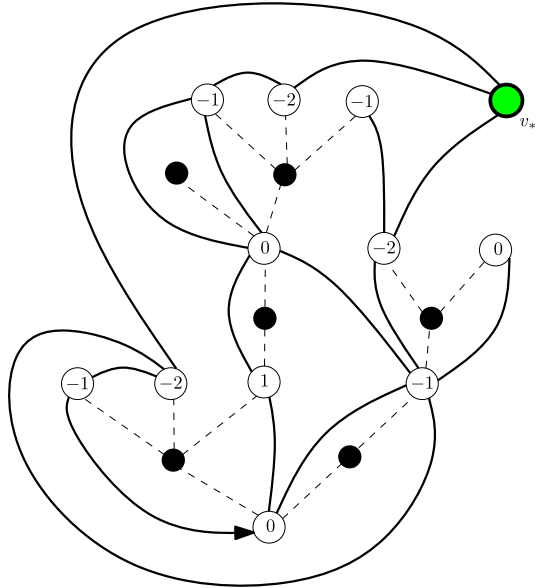}

\caption{The Bouttier--Di Francesco--Guitter construction for the mobile
of Figure \protect\ref{fig:mobile}.}
\label{fig:BDFG}
\end{figure}

See Figure \ref{fig:BDFG} for an example and Section 2 of \cite{BdFGmobiles} for a
more detailed
description.
\begin{prp}[(BDG bijection)]\label{sec:boutt-di-franc}
The preceding construction yields
a bijection from $\Theta$ onto $\m_+^*$. This bijection enjoys the
following two
properties:
\begin{enumerate}
\item
each face $f$ of $\mm$ contains exactly one vertex $v$ of
$\mathcal{T}^\bullet$, with $\deg(f)=2(k_v(\mathcal{T})+1)$;
\item
the graph distances in $\mm$ to the distinguished vertex $v_*$ are
linked to the labels of the mobile in the following way: for every
$v\in\mathcal{T}^\circ=V(\mm)\setminus\{v_*\}$,
\[
\dgr(v_*,v)=\ell(v)-\min_{v'\in\mathcal{T}^\circ}\ell(v')+1 .
\]
\end{enumerate}
\end{prp}

In our study of scaling limits of random planar maps, it will be important
to derive asymptotics for the random mobiles associated with these maps via
the BDG bijection. These asymptotics are
more conveniently stated in terms of random processes coding the mobiles.
Let us introduce such coding functions.

Let $\theta=(\mathcal{T},(\ell(v))_{v\in\mathcal{T}^\circ})$ be a
mobile with
$n+1$ vertices (so that $n=\#\mathcal{T}-1$) and let
$v^\circ_0,\ldots,v^\circ_n$ be, as previously, the white contour
sequence of
$\theta$. We set
%
\begin{equation}\label{eq:22}
C^{\theta}_i=\tfrac{1}{2}|v^\circ_i|\qquad
\mbox{for } 0\leq i\leq n ,\qquad
C^\theta_i=0 \qquad\mbox{for } i>n.
\end{equation}
We call $(C^{\theta}_i,0\leq i\leq n)$ the \textit{contour process} of
the mobile
$\theta$. It is a simple exercise to check that the contour process
$C^{\theta}$ determines
the tree $\mathcal{T}$. Similarly, we set
%
\begin{equation}\label{eq:22bis}
\Lambda^{\theta}_i=\ell(v^\circ_i)\qquad
\mbox{for } 0\leq i\leq n ,\qquad
\Lambda^\theta_i=0 \qquad\mbox{for } i>n
\end{equation}
and call $\Lambda^\theta$ the \textit{contour label process} of $\theta
$. The pair
$(C^\theta,\Lambda^\theta)$ determines the mobile $\theta$.

For technical reasons, we introduce variants of the preceding contour functions.
Let $n_\circ=\#\mathcal{T}^\circ-1$ and let
$w^\circ_0=\varnothing,w^\circ_1,\ldots,w^\circ_{n_\circ}$, be the
list of vertices of
$\mathcal{T}^\circ$ in lexicographical order. The
\textit{height process} of $\theta$ is defined by
\[
H^{\theta}_i=\tfrac{1}{2}|w^\circ_i| \qquad\mbox{for } 0\leq i\leq
n_\circ,\qquad
H^\theta_i=0 \qquad\mbox{for } i> n_\circ.
\]
Similarly, we introduce the \textit{label process}, which is defined by
\[
L^\theta_i=\ell(w^\circ_i) \qquad\mbox{for } 0\leq i\leq n_\circ,\qquad
L^\theta_i=0 \qquad\mbox{for } i > n_\circ.
\]

We will also\vspace*{1pt} need the Lukasiewicz path of $\mathcal{T}^\circ$. This is the
sequence $S^\theta=(S^\theta_0,S^\theta_1,\ldots)$, defined as
follows. First, $S^\theta_0=0$. Then, for every
$i\in\{0,1,\ldots,n_\circ\}$, $S^\theta_{i+1}-S^\theta_{i}+1$ is
the number of
(white) grandchildren
of $w^\circ_{i}$ in $\mathcal{T}$. Finally, $S^\theta_i=S^\theta
_{n_\circ+1}=-1$
for every $i>n_\circ$. It is easy to see that $S^\theta_i\geq0$ for
every $i\in\{0,1,\ldots,n_\circ\}$
so that
\[
\#\mathcal{T}^\circ=n_\circ+1=\inf\{i\geq0\dvtx S^\theta_i=-1\}.
\]

Let us briefly comment on the reason for introducing these
different processes.
In our applications to random planar maps, asymptotics for the
pair $(C^{\theta},\Lambda^\theta)$, which is directly linked to the
white contour sequence of $\theta$, turn out to be most useful. On the
other hand, in order to derive these asymptotics, it will be more convenient
to consider first the pair
$(H^{\theta},L^\theta)$.

In the following, the generic
element of $\Theta$ will be denoted by $(\theta,(\ell(v))_{v\in\mathcal{T}
^\circ}),$ as
previously.

\subsection{Boltzmann distributions and Galton--Watson trees}\label{sec:boltzm-maps-branch}

Let $q$ be an admissible sequence, in the sense of Section \ref{sec:discrete-picture}, and let
$M$ be a random element of
$\m_+^*$ with distribution $\mathbf{P}_q$. Our goal is to describe the
distribution of the
random mobile associated with $M$ via the
BDG bijection. We closely follow Section 2.2 of \cite{jfmgm05}.

We first need the notion of an alternating two-type Galton--Watson tree.
Recall that white vertices are those of even generation and
black vertices are those of odd generation.
Informally, an alternating two-type Galton--Watson
tree is just a Galton--Watson tree where white and black vertices
have a different offspring distribution. More precisely, if $\mu_0$
and $\mu_1$
are two probability distributions on the nonnegative integers,
the associated (alternating) two-type Galton--Watson tree is the random
plane tree whose distribution is specified
by saying that the numbers of children of the different vertices are
independent, the offspring distribution of each white vertex is $\mu_0$
and the offspring distribution of each black vertex is $\mu_1$; see
\cite{jfmgm05}, Section 2.2,
for a more rigorous presentation.

We also need to introduce the notion of a discrete bridge. Consider an
integer $p\geq1$
and the set
\[
E_{p}:= \Biggl\{(x_1,\ldots,x_{p})\in\{-1,0,1,2,\ldots\}^{p}\dvtx
\sum_{i=1}^{p} x_i=0 \Biggr\}.
\]
Note that $E_p$ is a finite set and, indeed, $\# E_p=N(p)$, with $N(p)$
as in (\ref{eq:18}). Let $(X_1,\ldots,X_p)$
be uniformly distributed over $E_p$. The sequence $(Y_0,Y_1,\ldots,Y_p)$
defined by $Y_0=0$ and
\[
Y_j=\sum_{i=1}^j X_i ,\qquad 1\leq j\leq p ,
\]
is called a \textit{discrete bridge of length} $p$.
\begin{prp}[(\cite{jfmgm05}, Proposition 7)]
\label{mobile-map}
Let $M$ be a random element of
$\m_+^*$ with distribution $\mathbf{P}_q$ and let $\theta=(\mathcal
{T},(\ell
(v),v\in
\mathcal{T}^\circ))$ be the
random mobile associated with $M$ via the
BDG bijection. Then:
\begin{enumerate}
\item the random tree $\mathcal{T}$ is an alternating two-type Galton--Watson
tree with
offspring distributions $\mu_0$ and $\mu_1$ given by
\[
\mu_0(k)=Z_q^{-1}f_q(Z_q)^k ,\qquad k\geq0 ,
\]
and
\[
\mu_1(k)=\frac{Z_q^kN(k+1)q_{k+1}}{f_q(Z_q)} ,\qquad k\geq0 ;
\]
\item conditionally given $\mathcal{T}$, the labels $(\ell(v),v\in
\mathcal{T}^\circ)$ are distributed uniformly over all possible choices
that satisfy
the constraints \textup{(a)} and \textup{(b)} in the definition of a mobile;
equivalently,
for every $v\in
\mathcal{T}^\bullet$,
with the notation introduced in property (b)
of the definition of a mobile, the sequence
$(\ell(v_{(j)})-\ell(v_{(0)}),0\leq j\leq k_v(\mathcal{T})+1)$ is
a discrete bridge of length $k_v(\mathcal{T})+1$
and these sequences are independent when $v$ varies over
$\mathcal{T}^\bullet$.
\end{enumerate}
\end{prp}

A random mobile having the distribution described in the
proposition will be called a $(\mu_0,\mu_1)$\textit{-mobile}.
The law $\Q$ of a $(\mu_0,\mu_1)$-mobile is a probability distribution
on $\Theta$.

Note that the respective means of
$\mu_0$ and $\mu_1$ are
\[
m_0:=\sum_{k\geq0}k\mu_0(k)=Z_qf_q(Z_q) ,\qquad m_1:=\sum_{k\geq
0}k\mu_1(k)=Z_qf_q'(Z_q)/f_q(Z_q)
\]
so that $m_0m_1=Z_q^2f'_q(Z_q)$ is
less than or equal to $1$ and equality holds if and only if $q$ is critical.

We now return to a weight sequence $q$ satisfying our basic Assumption
\ref{assumpA}.
Recall that the sequence $q$, which is both admissible and critical, is
given in terms of
the sequence $q^\circ$ by (\ref{formula-q}) and that we
have
$q^\circ_k\sim k^{-\alpha-1/2}$
as $k\to\infty$, with
$\alpha\in(1,2)$.

Then, $\mu_0$ is the geometric distribution with parameter
$f_q(Z_q)=1-\beta$ and
\[
\mu_1(k)=\frac{c}{1-\beta} 4^{-k} N(k+1) q^\circ_{k+1} ,\qquad
k=0,1,\ldots.
\]
From the asymptotic behavior of $q^\circ_k$, we obtain
\[
\mu_1(k)\mathop{\sim}_{k\to\infty}
\frac{2c}{(1-\beta)\sqrt{\pi}}
k^{-\alpha-1}.
\]
In particular, if we set
$\ov{\mu}_1(k)=\mu_1([k,\infty))$, this yields
%
\begin{equation}
\label{tail1}
\ov{\mu}_1(k) \mathop{\sim}_{k\to\infty}
\frac{2c}{\alpha(1-\beta)\sqrt{\pi}} k^{-\alpha}.
\end{equation}

Let $\mu$ be the probability distribution on the nonnegative integers
which is the law of
\[
\sum_{i=1}^U V_i,
\]
where $U$ is distributed according to $\mu_0$, $V_1,V_2,\ldots$
are distributed according to $\mu_1$ and the variables
$U,V_1,V_2,\ldots$ are independent. Then, $\mu$ is critical,
in the sense that
\[
\sum_{k=0}^\infty k \mu(k)=m_0m_1=1.
\]
Notice that $\mu$ is just
the distribution of the number of individuals at the second generation
of a $(\mu_0,\mu_1)$-mobile. It will be important to have information
on the tail
$\ov{\mu}(k):=\mu([k,\infty))$ of $\mu$. This follows easily from the
estimate (\ref{tail1}) and the definition of $\mu$. First, note that
\[
\ov{\mu}(k)
=\mathbb{P}\Biggl[\sum_{i=1}^U V_i\geq k \Biggr]
\geq\mathbb{P}[\exists i\in\{1,\ldots,U\}\dvtx V_i\geq k]
=1-\E\bigl[\bigl(1-\ov{\mu}_1(k)\bigr)^U \bigr].
\]
Then,
\[
1-\E\bigl[\bigl(1-\ov{\mu}_1(k)\bigr)^U \bigr]
=1-\frac{\beta}{1-(1-\ov{\mu}_1(k))(1-\beta)}
\mathop{\sim}_{k\to\infty} \frac{1-\beta}{\beta} \ov{\mu}_1(k).
\]
Using (\ref{tail1}), we get
\[
\ov{\mu}(k)\geq\frac{2c}{\alpha\beta\sqrt{\pi}} k^{-\alpha} +
o(k^{-\alpha}).
\]
A corresponding upper bound is easily obtained by
writing, for every $\vep\in(0,1/2)$,
\begin{eqnarray*}
\ov{\mu}(k)
&\leq& \mathbb{P}[\exists i\in\{1,\ldots,U\}\dvtx V_i\geq(1-\vep)k]
\\
&&{} +\mathbb{P}\Biggl[ \Biggl\{\sum_{i=1}^U V_i\geq k \Biggr\}\cap
\bigl\{\forall i\in\{1,\ldots,U\}\dvtx V_i\leq(1-\vep)k \bigr\} \Biggr]
\end{eqnarray*}
and checking that the second term in the right-hand side
is $o(k^{-\alpha})$ as $k\to\infty$. To see this, first note that
the probability
of the event $\{U>K\log k\}$ is $o(k^{-\alpha})$ if the constant $K$
is chosen sufficiently large. If $U\leq K\log k$, then the event in the
second term
may hold only if there are two distinct values of $i\in\{1,2,\ldots
,[K\log k]\}$
such that $V_i\geq\vep k/( K\log k)$. The desired estimate then follows
from (\ref{tail1}).

We have thus obtained
\[
\ov{\mu}(k)\mathop{\sim}_{k\to\infty}
\frac{2c}{\alpha\beta\sqrt{\pi}} k^{-\alpha},
\]
which we can rewrite in the form
%
\begin{equation}
\label{tail2}
\ov{\mu}(k)\mathop{\sim}_{k\to\infty}
\frac{\alpha-1}{\Gamma(2-\alpha)}c_0^\alpha k^{-\alpha}
\end{equation}
with the constant $c_0$ defined in (\ref{formula-c-0}).
The reason for introducing the constant $c_0$ and writing the
asymptotics (\ref{tail2})
in this form becomes clear when discussing scaling limits.
Recall that $1<\alpha<2$ by our assumption that
$\frac{3}{2}<a<\frac{5}{2}$. By (\ref{tail1}) or (\ref{tail2}),
$\mu$ is then in the domain of attraction of a stable law with index
$\alpha$. Recalling that $\mu$ is critical, we have the following, more
precise, result.

Let $\nu$ be the probability distribution on $\Z$
obtained by setting $\nu(k)=\mu(k+1)$ for every $k\geq-1$ [and
$\nu(k)=0$ if $k<-1$]. Let $S=(S_n)_{n\geq0}$ be a random walk on the
integers with jump distribution $\nu$. Then,
%
\begin{equation}
\label{Levyconv}
\bigl(n^{-1/\alpha}S_{[nt]} \bigr)_{t\geq0}
\mathop{\la}_{n\to\infty}^{\mathrm{(d)}} (c_0 X_t)_{t\geq0},
\end{equation}
where the convergence holds in distribution in the Skorokhod sense
and $X$ is a centered stable L\'evy process with index $\alpha$ and no
negative jumps, with Laplace transform given by
%
\begin{equation}
\label{Levytransform}
\E[\exp(-uX_t)]=\exp(tu^\alpha),\qquad t,u\geq0 .
\end{equation}
See, for instance, Chapter VII of Jacod and Shiryaev \cite{JS} for a thorough
discussion of the convergence of rescaled random walks toward
L\'evy processes.

\subsection{Discrete bridges}
\label{sec:bridges}

Recall from Proposition \ref{mobile-map} that the sequence of
labels of white vertices adjacent to a given black vertex in a
$(\mu_0,\mu_1)$-mobile is distributed as a discrete bridge. In this
section, we
collect some estimates for discrete bridges that will be used in
the proofs of our main results.

We consider
a random walk $(Y_n)_{n\geq0}$ on $\Z$ starting from $0$ and with jump
distribution
\[
\nu_*(k)= 2^{-k-2},\qquad k=-1,0,1,\ldots.
\]
Fix an integer $p\geq
1$ and let $(Y^{(p)}_n)_{0\leq n\leq p}$ be a vector whose
distribution is the conditional law of $(Y_n)_{0\leq n\leq p}$ given
that $Y_p=0$. Then,
the process $(Y^{(p)}_n)_{0\leq n\leq p}$ is a
discrete bridge with length $p$. Indeed, a simple calculation shows
that
\[
\bigl(Y^{(p)}_1,Y^{(p)}_2-Y^{(p)}_1,\ldots, Y^{(p)}_p-Y^{(p)}_{p-1}\bigr)
\]
is
uniformly distributed over the set
$E_{p}$.
\begin{lmm}
\label{estimatebridge}
For every real $r\geq1$, there exists a constant $K_{(r)}$ such that
for every
integer $p\geq1$ and $k,k'\in\{0,1,\ldots,p\}$,
\[
\E\bigl[\bigl(Y^{(p)}_k-Y^{(p)}_{k'}\bigr)^{2r}\bigr]\leq K_{(r)}|k-k'|^{r} .
\]
\end{lmm}
\begin{pf}
We may, and will, assume that $p\geq2$. Let us first suppose
that $k\leq k'\leq2p/3$. By the definition of $Y^{(p)}$ and then the
Markov property of $Y$, we have
\[
\E\bigl[\bigl(Y^{(p)}_k-Y^{(p)}_{k'}\bigr)^{2r}\bigr]=\frac{\E[|Y_k-Y_{k'}|^{2r}
\ind_{\{Y_p=0\}}]}{\mathbb{P}(Y_p=0)}= \E\biggl[|Y_k-Y_{k'}|^{2r}
\frac{\pi_{p-k'}(-Y_{k'})}{\pi_p(0)} \biggr] ,
\]
where
$\pi_n(x)=\mathbb{P}(Y_n=x)$ for every integer $n\geq0$ and $x\in\Z$. A
standard local limit theorem (see, e.g., Section 7 of \cite{spitzer})
shows that if $g(x)=(4\pi)^{-1/2}e^{-x^2/4}$, then we have
\[
\sqrt{n} \pi_{n}(x)=g\bigl(x/\sqrt{n}\bigr)+\eps_n(x) \qquad \mbox{where
}\sup_{x\in\Z}|\eps_n(x)|\mathop{\to}_{n\to\infty}0 .
\]
Then,
\[
\frac{\pi_{p-k'}(-Y_{k'})}{\pi_p(0)}\leq\sqrt{3} \frac{\sqrt
{p-k'}\pi_{p-k'}(-Y_{k'})}{\sqrt{p}\pi_p(0)}\leq K,
\]
where
\[
K=\sqrt{3} \frac{(4\pi)^{-1/2}+\sup_{n\geq
1}\sup_{x\in\Z}|\eps_n(x)|}{\inf_{n\geq1}\sqrt{n}\pi
_n(0)}<\infty.
\]
It follows that
\[
\E\bigl[\bigl(Y^{(p)}_k-Y^{(p)}_{k'}\bigr)^{2r}\bigr]\leq K\E[|Y_k-Y_{k'}|^{2r}] .
\]

Then,\vspace*{1pt} the bound $\E[|Y_k-Y_{k'}|^{2r}]\leq K'_{(r)}|k-k'|^{r},$
with a finite constant $K'_{(r)}$ depending only on $r$, is a
consequence of Rosenthal's inequality for i.i.d. centered random
variables \cite{petrov75}, Theorem 2.10. We have thus obtained the
desired estimate under the restriction $k\leq k'\leq2p/3$.

If $p/3\leq k\leq k'\leq
p$, the same estimate is readily obtained by observing
that $(-Y^{(p)}_{p-n},0\leq n\leq p)$ has the same distribution as
$Y^{(p)}$. Finally, in the case $k\leq p/3\leq2p/3\leq k'$, we
apply the preceding bounds successively to
$\E[|Y_k-Y_{[p/2]}|^{2r}]$ and to $\E[|Y_{[p/2]}-Y_{k'}|^{2r}]$.
\end{pf}

An immediate consequence of the lemma (applied with $r=1$) is the bound
%
\begin{equation}
\label{momentbridge2}
\E\bigl[\bigl(Y^{(p)}_j\bigr)^2\bigr]\leq K_{(1)} \min\{j,p-j\}\leq2K_{(1)} \frac{j(p-j)}{p}
\end{equation}
for every integer $p\geq2$ and $j\in\{0,1,\ldots,p\}$.

Finally, a conditional version of Donsker's theorem gives
%
\begin{equation}
\label{Donsker}
\biggl(\frac{1}{\sqrt{2p}}Y^{(p)}_{[pt]} \biggr)_{0\leq t\leq1}
\mathop{\la}_{p\to\infty}^{\mathrm{(d)}} (\gamma_t)_{0\leq t\leq1},
\end{equation}
where $\gamma$ is a standard Brownian bridge. Such results are part
of the folklore of the subject; see Lemma 10 in \cite{betti} for a detailed
proof of a more general statement.

\section{The continuous distance process}\label{sec:continuous-object}

Our goal in this section is to discuss the so-called continuous distance
process, which will appear as the scaling limit of the label processes
$L^\theta$ and $\Lambda^\theta$ of Section \ref
{sec:mobiles-bouttier-di} when $\theta$
is a $(\mu_0,\mu_1)$-mobile conditioned to be large in some sense.

\subsection{Definition and basic properties}\label{sec:defin-first-prop}

We consider the centered stable L\'evy process $X$
with no negative jumps and index $\alpha$, and Laplace exponent as in
(\ref{Levytransform}). The canonical filtration associated with $X$ is
defined, as usual, by
\[
\FF_t=\sigma\{X_s,0\leq s\leq
t\}
\]
for every $t\geq0$.
We let $(t_i)_{i\in\N}$ be a measurable enumeration of the jump times of
$X$ and set $x_i=\Delta X_{t_i}$ for every $i\in\N$. Then, the point measure
\[
\sum_{i\in\N} \delta_{(t_i,x_i)}
\]
is Poisson on $[0,\infty)\times[0,\infty)$ with intensity
\[
\frac{\alpha(\alpha-1)}{\Gamma(2-\alpha)} \,{d}t\, \frac{{d}x}{x^{\alpha+1}}.
\]

For $s\leq t$, we set
\[
I_t^s=\inf_{s\leq r\leq t}X_r
\]
and $I_t=I^0_t$. For every $x\geq0$, we set
\[
T_x=\inf\{t\geq0\dvtx -I_t>x\}.
\]
We recall
that the process $(T_x,x\geq0)$ is a stable subordinator
of index $1/\alpha$ with Laplace transform
%
\begin{equation}\label{eq:2}
\E[\exp(-u T_x)]=\exp(-x u^{1/\alpha}) ;
\end{equation}
see, for example, Theorem 1 in \cite{bertlev96}, Chapter VII.

Suppose that, on the same probability space, we are given
a sequence $(b_i)_{i\in\N}$ of independent (one-dimensional) standard Brownian
bridges over the time interval $[0,1]$ starting and ending at the
origin. Assume that
the sequence $(b_i)_{i\in\N}$ is independent of the L\'evy process $X$.
Then, for every $i\in\N$, we introduce the
rescaled bridge
\[
\wt{b}_i(r)=x_i^{1/2}b_i(r/x_i) ,\qquad 0\leq r\leq x_i ,
\]
which, conditionally on $\FF_\infty$, is a standard Brownian
bridge with duration $x_i$.

Recall that $X_{s-}$ denotes the left limit of
$X$ at $s$ for every $s>0$.
\begin{prp}
\label{defdistance}
For every $t\geq0$, the series
%
\begin{equation}\label{eq:1}
\sum_{i\in\N}\wt{b}_i(I_t^{t_i}-X_{t_i-})\ind_{\{X_{t_i-}\leq
I_t^{t_i}\}} \ind_{\{t_i\leq t\}}
\end{equation}
converges in $L^2$-norm. The sum of this series is
denoted by $D_t$. The process $(D_t,t\geq0)$
is called the continuous distance process.
\end{prp}
\begin{Remark*}
In a more compact form, we can write
\[
D_t=\sum_{i\in\N\dvtx t_i\leq t}\wt{b_i}\bigl((I_t^{t_i}-X_{t_i-})^+\bigr) .
\]
\end{Remark*}
\begin{pf*}{Proof of Proposition \ref{defdistance}}
Note that in (\ref{eq:1}), the summands are
well defined since, obviously, $I_t^{t_i}\leq X_{t_i}$ for every
$t_i\leq t$ so
that $I_t^{t_i}-X_{t_i-}\leq\Delta
X_{t_i}=x_i$. The nonzero summands in (\ref{eq:1}) correspond to those
values of $i$
for which $t_i\leq t$ and $X_{t_i-}\leq
I_t^{t_i}$. Conditionally on $\FF_\infty$, these summands are independent
centered Gaussian random variables
with respective variances
\[
\E[\wt{b}_{i}(I_t^{t_i}-X_{t_i-})^2 |\FF_\infty]
=\frac{(I_t^{t_i}-X_{t_i-})(X_{t_i}-I_t^{t_i})}{x_i}\leq
I_t^{t_i}-X_{t_i-} .
\]
The equality in the previous
display follows from the fact that $\var b_{(a)}(t)=\frac{t(a-t)}{a}$
whenever $b_{(a)}$ is a
Brownian bridge with duration $a>0$ and $0\leq t\leq a$.

We then have
\begin{eqnarray*}
&&\E\biggl[ \sum_{i\in\N}\wt{b}_i(I_t^{t_i}-X_{t_i-})^2\ind_{\{
X_{t_i-}\leq
I_t^{t_i}\}} \ind_{\{t_i\leq t\}} \biggr]\\
&&\qquad \leq\E\biggl[ \sum_{i\in\N}(I_t^{t_i}-X_{t_i-}) \ind_{\{X_{t_i-}\leq
I_t^{t_i}\}} \ind_{\{t_i\leq t\}} \biggr]\\
&&\qquad = \E\biggl[ \sum_{t_i\leq t}(I_t^{t_i}-I_t^{t_i-}) \biggr]
\leq\E[X_t-I_t]= \E[-I_t],
\end{eqnarray*}
where the last equality holds because $X$ is centered. It is well known
that \mbox{$\E[-I_t]<\infty$}.
Indeed, $-I_t$ even has exponential moments; see Corollary 2 in
\cite{bertlev96}, Chapter VII. Since
the summands in (\ref{eq:1}) are centered and orthogonal in $L^2$, the
desired convergence readily follows
from the preceding estimate.
\end{pf*}

In order to simplify the presentation, it will be convenient to adopt a
point process notation, by letting $(x_s,b_s)=(x_i,b_i)$ whenever
$t_i=s$ for some $i\in\N$ and, by convention, $x_s=0$, $b_s=0$
(i.e.,\vspace*{1pt} the path with duration zero started from the
origin) when $s\notin\{t_i,i\in\N\}$. The process $\wt{b}_s$ is defined
accordingly and is equal to $0$ when $b_s=0$. We can thus rewrite
%
\begin{equation}\label{eq:19}
D_t=\sum_{0<s\leq t}\wt{b}_s\bigl((I_t^s-X_{s-})^+\bigr) .
\end{equation}

Let us conclude this section with a useful scaling property.
For every $r>0$, we have
%
\begin{equation}\label{eq:20}
(r^{-1/\alpha}X_{rt},r^{-1/2\alpha}D_{r t} )_{t\geq
0}\stackrel{\mathrm{(d)}}{=}(X_t,D_t)_{t\geq0} .
\end{equation}
This easily follows from our construction and the scaling property of $X$.

\subsection{\texorpdfstring{H\"older regularity}{Holder regularity}}\label{sec:holder-regularity}

In this subsection, we prove the following
regularity property of $D$.
\begin{prp}\label{sec:continuous-object-4}
The process $(D_t,t\geq0)$ has a modification
that is locally H\"older continuous with any exponent $\eta\in
(0,1/2\alpha)$.
\end{prp}

We start with a few preliminary lemmas.
\begin{lmm}\label{sec:continuous-object-2}
For every real $t>0$ and $r>-1$, we have $\E[(-I_t)^r]<\infty$.
\end{lmm}
\begin{pf}
By scaling, it is enough to consider $t=1$. As mentioned in the last proof,
the case $r\geq0$ is a consequence of Corollary 2 in
\cite{bertlev96}, Chapter VII.
To handle the case $r<0$,
we use a scaling argument to write
\[
\mathbb{P}(-I_1> x)=\mathbb{P}(T_x < 1)
=\mathbb{P}(x^{\alpha}T_1 <1)=\mathbb{P}\bigl((T_1)^{-1/\alpha}>x\bigr),
\]
so $-I_1$ has the same distribution as
$T_1^{-1/\alpha}$. We have already observed that the
process $(T_x,x\geq0)$ is a stable subordinator with index
$1/\alpha$. This implies that
$\E[(T_1)^s]<\infty$ for every $0\leq s< 1/\alpha$,
from which the desired result follows.
\end{pf}
\begin{lmm}\label{sec:continuous-object-1}
For every real $t\geq0$ and $r>0$, we have $\E[|D_t|^r]<\infty$.
\end{lmm}
\begin{pf}
Again by scaling, we may concentrate on the case $t=1$. Arguing as in
the proof of Proposition
\ref{defdistance}, we get that, conditionally on $\FF_\infty$, the
random variable $D_1$
is a centered Gaussian variable with variance
\[
\sum_{0<s\leq1}\frac{(I_1^s-X_{s-})(X_s-I_1^s)}{\Delta
X_s}\ind_{\{X_{s-}<I_1^s\}}\leq\sum_{0<s\leq
1}(X_s-I_1^s)\ind_{\{X_{s-}<I_1^s\}} .
\]
Note that this time, we
chose the upper bound $X_s-I_1^s$ rather than $I_1^s-X_{s-}$ for the
summands as the latter is ineffective for getting finiteness of high
moments. Thus, if ${\mathcal N}$ denotes a standard normal variable
and $K_r= \E[|{\mathcal N}|^r]$, we have
%
\begin{eqnarray}\label{eq:6}\quad
\E[|D_1|^r] &=& \E[|{\mathcal N}|^r]\times\E\biggl[ \biggl(\sum_{0<s\leq1}\frac
{(I_1^s-X_{s-})(X_s-I_1^s)}{\Delta
X_s}\ind_{\{X_{s-}<I_1^s\}} \biggr)^{r/2} \biggr]\nonumber\\[-8pt]\\[-8pt]
&\leq& K_r \E\biggl[ \biggl(\sum_{0<s\leq
1}(X_s-I_1^s)\ind_{\{X_{s-}<I_1^s\}} \biggr)^{r/2} \biggr] .\nonumber
\end{eqnarray}
By a standard time-reversal
property of L\'evy processes, the process $(X_1-X_{(1-s)-},0\leq s<
1)$ has the same distribution as $(X_s,0\leq s< 1)$, which entails that
%
\begin{equation}\label{eq:7}
\sum_{0<s\leq
1}(X_s-I_1^s)\ind_{\{X_{s-}<I_1^s\}}\stackrel{\mathrm{(d)}}{=}\sum_{0<s\leq
1}(\ov X_{s-}-X_{s-})\ind_{\{\ov X_{s-}<X_s\}} ,
\end{equation}
where $\ov X_s=\sup_{0\leq r\leq s}X_r$. For every integer $k\geq0$,
we introduce the
process
\[
A^{(k)}_t=\sum_{0<s\leq
t}(\ov X_{s-}-X_{s-})^{2^k}\ind_{\{\ov X_{s-}<X_s\}} ,\qquad t\geq
0 ,
\]
which is an increasing c\`adl\`ag process adapted to the
filtration $(\FF_t)$, with compensator
\begin{eqnarray*}
\wt{A}^{(k)}_t&=&\frac{\alpha(\alpha-1)}{\Gamma(2-\alpha)} \int
_0^t{d}s(\ov X_s-X_s)^{2^k}\int_0^\infty\frac{{d}
x}{x^{\alpha+1}}\ind_{\{\ov X_s<X_s+x\}}\\
&=&\frac{\alpha-1}{
\Gamma(2-\alpha)}\int_0^t(\ov X_s-X_s)^{2^k-\alpha} \,{d}s .
\end{eqnarray*}
Note that
$\E[\wt{A}^{(k)}_t]<\infty$ since this expectation is
\[
\frac{\alpha-1}{\Gamma(2-\alpha)}\E[(\ov X_1-X_1)^{2^{k}-\alpha}]
\int_0^ts^{2^{k}/\alpha-1}\,{d}s
\]
and time reversal shows that $\E[(\ov X_1-X_1)^{2^{k}-\alpha
}]=E[(-I_1)^{2^{k}-\alpha}]<\infty$,
by Lem\-ma~\ref{sec:continuous-object-2}, since
$2^{k}-\alpha\geq1-\alpha>-1$. In order to complete the
proof of Lemma~\ref{sec:continuous-object-1}, we will need the
following, stronger, fact.
\begin{lmm}\label{sec:continuous-object-3}
For all integers $k,p\geq0$, we have
$\E[(\wt{A}^{(k)}_1)^p]<\infty$.
\end{lmm}
\begin{pf}
We must show that
%
\begin{equation}\label{eq:5}
\int_{[0,1]^p}{d}s_1\cdots{d}s_p
\E\Biggl[\prod_{i=1}^p(\ov X_{s_i}-X_{s_i})^{2^k-\alpha} \Biggr]<\infty.
\end{equation}
When $k\geq1$, we have $2^k-\alpha>0$ and the result easily follows
from H\"older's inequality, using a scaling argument, then time reversal
and Lemma \ref{sec:continuous-object-2}, just as we did to verify that
$\E[\wt{A}^{(k)}_t]<\infty$. The case $k=0$ is slightly more
delicate. We rewrite the left-hand side of (\ref{eq:5}) as
\[
p!\int_{0\leq s_1\leq\cdots\leq s_p\leq1}{d}s_1\cdots{d}s_p\,
\E\Biggl[\prod_{i=1}^p(\ov X_{s_i}-X_{s_i})^{1-\alpha} \Biggr] .
\]
By Proposition 1 in
\cite{bertlev96}, Chapter VI, the reflected process $\ov X-X$ is
Markov with respect to the filtration $(\FF_t)$. When started from a value
$x\geq0$, this Markov process has the same distribution as $x\vee\ov
X-X$ under
$\mathbb{P}$ and thus stochastically dominates $\ov X-X$ (started from $0$).
Consequently, since $1-\alpha<0$, we get,
for $0=s_0\leq s_1\leq\cdots\leq s_p\leq1$,
that
\begin{eqnarray*}
&&\E\Biggl[\prod_{i=1}^p(\ov X_{s_i}-X_{s_i})^{1-\alpha} \Biggr]\\
&&\qquad =
\E\Biggl[(\ov X_{s_1}-X_{s_1})^{1-\alpha}
\E\Biggl[\prod_{i=2}^p(\ov X_{s_i}-X_{s_i})^{1-\alpha} \Big|
\ov X_{s_1}-X_{s_1} \Biggr] \Biggr]\\
&&\qquad \leq
\E\Biggl[(\ov X_{s_1}-X_{s_1})^{1-\alpha}
\E\Biggl[\prod_{i=2}^p(\ov X_{s_i-s_1}-X_{s_i-s_1})^{1-\alpha}
\Biggr] \Biggr]\\
&&\qquad \leq
\prod_{i=1}^{p}\E[(\ov X_{s_i-s_{i-1}}-X_{s_i-s_{i-1}})^{1-\alpha}]
\end{eqnarray*}
by induction. Finally, by scaling and a simple change of variables, we
get that
(\ref{eq:5}) is bounded above by
\[
p! \E[(\ov X_1-X_1)^{1-\alpha} ]^p\int_{[0,1]^p}\prod_{i=1}^p
s_i^{1/\alpha-1}\,{d}s_i ,
\]
which is finite by Lemma \ref{sec:continuous-object-2}
since $\ov X_1-X_1\stackrel{\mathrm{(d)}}{=} -I_1$, by time reversal.
\end{pf}

We now complete the proof of Lemma \ref{sec:continuous-object-1}. Note
that $A^{(k+1)}$ is the square bracket of the compensated martingale
$A^{(k)}-\wt{A}^{(k)}$ for every $k\geq0$. For any real $r\geq1$, the
Burkholder--Davis--Gundy inequality \cite{DeMe5-8}, Chapter VII.92,
gives the existence of a finite constant
$K'_r$, depending only on $r$, such that
\[
\E\bigl[\bigl|A_1^{(k)}-\wt{A}_1^{(k)}\bigr|^r \bigr]\leq K'_r
\E\bigl[ \bigl(A^{(k+1)}_1 \bigr)^{r/2} \bigr] .
\]
Since $\wt{A}^{(k)}_1$
has moments of arbitrarily high order by Lemma
\ref{sec:continuous-object-3}, and $\E[A^{(k)}_1]=\E[\wt
{A}^{(k)}_1]<\infty$, a
repeated use of the last inequality shows that $\E
[(A^{(k-i)}_1)^{2^i}]<\infty$
for every $i=0,\ldots,k$. In particular,
$\E[(A^{(0)}_1)^{2^k}]<\infty$ for every integer $k\geq0$. The
desired result
now follows from
(\ref{eq:6}) and (\ref{eq:7}).
\end{pf}
\begin{pf*}{Proof of Proposition \ref{sec:continuous-object-4}}
Fix $s\geq0$ and $t>0$. Let $u=\sup\{r\in(0,s]\dvtx\break X_{r-}<I_{s+t}^s\}$
with the convention that
$\sup\varnothing=0$. Then, $I^r_{s+t}=I^r_s$ for every $r\in[0,u)$,
whereas $I^r_{s+t}=I^s_{s+t}$ for $r\in[u,s]$. By splitting the sum
(\ref{eq:19}), we get
\[
D_s=\sum_{0<r<u}\wt{b}_r\bigl((I_s^r-X_{r-})^+\bigr)+\wt{b}_u\bigl((I_s^u-X_{u-})^+\bigr)
+\sum_{u<r\leq s}\wt{b}_r\bigl((I_s^r-X_{r-})^+\bigr)
\]
and, similarly,
\begin{eqnarray*}
D_{s+t}&=&\sum_{0<r<u}\wt{b}_r\bigl((I_{s+t}^r-X_{r-})^+\bigr)+\wt
{b}_u\bigl((I_{s+t}^u-X_{u-})^+\bigr)\\
&&{}+
\sum_{s<r\leq s+t}\wt{b}_r\bigl((I_{s+t}^r-X_{r-})^+\bigr) .
\end{eqnarray*}
In the last display, we should also have considered the sum over
$r\in(u,s]$, but, in fact, this term gives no contribution because we
have $X_{r-}\geq
I_{s+t}^s=I_{s+t}^r$ for these values of $r$, by the definition of $u$.
Moreover, as
$I^r_s=I^r_{s+t}$ for $r\in[0,u)$, we have
\[
\sum_{0<r<u}\wt{b}_r\bigl((I_s^r-X_{r-})^+\bigr)=
\sum_{0<r<u}\wt{b}_r\bigl((I_{s+t}^r-X_{r-})^+\bigr) .
\]
Also, a simple translation argument shows that we may write
\[
\sum_{s<r\leq s+t}\wt{b}_r\bigl((I_{s+t}^r-X_{r-})^+\bigr) = D^{(s)}_t,
\]
where the process $D^{(s)}$ has the same distribution as $D$ and, in
particular, $D^{(s)}_t$ has the same
distribution as $D_t$. By combining the preceding remarks, we get
\begin{eqnarray*}
D_{s+t}-D_s-D^{(s)}_t &=& -\sum_{u<r\leq s}\wt{b}_r\bigl((I_s^r-X_{r-})^+\bigr)\\
&&{} + \bigl(\wt{b}_u\bigl((I_{s+t}^u-X_{u-})^+\bigr)-\wt{b}_u\bigl((I_s^u-X_{u-})^+\bigr) \bigr) .
\end{eqnarray*}
Conditionally on $\FF_{\infty}$, the right-hand side of the last
display is
distributed as a centered
Gaussian variable with variance bounded above by
\begin{eqnarray*}
\sum_{u<r\leq s}(I^r_s-X_{r-})^++(I_s^u-I_{s+t}^u)&=&
\sum_{u<r\leq s}(I^r_s-I^{r-}_s)+(I_s^u-I_{s+t}^u)\\
&\leq& X_s-I_{s+t}^u= X_s-I_{s+t}^s .
\end{eqnarray*}
Furthermore, $X_s-I_{s+t}^s$ has the same distribution as
$-I_t$, by the Markov property of $X$.

Now, let $p\geq1$. From previous considerations, we obtain
\begin{eqnarray*}
\E[|D_{s+t}-D_s|^p]&\leq&
2^{p} \bigl(\E\bigl[\bigl|D^{(s)}_t\bigr|^p\bigr]+\E\bigl[\bigl|D_{s+t}-D_s-D^{(s)}_t\bigr|^p\bigr] \bigr)\\
&\leq& 2^{p} \bigl(\E[|D_t|^p]+K_p\E[(-I_t)^{p/2}]\bigr)\\
&=& 2^{p} \bigl(\E[|D_1|^p]+K_p\E[(-I_1)^{p/2}] \bigr)
t^{p/2\alpha} ,
\end{eqnarray*}
where we have made
further use of the scaling properties of $X$ and $D$. The constant in front
of $t^{p/2\alpha}$ is finite, by Lemmas \ref{sec:continuous-object-2}
and \ref{sec:continuous-object-1}. The classical Kolmogorov
continuity criterion then yields the desired result.
\end{pf*}

In what follows, we will always consider the continuous modification of
$(D_t$, $t\geq0)$.
\begin{Remark*}
The process $D$ is closely related to the so-called
exploration process associated with $X$, as defined in the monograph
\cite{duqleg02}. The latter is a measure-valued strong Markov process
$(\rho_t,t\geq0)$ such that, for every $t\geq0$,
$\rho_t$ is an atomic measure on $[0,\infty)$ and the masses of the
atoms of $\rho_t$ are precisely the
quantities $(I_t^s-X_{s-})^+$, $s\leq t$, that are involved in the
definition of $D_t$ (see the proof of Theorem \ref{GHconv} below for more
information on this exploration process). As a matter of fact, part of
the proof of
Proposition \ref{sec:continuous-object-4} resembles the proof of the
Markov property for $(\rho_t,t\geq0)$; see \cite{duqleg02}, Proposition
1.2.3. However, the definition of $\rho_t$ requires the
introduction of the continuous-time height process (see the next
section), which is not needed in the definition of
$D_t$.
\end{Remark*}

\subsection{Excursion measures}\label{sec43}

It will be useful to consider the distance process $D$
under the excursion measure of $X$ above its minimum process $I$.
Recall that $X-I$ is a strong Markov process, that $0$ is a regular recurrent
point for this Markov process and that $-I$ provides a local time for $X-I$
at level $0$ (see \cite{bertlev96}, Chapters VI and VII).
We write $\bN$ for the excursion measure of $X-I$ away from
$0$ associated with this choice of local time. This excursion measure
is defined on the Skorokhod space
$\D(\R)$ and,
without risk of confusion, we will also use the notation $X$ for the
canonical process
on the space $\D(\R)$.
The duration of the excursion under $\mathbf{N}$ is
$\sigma=\inf\{t>0\dvtx X_t=0\}$. For every $a>0$, we have
\[
\mathbf{N}(\sigma\in{d}a)=\frac{{d}
a}{\alpha\Gamma(1-1/\alpha)a^{1+1/\alpha}} .
\]
This easily follows from formula (\ref{eq:2}) for the
Laplace transform of $T_x$.

In order to assign an independent bridge to each jump of $X$,
we consider an auxiliary probability space $(\Omega^*,\FF^*,\mathbb
{P}^*)$ which
supports a countable collection of independent Brownian bridges
$(b_i)_{i\in\N}$.
We then argue on the product space $\D(\R)\times\Omega^*$, which is equipped
with the product measure $\bN\otimes\mathbb{P}^*$. With a slight abuse
of notation,
we will write $\bN$ instead of $\bN\otimes\mathbb{P}^*$ in what follows.

The construction of the distance process under $\bN$ is then similar
to the constructions in the preceding subsections.
The process $X$ has a countable number of jumps under $\mathbf{N}$
and these jumps can be enumerated, for instance, by
decreasing size, as a sequence $(t_i)_{i\in\N}$. The same formula
(\ref{eq:1}) can be used to define the distance process $D_t$
under $\bN$. It is again easy to check that the series (\ref{eq:1})
converges, say in $\bN$-measure. Note that $D_t=0$
on $\{\sigma\leq t\}$.

To connect this construction with the previous subsections, we
may consider, under the probability measure $\mathbb{P}$, the first
excursion interval
of $X-I$ (away from~$0$) with length greater than $a$, where $a>0$ is fixed.
We denote this interval by $(g_{(a)},d_{(a)})$.
Then, the distribution of $(X_{(g_{(a)}+t)\wedge d_{(a)}},t\geq0)$
under $\mathbb{P}$
coincides with that of $(X_t,t\geq0)$ under $\bN(\cdot\mid\sigma
>a)$. Furthermore,
it is easily checked that the finite-dimensional marginals of the process
$(D_{(g_{(a)}+t)\wedge d_{(a)}},t\geq0)$ under $\mathbb{P}$ also coincide
with those of $(D_t,t\geq0)$ under $\bN(\cdot\mid\sigma>a)$. The
point here is that the
only jumps that may give a nonzero contribution in formula (\ref{eq:1})
are those that belong to the excursion interval of $X-I$
that straddles $t$. From the previous observations and Proposition
\ref{sec:continuous-object-4}, we deduce that the process $(D_t,t\geq0)$
also has a H\"older continuous modification under $\bN$ and, from
now on, we will deal with this modification.

Finally, it is well known that the scaling properties of stable processes
allow one to make sense of the conditioned measure $\bN(\cdot\mid
\sigma=a)$ for any choice of $a>0$.
Using the scaling property (\ref{eq:20}), it is then a simple matter
to define
the distance process $D$ also under this conditioned measure.
Furthermore, the H\"older continuity properties of $D$ still hold under
$\bN(\cdot\mid\sigma=a)$.

\section{Convergence of labels in a forest of mobiles}
\label{forest-mobiles}

We now consider a sequence $\mathbf{F}=(\theta_1,\theta_2,\ldots)$ of
independent random mobiles.
We assume that, for every $i\in\N$,
$\theta_i=(\mathcal{T}_i,(\ell_i(v),v\in\tc_i))$
is a $(\mu_0,\mu_1)$-mobile. We will call $\mathbf{F}$
a (\textit{random}) \textit{labeled forest}. It will also be useful to consider
the (unlabeled) forest
$\F$, defined as the sequence $(\mathcal{T}_1,\mathcal{T}_2,\ldots)$.

For our purposes, it will be important to
distinguish the vertices of the different trees in the forest $\F$.
This can be achieved
by a minor modification of the formalism of Section \ref
{sec:mobiles-bouttier-di},
letting $\mathcal{T}_1$ be a (random) subset of $\{1\}\times\mathcal
{U}$, $\mathcal{T}
_2$ be a subset of
$\{2\}\times\mathcal{U}$ and so on. Whenever we deal with a sequence
of trees or of mobiles,
we will tacitly assume that this modification has been made.

Our goal is to study the
scaling limit of the collection of labels in the forest $\mathbf{F}$.

\subsection{Statement of the result}\label{sec:conv-height-proc}

We first recall known results about scaling limits of the height process.
We let $(H^\circ_n)_{n\geq0}$ denote the height process of the forest
$\mathbf F$. This means that
the process $H^\circ$ is obtained by concatenating the height
processes $[H^{\theta_i}(n),0\leq n\leq\#\tc_i-1]$ of the mobiles
$\theta_i$.
Equivalently, let $u_0,u_1,\ldots$ be the sequence of all white
vertices of the
forest $\F$, listed one tree after another and in lexicographical
order for each tree. Then,
$H^\circ_n$ is equal to half the generation of $u_n$.

Scaling limits of $(H^\circ_n)_{n\geq0}$ are better understood, thanks
to the connection between the height process and the Lukasiewicz path
of the forest $\mathbf F$. We denote this Lukasiewicz path by
$(S^\circ_n)_{n\geq0}$. This means that $S^\circ_0=0$ and, for every
integer $n\geq0$, $S^\circ_{n+1}-S^\circ_{n}+1$ is the number of
(white) grandchildren of $u_n$ in $\F$. Then, $(S^\circ_n)_{n\geq0}$
is a random walk with jump distribution $\nu$, as defined before
(\ref{Levyconv}). To see this, note that for every $i\in\N$, the set
$\tc_i$ of all white vertices of $\mathcal{T}_i$ can be viewed as a
plane tree,
simply by saying that a white vertex of $\mathcal{T}_i$ is a child in
$\tc_i$ of
another white vertex of $\mathcal{T}_i$ if and only if it is a
grandchild of
this other vertex in the tree $\mathcal{T}_i$. Modulo this identification,
$\tc_1,\tc_2,\ldots$ are independent Galton--Watson trees with
offspring distribution $\mu$. The fact that $(S^\circ_n)_{n\geq0}$ is
a random walk with jump distribution $\nu$ is then a consequence of
well-known results for forests of i.i.d. Galton--Watson trees; see, for
example, Section 1 of \cite{Trees}.

Moreover,
the height process $(H^\circ_n)_{n\geq0}$
is related to the random walk $(S^\circ_n)_{n\geq0}$ by the formula
%
\begin{equation}
\label{heightLuka}
H^\circ_n=\#\Bigl\{k\in\{0,1,\ldots,n-1\}\dvtx S^\circ_k=\min_{k\leq j\leq
n} S^\circ_j\Bigr\}.
\end{equation}
The integers $k$ that appear in the right-hand side of
(\ref{heightLuka}) are exactly those for which $u_k$ is an ancestor of $u_n$
distinct from $u_n$. For each such integer $k$, the quantity
%
\begin{equation}
\label{rankchild}
S^\circ_{k+1}-\min_{k+1\leq j\leq n} S^\circ_j+1
\end{equation}
is the rank of $u_{k+1}$ among the grandchildren of $u_k$ in $\F$.
We again refer to Section 1 of \cite{Trees} for a thorough discussion of
these results and related ones.
For every integer $k$ such that
$u_k$ is a strict ancestor of $u_n$, it will also be of interest to
consider the (black) parent of $u_{k+1}$ in the forest $\F$.
As a consequence of the preceding remarks,
the number of children of this black vertex is less than or equal to
$S^\circ_{k+1}-S^\circ_k+1$ and the rank of $u_{k+1}$ among these
children is less than or equal to the
quantity~(\ref{rankchild}).

Let us now discuss scaling limits. We can apply the convergence (\ref
{Levyconv}) to the random walk
$(S^\circ_n)_{n\geq0}$. As a consequence of the results in Chapter 2
of \cite{duqleg02} (see, in particular, Theorem 2.3.2 and Corollary
2.5.1), we have
the joint convergence
%
\begin{equation}
\label{heightconv}
\bigl(n^{-1/\alpha}S^\circ_{[nt]},n^{-(1-1/\alpha)}H^{\circ}_{[nt]}
\bigr)_{t\geq0}
\mathop{\la}_{n\to\infty}^{\mathrm{(d)}} (c_0X_t,c_0^{-1}H_t)_{t\geq0},
\end{equation}
where the convergence holds in distribution, in the Skorokhod sense, and
$(H_t)_{t\geq0}$ is the so-called continuous-time height process
associated with $X$, which may be defined by the limit in probability
\[
H_t=\lim_{\eps\to0} \frac{1}{\eps} \int_0^t \ind_{\{
X_s>I^s_t-\eps\}} \,{d}s.
\]
Note that the preceding approximation of $H_t$ is a continuous analog
of (\ref{heightLuka}). The
process $(H_t)_{t\geq0}$ has continuous sample paths and satisfies the scaling
property
\[
(H_{r t})_{t\geq0}\stackrel{\mathrm{(d)}}{=} (r^{1-1/\alpha}
H_t)_{t\geq0}
\]
for every $r>0$. We refer to \cite{duqleg02}
for a thorough analysis of the continuous-time height process.

We aim to establish a version of (\ref{heightconv}) that
includes the convergence of rescaled labels.
The label process $(L^\circ_n,n\geq0)$ of the forest $\mathbf{F}$
is obtained by concatenating the label processes
$L^{\theta_1}, L^{\theta_2},\ldots$ of the mobiles $\theta_1,\theta
_2,\ldots$ (cf. Section \ref{sec:mobiles-bouttier-di}).
Our goal is to prove the following theorem.
\begin{theorem}
\label{functional}
We have
\[
\bigl( n^{-1/\alpha}S^\circ_{[nt]},n^{-(1-1/\alpha)}H^{\circ
}_{[nt]},n^{-1/2\alpha} L^\circ_{[nt]} \bigr)_{t\geq0}
\mathop{\la}_{n\to\infty}^{{(d)}} \bigl(c_0X_t,c_0^{-1}H_t,\sqrt
{2c_0}D_t\bigr)_{t\geq0},
\]
where the convergence holds in the sense of weak convergence of the
laws in
the Skorokhod space $\D(\R^3)$.
\end{theorem}

The proof of Theorem \ref{functional} is rather long and occupies the
remaining part
of this section. We will first establish the convergence of finite-dimensional
marginals of the rescaled label process and then complete the proof by using
a tightness argument.

\subsection{Finite-dimensional convergence}

\begin{prp}
\label{finitedim}
For every choice of $0\leq t_1<t_2<\cdots<t_p$, we have
\[
n^{-1/2\alpha} \bigl(L^\circ_{[nt_1]},L^\circ_{[nt_2]},\ldots,L^\circ
_{[nt_p]} \bigr)
\mathop{\la}_{n\to\infty}^{{(d)}} \sqrt{2c_0}
(D_{t_1},D_{t_2},\ldots,D_{t_p} ).
\]
Furthermore, this convergence holds jointly with the convergence
(\ref{heightconv}).
\end{prp}
\begin{pf}
In order to write the subsequent arguments in a simpler form, it will
be convenient to use the Skorokhod representation theorem to
replace the convergence in distribution (\ref{heightconv}) by an
almost sure convergence. For every $n\geq1$,
we can construct a labeled forest $\mathbf{F}^{(n)}$ having the same
distribution as~$\mathbf{F}$,
in such a way that if $S^{(n)}$ is the
Lukasiewicz path of $\mathbf{F}^{(n)}$
and $H^{(n)}$ is the height process of $\mathbf{F}^{(n)}$, then we have
the almost sure convergence
%
\begin{equation}
\label{asconv}
\bigl(n^{-1/\alpha}S^{(n)}_{[nt]},
n^{-(1-1/\alpha)}H^{(n)}_{[nt]} \bigr)_{t\geq0}
\mathop{\la}_{n\to\infty}^{\mathrm{(a.s.)}} (c_0 X_t,c_0^{-1}
H_t)_{t\geq0},
\end{equation}
in the sense of the Skorokhod topology. We also use the notation $\F^{(n)}$
for the unlabeled forest associated with $\mathbf{F}^{(n)}$.

We denote by
$u^{(n)}_0,u^{(n)}_1,\ldots$ the white vertices of the forest
$\F^{(n)}$ listed in lexicographical order. For every $k\geq0$, we
denote the label
of $u^{(n)}_k$ by
$L^{(n)}_k=\ell^{(n)}(u^{(n)}_k)$. In order to get the convergence of
one-dimensional marginals in Proposition \ref{finitedim}, we need to
verify that for every $t>0$,
\[
n^{-1/2\alpha} L^{(n)}_{[nt]}\mathop{\la}_{n\to\infty}^{\mathrm{(d)}}
\sqrt{2c_0} D_t.
\]

We fix $t>0$ and $\vep\in(0,1)$.
We denote by $s_i$, $i=1,2,\ldots,$ the sequence consisting of all
times $s\in[0,t]$ such that
\[
X_{s-}<I^s_t.
\]
The times $s_i$ are ranked in such a way that $\Delta X_{s_i}
< \Delta X_{s_j}$ if $i>j$.

On the other hand, let $\mathcal{J}^{(n)}_t$ be the set of all integers
$k\in\{0,1,\ldots,[nt]-1\}$ such that
\[
S^{(n)}_k=\min_{k\leq p\leq[nt]} S^{(n)}_p.
\]
We list the elements of $\mathcal{J}^{(n)}_t$ as $\mathcal{J}^{(n)}_t=\{
a^{(n)}_1,a^{(n)}_2,\ldots,
a^{(n)}_{k_n}\}$, in such a way that
\[
S^{(n)}_{a^{(n)}_i+1}-S^{(n)}_{a^{(n)}_i} \leq
S^{(n)}_{a^{(n)}_j+1}-S^{(n)}_{a^{(n)}_j} \qquad\mbox{if } 1\leq j\leq
i\leq k_n.
\]

The convergence (\ref{asconv}) ensures that almost surely, for every
$i\geq1$,
%
\begin{eqnarray}
\label{convjumptimes}
\frac{1}{n} a^{(n)}_i&\mathop{\la}\limits_{n\to\infty}&
s_i ,\nonumber\\ \frac{1}{c_0
n^{1/\alpha}} \bigl(S^{(n)}_{a^{(n)}_i+1}-S^{(n)}_{a^{(n)}_i} \bigr)
&\mathop{\la}\limits_{n\to\infty}& \Delta X_{s_i} ,\\
\frac{1}{c_0n^{1/\alpha}} \Bigl( \min_{a^{(n)}_i+1\leq k\leq[nt]} S^{(n)}_k
-S^{(n)}_{a^{(n)}_i} \Bigr)&\mathop{\la}\limits_{n\to\infty}& I^{s_i}_t -
X_{s_i-} .\nonumber
\end{eqnarray}

By the observations following (\ref{heightLuka}), we know that the
(white) ancestors of
$u^{(n)}_{[nt]}$ are the vertices $u^{(n)}_k$ for all $k\in
\mathcal{J}^{(n)}_t$. In particular, the generation of $u^{(n)}_{[nt]}$ is
(twice) $H^{(n)}_{[nt]}=\#\mathcal{J}^{(n)}_t$, in
agreement with (\ref{heightLuka}). We can then write
%
\begin{equation}
\label{incredist}
L^{(n)}_{[nt]}= \ell^{(n)}\bigl(u^{(n)}_{[nt]}\bigr) =\sum_{j\in\mathcal{J}^{(n)}_t}
\bigl( \ell^{(n)}\bigl(u^{(n)}_{\vf_n(j)}\bigr)-\ell^{(n)}\bigl(u^{(n)}_j\bigr) \bigr),
\end{equation}
where, for $j\in\mathcal{J}^{(n)}_t$, $\vf_n(j)$ is the smallest
element of $(\{j+1,\ldots,[nt]-1\}\cap\mathcal{J}^{(n)}_t)\cup\{[nt]\}$.
Equivalently,
$u^{(n)}_{\vf_n(j)}$ is the unique (white)
grandchild of $u^{(n)}_j$ that is also an ancestor of
$u^{(n)}_{[nt]}$.

Now, consider the L\'evy process $X$. As a consequence of classical
results of fluctuation theory (see, e.g., Lemma 1.1.2 in \cite{duqleg02}),
we know that the ladder height process of $X$ is a subordinator
without drift, hence a pure jump process. By applying this to the dual process
$(X_{(t-r)-}-X_t,0\leq r<t)$, we obtain that
\[
X_t-I_t=\sum_{i=1}^\infty(I^{s_i}_t - X_{s_i-}).
\]
It follows that we can fix an integer $N\geq1$ such that, with
probability greater
than $1-\vep,$ we have
%
\begin{equation}
\label{bound0}
X_t-I_t-\sum_{i=1}^N (I^{s_i}_t - X_{s_i-})=\sum_{i>N} (I^{s_i}_t -
X_{s_i-}) \leq\frac{\vep}{2} .
\end{equation}

Now, note that
\[
\frac{1}{c_0 n^{1/\alpha}} \Bigl(S^{(n)}_{[nt]}-\min_{k\leq[nt]}
S^{(n)}_k \Bigr) \mathop{\la}_{n\to\infty}^{\mathrm{a.s.}} X_t - I_t
\]
and
recall the convergences (\ref{convjumptimes}). Using (\ref{bound0}),
it follows that we can find $n_0$ sufficiently large such that for
every $n\geq n_0$, with probability greater than $1-2\vep$, we have
\[
\frac{1}{c_0 n^{1/\alpha}} \Biggl( \Bigl(S^{(n)}_{[nt]}-\min_{k\leq
[nt]} S^{(n)}_k \Bigr) -\sum_{i=1}^{N\wedge k_n} \Bigl( \min_{a^{(n)}_i+1\leq
k\leq
[nt]} S^{(n)}_k -S^{(n)}_{a^{(n)}_i} \Bigr) \Biggr) < \vep.
\]
Since
\[
\sum_{i=1}^{k_n} \Bigl(
\min_{a^{(n)}_i+1\leq k\leq[nt]} S^{(n)}_k
-S^{(n)}_{a^{(n)}_i} \Bigr) = S^{(n)}_{[nt]}-\min_{k\leq[nt]} S^{(n)}_k ,
\]
we get that, for every $n\geq n_0$, with probability greater than
$1-2\vep$,
%
\begin{equation}
\label{bound1}
\frac{1}{c_0 n^{1/\alpha}}\sum_{i>N} \Bigl( \min_{a^{(n)}_i+1\leq
k\leq[nt]} S^{(n)}_k -S^{(n)}_{a^{(n)}_i} \Bigr) < \vep.
\end{equation}

Now, recall (\ref{incredist}). By Proposition \ref
{sec:boutt-di-franc} and the observations
following (\ref{heightLuka}), we know that, conditionally on the forest
$\F^{(n)}$, for every $j\in\mathcal{J}^{(n)}_t$, the quantity
\[
\ell^{(n)}\bigl(u^{(n)}_{\vf_n(j)}\bigr)-\ell^{(n)}\bigl(u^{(n)}_j\bigr)
\]
is distributed as the value of a discrete bridge with length
$p_j\leq S^{(n)}_{j+1}- S^{(n)}_j+2$, at a time
$k_j\leq S^{(n)}_{j+1}-\min_{j+1\leq k\leq[nt]} S^{(n)}_k+1$
such that $p_j-k_j\leq\min_{j+1\leq k\leq[nt]} S^{(n)}_k-
S^{(n)}_j+1$. Thanks
to our estimate (\ref{momentbridge2}) on discrete bridges, we thus
have
\begin{eqnarray*}
\E\bigl[\bigl(\ell^{(n)}\bigl(u^{(n)}_{\vf_n(j)}\bigr)-\ell^{(n)}\bigl(u^{(n)}_j\bigr)\bigr)^2\mid
\F^{(n)}\bigr] &\leq& K \frac{k_j(p_j-k_j)}{p_j}\\
&\leq& K \Bigl(\min_{j+1\leq
k\leq[nt]} S^{(n)}_k- S^{(n)}_j+1 \Bigr).
\end{eqnarray*}
Furthermore, still
conditionally on the forest $\F^{(n)}$, the random variables\break
$\ell^{(n)}(u^{(n)}_{\vf_n(j)})-\ell^{(n)}(u^{(n)}_j)$ are independent
and centered. It follows that for $n\geq n_0$,
\begin{eqnarray*}
&&\E\biggl[ \biggl(n^{-1/2\alpha}\sum_{j\in\mathcal{J}^{(n)}_t\setminus
\{a^{(n)}_1,\ldots,a^{(n)}_N\}}
\bigl(\ell^{(n)}\bigl(u^{(n)}_{\vf_n(j)}\bigr)-\ell^{(n)}\bigl(u^{(n)}_j\bigr)\bigr) \biggr)^2
\Big| \F^{(n)} \biggr]\\
&&\qquad \leq K n^{-1/\alpha}
\sum_{j\in\mathcal{J}^{(n)}_t\setminus\{a^{(n)}_1,\ldots,a^{(n)}_N\}}
\Bigl(\min_{j+1\leq k\leq[nt]} S^{(n)}_k-
S^{(n)}_j+1 \Bigr)\\
&&\qquad \leq K \bigl(c_0\vep+ n^{-1/\alpha}\#
\mathcal{J}^{(n)}_t\bigr),
\end{eqnarray*}
the last bound holding on a set of probability greater than $1-2\vep$,
by (\ref{bound1}). Since $\# \mathcal{J}^{(n)}_t=H^{(n)}_{[nt]}$, we have
$n^{-1/\alpha}\# \mathcal{J}^{(n)}_t \la0$ a.s. as $n\to\infty$, by
(\ref{asconv}).

From (\ref{incredist}) and the preceding considerations, the limiting
behavior of $n^{-1/2\alpha}\times L^{(n)}_{[nt]}$
will follow from that of
\[
n^{-1/2\alpha}\sum_{j\in\{a^{(n)}_1,\ldots,a^{(n)}_N\}} \bigl(
\ell^{(n)}\bigl(u^{(n)}_{\vf_n(j)}\bigr)-\ell^{(n)}\bigl(u^{(n)}_j\bigr) \bigr).
\]
Recall that for every
$j\in\{a^{(n)}_1,\ldots,a^{(n)}_N\}$, the number of white grandchildren
of $u^{(n)}_j$ in the forest $\F^{(n)}$ is
$m^{(n)}_j=S^{(n)}_{j+1}-S^{(n)}_j+1$. Moreover, $u^{(n)}_{\vf_n(j)}$
appears at the rank
\[
r^{(n)}_j=S^{(n)}_{j+1}-\min_{j+1\leq k\leq[nt]} S^{(n)}_k +1
\]
in the list of these grandchildren.
The next lemma will imply that $u^{(n)}_{\vf_n(j)}$ is the child
of a black vertex whose number of children is also close to $m^{(n)}_j$.
\begin{lmm}
\label{keyestimate}
We can choose $\delta>0$ small enough so that, for every fixed $\eta
>0$, the following holds with probability close to $1$ when $n$ is
large. For every white vertex belonging to
$\{u^{(n)}_0,u^{(n)}_1,\ldots,u^{(n)}_{[nt]}\}$ that has more than
$\eta n^{1/\alpha}$ white grandchildren in the forest $\F^{(n)}$, all these
grandchildren have the same (black) parent in the forest $\F^{(n)}$,
except for
at most $n^{1/\alpha-\delta}$ of them.
\end{lmm}
\begin{pf}
Recall that $\mu_0(k)=\beta(1-\beta)^k$
for every $k\geq0$. We choose $\delta>0$ such that $2\delta\alpha
<1$ and
take $n$ sufficiently large so that $\eta n^{1/\alpha}>2n^{1/\alpha
-\delta}$.
Let us fix $i\in\{0,1,\ldots,[nt]\}$.
The number of black children of $u^{(n)}_i$ is
distributed according to $\mu_0$ and each of these black children
has a number of white children distributed according to $\mu_1$.
Supposing that $u^{(n)}_i$ has $k$ black children,
if it has a number $M\geq\eta n^{1/\alpha}$ of grandchildren
and simultaneously none of its black children has more than
$M-n^{1/\alpha-\delta}$ white children, this implies that at least two
among its black children will have more than $n^{1/\alpha-\delta}/k$
white children. The probability that this occurs is bounded above by
\[
\beta\sum_{k=2}^\infty(1-\beta)^k \pmatrix{k\cr2} \ov\mu
_1(n^{1/\alpha-\delta}/k)^2.
\]
From (\ref{tail1}), there is a constant $K$ such that $\ov\mu
_1(k)\leq K k^{-\alpha}$
for every $k\geq1$. Hence, the last displayed quantity is bounded by
\[
K^2\beta\Biggl(\sum_{k=2}^\infty(1-\beta)^k \pmatrix{k\cr2} k^{2\alpha} \Biggr)
n^{-2+2\delta\alpha} =o(n^{-1}).
\]
The desired result follows by summing this estimate over $i\in\{
0,1,\ldots,[nt]\}$.
\end{pf}

We return to the proof of Proposition \ref{finitedim}.
We fix $\delta>0,$ as in the lemma. We first observe that for every
$j\in\{a^{(n)}_1,\ldots,a^{(n)}_N\}$,
(\ref{convjumptimes}) implies that
\[
\lim_{n\to\infty} n^{-1/\alpha}r^{(n)}_j= c_0(X_{s_j} - I^{s_j}_t) >0.
\]
We then
deduce from Lemma \ref{keyestimate} that, with a probability close to
$1$ when
$n$ is large, for every $j\in\{a^{(n)}_1,\ldots,a^{(n)}_N\}$,
$u^n_{\vf_n(j)}$ is the child of a black child of $u^{(n)}_j$, whose
number of white children is $\wt m^{(n)}_j$ such that
%
\begin{equation}
\label{estim1}
m^{(n)}_j\geq\wt m^{(n)}_j\geq m^{(n)}_j - n^{1/\alpha-\delta}.
\end{equation}
Moreover, the rank $\wt r^{(n)}_j$ of $u^n_{\vf_n(j)}$ among the
children of its (black) parent satisfies
%
\begin{equation}
\label{estim2}
r^{(n)}_j\geq\wt r^{(n)}_j\geq r^{(n)}_j - n^{1/\alpha-\delta}.
\end{equation}
On the other hand, we know that, conditionally on $\F^{(n)}$, the difference
\[
\ell^{(n)}\bigl(u^{(n)}_{\vf_n(j)}\bigr)-\ell^{(n)}\bigl(u^{(n)}_j\bigr)
\]
is
distributed as the value of a discrete bridge with length $\wt
m^{(n)}_j +1$ at time $\wt r^{(n)}_j$. Thus, conditionally on
$\F^{(n)}$,
\[
\sum_{j\in\{a^{(n)}_1,\ldots,a^{(n)}_N\}} \bigl(
\ell^{(n)}\bigl(u^{(n)}_{\vf_n(j)}\bigr)-\ell^{(n)}\bigl(u^{(n)}_j\bigr) \bigr)
\stackrel{\mathrm{(d)}}{=} \sum_{i=1}^N b_i^{(n)}\bigl(\wt r^{(n)}_{a_i^{(n)}}\bigr),
\]
where, for every $i\in\{1,\ldots,N\}$, $b^{(n)}_i$ is a discrete
bridge with length $\wt m^{(n)}_{a^{(n)}_i}+1$ and the bridges
$b^{(n)}_i$ are independent.

Using Donsker's theorem for bridges (\ref{Donsker}), the convergences
(\ref{asconv}) and (\ref{convjumptimes}) and the bounds (\ref{estim1})
and (\ref{estim2}), together with scaling properties of Brownian
bridges, it is then a simple matter to obtain that, for every
$i\in\{1,\ldots,N\}$,
%
\begin{equation}
\label{finitedimtech0}
n^{-1/2\alpha} b^{(n)}_i\bigl(\wt r^{(n)}_{a^{(n)}_i}\bigr)
\mathop{\la}_{n\to\infty}^{\mathrm{(d)}}
\sqrt{2c_0} \gamma_i(X_{s_i}-I^{s_i}_t),
\end{equation}
where, conditionally on
$X$, $\gamma_i=(\gamma_i(r))_{0\leq r\leq\Delta X_{s_i}}$ is a Brownian
bridge with length $\Delta X_{s_i}$. The preceding convergences hold
jointly when $i$ varies in $\{1,\ldots, N\}$ with Brownian bridges
$\gamma_1,\ldots,\gamma_N$ that are independent conditionally on
$X$. Finally, it follows that
\[
n^{-1/2\alpha} \sum_{j\in\{a^{(n)}_1,\ldots,a^{(n)}_N\}} \bigl(
\ell^{(n)}\bigl(u^{(n)}_{\vf_n(j)}\bigr)-\ell^{(n)}(u^{(n)}_j) \bigr)
\mathop{\la}_{n\to\infty}^{\mathrm{(d)}}
\sqrt{2c_0} \sum_{i=1}^N \gamma_i(X_{s_i}-I^{s_i}_t).
\]
From Proposition \ref{defdistance}, the limit is close to $\sqrt
{2c_0}D_t$ when $N$ is
large. This completes the proof of the convergence of one-dimensional
marginals. It is also clear from our argument that the convergences
(\ref{finitedimtech0}) hold jointly with (\ref{asconv}), so the convergence
of $n^{-1/2\alpha}L^\circ_{[nt]}$ must hold jointly with (\ref{heightconv}).

The same arguments yield the convergence of finite-dimensional
marginals. It would be tedious to write a detailed proof,
but we sketch the method in the case of two-dimensional
marginals. So, fix $0<s<t$. We aim to prove that
\[
n^{-1/2\alpha}
\bigl(L^{(n)}_{[ns]},L^{(n)}_{[nt]}\bigr)\mathop{\la}_{n\to\infty}^{\mathrm{(d)}}
\sqrt{2c_0} (D_s,D_t).
\]
It is convenient to argue separately on
the events $\{I_s>I_t\}$ and $\{I_s=I_t\}$. Discarding sets of
probability zero,
the first event corresponds to the case where $s$ and $t$ belong to
different excursion intervals of $X-I$ away from $0$ and the second
event corresponds to the case where $s$ and $t$ are in the same
excursion interval
of $X-I$.

On the event $\{I_s>I_t\}$, things are easy. We first note that,
conditionally on $X$, $D_s$ and $D_t$ are independent on that
event. This is the case because the jumps $t_i$ that give
a nonzero contribution in (\ref{eq:1}) belong to the
excursion interval of $X-I$ that straddles $t$. Similarly, $L^n_{[ns]}$
and $L^n_{[nt]}$
are independent, conditionally given the forest $\F^{(n)}$, on the event
\[
\min_{k\leq[ns]}S^{(n)}_{k} >\min_{k\leq[nt]} S^{(n)}_{k}.
\]
Furthermore, the latter event converges to $\{I_s>I_t\}$ as
$n\to\infty$. Thus, the very same arguments as in the case of
one-dimensional marginals yield that the conditional distribution
of the pair $n^{-1/2\alpha}
(L^{(n)}_{[ns]},L^{(n)}_{[nt]})$ given
$\{I_s>I_t\}$ converges to the conditional distribution of
$\sqrt{2c_0} (D_s,D_t)$ given the same event.

On the event $\{I_s=I_t\}$, we need to be a little more careful. Set
\begin{eqnarray*}
\mathcal{J}_s&=&\{r\in[0,s]\dvtx X_{r-}<I^r_s\} ,\\
\mathcal{J}_t&=&\{r\in[0,t]\dvtx X_{r-}<I^r_t\} .
\end{eqnarray*}
Then, a.s. there exists a unique $r_0\in\mathcal{J}_s$ such that
\[
I^s_t\in(X_{r_0-},I^{r_0}_s).
\]
Furthermore, we have $\mathcal{J}_s\cap
\mathcal{J}_t=\mathcal{J}_s\cap[0,r_0]=\mathcal{J}_t\cap[0,r_0]$ and
$I^r_s=I^r_t$ for every
$r\in\mathcal{J}_s\cap[0,r_0)$. Using the convergence (\ref{asconv}),
we get
that, a.s. on the event $\{I_s=I_t\}$, for $n$ sufficiently large,
there exists a time $j_0(n)\in\mathcal{J}^{(n)}_{s}\cap\mathcal
{J}^{(n)}_t$ such that
\[
S^{(n)}_{j_0(n)}<\min_{[ns]\leq k\leq[nt]} S^{(n)}_k < \min
_{j_0(n)+1\leq
k\leq[ns]} S^{(n)}_k < S^{(n)}_{j_0(n)+1}
\]
and, furthermore,
$\mathcal{J}^{(n)}_s\cap\mathcal{J}^{(n)}_t=\mathcal{J}^{(n)}_s\cap
[0,j_0(n)]=\mathcal{J}^{(n)}_t\cap
[0,j_0(n)]$. The white vertices that are common ancestors to
$u^{(n)}_{[ns]}$ and $u^{(n)}_{[nt]}$ are exactly the vertices
$u^{(n)}_k$ for $k\in\mathcal{J}^{(n)}_s\cap[0,j_0(n)]$. Also, note that
$n^{-1}j_0(n)$ converges to $r_0$, a.s. on the event $\{I_s=I_t\}$.

Write $\psi_n\dvtx\mathcal{J}^{(n)}_s\la\mathcal{J}^{(n)}_s\cup\{[ns]\}$
for the function
analogous to $\vf_n$ when $t$ is replaced by $s$. Analogously to
(\ref{incredist}), we have
\begin{eqnarray*}
L^{(n)}_{[ns]}&=&\sum_{j\in\mathcal{J}^{(n)}_s} \bigl(
\ell^{(n)}\bigl(u^{(n)}_{\psi_n(j)}\bigr)-\ell^{(n)}\bigl(u^{(n)}_j\bigr) \bigr),\\
L^{(n)}_{[nt]}&=&\sum_{j\in\mathcal{J}^{(n)}_t} \bigl(
\ell^{(n)}\bigl(u^{(n)}_{\vf_n(j)}\bigr)-\ell^{(n)}\bigl(u^{(n)}_j\bigr) \bigr).
\end{eqnarray*}
The terms corresponding to $j\in\mathcal{J}^{(n)}_s\cap
[0,j_0(n))=\mathcal{J}
^{(n)}_t\cap
[0,j_0(n))$ are the same in both sums of the preceding display. On the
other hand, conditionally on $\F^{(n)}$, the terms corresponding
to $j\in\mathcal{J}^{(n)}_s\cap(j_0(n),[ns])$ in the first sum are
independent of the terms of the second sum and similarly for the
terms corresponding to $j\in\mathcal{J}^{(n)}_t\cap(j_0(n),[nt])$ in the
second sum. As for the term corresponding to $j_0(n)$, the same
arguments as in the proof of the convergence of one-dimensional
marginals, using Lemma \ref{keyestimate} in particular, show that
\begin{eqnarray*}
&&n^{-1/2\alpha} \bigl(\ell^{(n)}\bigl(u^{(n)}_{\psi_n(j_0(n))}\bigr)-
\ell^{(n)}\bigl(u^{(n)}_{j_0(n)}\bigr),
\ell^{(n)}\bigl(u^{(n)}_{\vf_n(j_0(n))}\bigr)-\ell^{(n)}\bigl(u^{(n)}_{j_0(n)}\bigr) \bigr)\\
&&\qquad
\mathop{\la}_{n\to\infty}^{\mathrm{(d)}}
\sqrt{2c_0} \bigl(\gamma(X_{r_0}-I^{r_0}_s),\gamma(X_{r_0}-I^{r_0}_t) \bigr),
\end{eqnarray*}
where, conditionally given $X$, $\gamma$ is a Brownian bridge with
length $\Delta X_{r_0}$.

Finally, let $(r_i)_{i\in\N}$ be a measurable enumeration
of $\mathcal{J}_s\cap
[0,r_0)=\mathcal{J}_t\cap[0,r_0)$, $(r'_i)_{i\in\N}$ a measurable
enumeration of
$\mathcal{J}_s\cap(r_0,s]$ and $(r''_i)_{i\in\N}$ a measurable
enumeration of
$\mathcal{J}_s\cap(r_0,t]$. Set
\begin{eqnarray*}
L^\infty_s&=&\sum_{i\in\N} \gamma_i(X_{r_i}-I^{r_i}_s) + \gamma
(X_{r_0}-I^{r_0}_s)
+ \sum_{i\in\N} \gamma'_i(X_{r'_i}-I^{r'_i}_s),\\
L^\infty_t&=&\sum_{i\in\N} \gamma_i(X_{r_i}-I^{r_i}_t) + \gamma
(X_{r_0}-I^{r_0}_t)
+ \sum_{i\in\N} \gamma''_i(X_{r''_i}-I^{r''_i}_t),
\end{eqnarray*}
where, conditionally given $X$, $(\gamma_i)_{i\in\N}$, $(\gamma
'_i)_{i\in\N}$, $(\gamma''_i)_{i\in\N}$
and $\gamma$ are independent Brownian bridges and the duration of
$\gamma_i$
(resp., $\gamma'_i$, $\gamma''_i$) is $\Delta X_{r_i}$ (resp., $\Delta
X_{r'_i}$,
$\Delta X_{r''_i}$).
Then, by following the lines of the proof of the convergence
of one-dimensional marginals, we obtain that the conditional distribution
of $n^{-1/2\alpha}
(L^{(n)}_{[ns]},L^{(n)}_{[nt]})$ given
$\{I_s=I_t\}$ converges to the conditional distribution of $\sqrt
{2c_0} (L^\infty_s,L^\infty_t)$ given the same event. However, the
latter conditional distribution clearly coincides with the
conditional distribution of $\sqrt{2c_0} (D_s,D_t)$ given
$\{I_s=I_t\}$.
So, we get the desired convergence for two-dimensional marginals and the
same argument as in the case of one-dimensional marginals gives a
joint convergence with (\ref{heightconv}).
This completes the proof.
\end{pf}

\subsection{Tightness of the rescaled label process}

The next proposition will allow us to complete the proof of Theorem
\ref{functional}.
\begin{prp}
\label{keybound}
There exists a constant $K_0$ such that, for all integers $i$, \mbox{$j\geq
0$},
\[
\E[(L^\circ_i-L^\circ_j)^4]\leq K_0 |i-j|^{2/\alpha}.
\]
\end{prp}

Theorem \ref{functional} is an easy consequence of this proposition
and Proposition \ref{finitedim}. To see this,
define $L^{\{n\}}_t=n^{-1/2\alpha} L^\circ_{nt}$ if $nt$ is an
integer and
use linear approximation to define $L^{\{n\}}_t$ for every real $t\geq
0$. By the bound of the proposition,
\[
\E\bigl[\bigl( L^{\{n\}}_s-L^{\{n\}}_t\bigr)^4\bigr]\leq K_0
|s-t|^{2/\alpha},
\]
if $ns$ and
$nt$ are both integers. It readily follows that the same bound holds
(possibly with a different constant) for all reals $s,t\geq
0$. Since $2/\alpha>1$, standard criteria entail that the sequence of
the distributions of the processes $L^{\{n\}}$ is tight in the space of
probability measures on $C(\R)$. Theorem \ref{functional} then
follows by using Proposition \ref{finitedim}.
\begin{pf*}{Proof of Proposition \ref{keybound}} We use the same
notation as
in Section
\ref{sec:conv-height-proc}. In particular, $u_0,u_1,u_2,\ldots$ are
the white vertices of the
forest $\F$ listed in lexicographical order and one tree after
another, so $L^\circ_i=\ell(u_i)$ is the label of $u_i$. We
also set
\[
\mathcal{J}(i)=\Bigl\{k\in\{0,1,\ldots, i-1\}\dvtx S^\circ_k\leq\min_{k+1\leq\ell
\leq i} S^\circ_\ell\Bigr\}
\]
in such a way that the vertices $u_k, k\in\mathcal{J}(i)$
are the white vertices of $\F$ that are strict ancestors of $u_i$.

We fix two nonnegative integers $i<j$. If $k\in\mathcal{J}(i)$, then we write
$\vf(k)$ for the
index such that $u_{\vf(k)}$ is the (unique) grandchild of $u_k$ that is
also an ancestor of $u_i$. We similarly define $\psi(k)$ for
$k\in\mathcal{J}(j)$ in such a way that $u_{\psi(k)}$ is the grandchild
of $u_k$
that is an ancestor of $u_j$.

In the case where $u_i$ and $u_j$ belong to the same tree of the
forest, we define $i_0$ by requiring that $u_{i_0}$ is the most recent white
common ancestor of $u_i$ and $u_j$ in $\F$. If $i_0<i$, then we have
%
\begin{equation}
\label{tech00}
S^\circ_{i_0}\leq\min_{i\leq k\leq j}S^\circ_k\leq S^\circ_{\vf(i_0)}.
\end{equation}
This easily follows from the relations
between the sequence $\mathcal{T}^\circ_1,\mathcal{T}^\circ_2,\ldots$
and the
Lukasiewicz path
$S^\circ$ (see, e.g., \cite{duqleg02}, Section 0.2, or \cite{Trees},
Section 1) and we leave the proof as an exercise for the reader.
It may happen that $i_0=i$
(but not that $i_0=j$) and, in that case, we set $\vf(i_0)=i_0$, by convention.

In the case where $u_i$ and $u_j$ belong to different trees of the
forest, we take $i_0=-\infty$, by convention, and we also agree that
$\vf(-\infty)$ [resp., $\psi(-\infty)$] is defined in such a way that
$u_{\vf(-\infty)}$ [resp., $u_{\psi(-\infty)}$] is the root of the tree
containing $u_i$ (resp., containing $u_j$).

We then have
%
\begin{eqnarray}
\label{decompos}
L^\circ_i-L^\circ_j&=&\ell(u_i)-\ell(u_j)\nonumber\\
&=&\sum_{k\in\mathcal{J}(i)\cap
(i_0,i)} \bigl(\ell\bigl(u_{\vf(k)}\bigr) -\ell(u_k) \bigr)\nonumber\\[-8pt]\\[-8pt]
&&{}-
\sum_{k\in\mathcal{J}(j)\cap(i_0,j)} \bigl( \ell\bigl(u_{\psi(k)}\bigr)
-\ell(u_k) \bigr)\nonumber\\
&&{} + \ell\bigl(u_{\vf(i_0)}\bigr) -\ell\bigl(u_{\psi
(i_0)}\bigr).\nonumber
\end{eqnarray}

As in the proof of
Proposition \ref{finitedim}, we can write
\[
\sum_{k\in\mathcal{J}(i)\cap(i_0,i)} \bigl( \ell\bigl(u_{\vf(k)}\bigr) -\ell(u_k) \bigr)
= \sum_{k\in\mathcal{J}(i)\cap(i_0,i)} b_k(r_k),
\]
where,
conditionally on $\F$, the
processes $b_k$ are independent discrete bridges, $b_k$ has length
$m_k\leq S^\circ_{k+1}-S^\circ_k+2$ and $r_k\in\{1,\ldots,m_k-1\}$
is such that
%
\begin{eqnarray}
\label{techbound1}
r_k&\leq& S^\circ_{k+1}- \min_{k+ 1\leq\ell\leq i} S^\circ_\ell+1
,\\
\label{techbound2}
m_k-r_k&\leq&\min_{k+1 \leq\ell\leq i} S^\circ_\ell-S^\circ_k+1 .
\end{eqnarray}

From the bound of Lemma \ref{estimatebridge} and
(\ref{techbound2}), we get, with some constant $K_1$,
\begin{eqnarray*}
\E\biggl[ \biggl(\sum_{k\in\mathcal{J}(i)\cap(i_0,i)} b_k(r_k) \biggr)^4 \Big| \F\biggr]
&\leq& K_1 \biggl(\sum_{k\in\mathcal{J}(i)\cap(i_0,i)}
(m_k-r_k) \biggr)^2\\
&\leq& K_1 \biggl(\sum_{k\in\mathcal{J}(i)\cap(i_0,i)}
\Bigl(\min_{k+1 \leq\ell\leq i} S^\circ_\ell-S^\circ_k+1 \Bigr) \biggr)^2\\
&\leq& 2K_1
\Bigl( \Bigl(S^\circ_i-\min_{i\leq\ell\leq j} S^\circ_\ell\Bigr)^2 +
\Bigl(H^\circ_i-\min_{i\leq\ell\leq j} H^\circ_\ell\Bigr)^2 \Bigr).
\end{eqnarray*}
In the last inequality, we have used the identity
\[
\#\{k\in\mathcal{J}(i)\cap(i_0,i)\}= H^\circ_i-\min_{i\leq\ell\leq j}
H^\circ_\ell
\]
and the bound
\[
\sum_{k\in\mathcal{J}(i)\cap(i_0,i)} \Bigl(\min_{k+1 \leq\ell\leq i} S^\circ
_\ell
-S^\circ_k \Bigr) \leq S^\circ_i-\min_{i\leq\ell\leq j} S^\circ_\ell,
\]
which follows from (\ref{tech00}) in the case $i_0<i$.

To simplify notation, set
\[
J_n=\min_{0\leq k\leq n} S^\circ_k
\]
and note that
\[
S^\circ_i-\min_{i\leq\ell\leq j} S^\circ_\ell\stackrel{\mathrm{(d)}}{=} -J_{j-i}.
\]
\begin{lmm}
\label{firstbound}
There exists a constant $K_2$ such that, for every integer $n\geq1$,
\[
\E[(J_n)^2]\leq K_2 n^{2/\alpha}.
\]
\end{lmm}
\begin{lmm}
\label{secondbound}
There exists a constant $K_3$, which does not depend on
the choice of $i$ and $j$, such that
\[
\E\Bigl[ \Bigl(H^\circ_i+H^\circ_j-2\min_{i\leq\ell\leq j} H^\circ_\ell
\Bigr)^2 \Bigr]\leq
K_3 |i-j|^{2(1-1/\alpha)}.
\]
\end{lmm}

The proof of these lemmas is postponed to the end of the section.
By combining Lemmas \ref{firstbound}, \ref{secondbound} and the
previous observations, we get, with a certain constant $K_4$,
\[
\E\biggl[ \biggl(\sum_{k\in\mathcal{J}(i)\cap(i_0,i)} \bigl( \ell\bigl(u_{\vf(k)}\bigr) -\ell(u_k)
\bigr) \biggr)^4 \biggr]
\leq K_4 |i-j|^{2/\alpha}.
\]

We still have to treat the other two terms in the right-hand side of
(\ref{decompos}). As previously, we have
\[
\sum_{k\in\mathcal{J}(j)\cap(i_0,j)} \bigl( \ell\bigl(u_{\psi(k)}\bigr) -\ell(u_k) \bigr)
=\sum_{k\in\mathcal{J}(j)\cap(i_0,j)} b_k(r_k),
\]
where, conditionally on
$\F$,
the processes $b_k$ are independent discrete bridges, $b_k$ has length
$m_k\leq S^\circ_{k+1}-S^\circ_k+2$ and $r_k\in\{1,\ldots,m_k-1\}$
satisfies the
bounds (\ref{techbound1}) and (\ref{techbound2}) with $i$ replaced by $j$.
Arguing as above, but now using the bound (\ref{techbound1}), we get
\begin{eqnarray*}
&&\E\biggl[ \biggl(\sum_{k\in\mathcal{J}(j)\cap(i_0,j)}
b_k(r_k) \biggr)^4 \Big| \F\biggr] \\
&&\qquad \leq
2K_1 \biggl( \biggl(\sum_{k\in\mathcal{J}(j)\cap(i_0,j)} \Bigl(S^\circ_{k+1}- \min_{k+
1\leq
\ell\leq j} S^\circ_\ell\Bigr) \biggr)^2 + \Bigl(H^\circ_j-\min_{i\leq\ell\leq j}
H^\circ_\ell\Bigr)^2 \biggr).
\end{eqnarray*}
The expected value of the second term in the right-hand side is
bounded by Lem\-ma~\ref{secondbound}. As for the first term, we
observe that $\mathcal{J}(j)\cap(i_0,j)=\mathcal{J}(j)\cap(i,j)$ and thus
\begin{eqnarray*}
&&\sum_{k\in\mathcal{J}(j)\cap(i_0,j)} \Bigl(S^\circ_{k+1}- \min_{k+ 1\leq\ell
\leq j} S^\circ_\ell\Bigr)\\
&&\qquad =\sum_{k\in(i,j)} \ind_{\{S^\circ_k\leq\min_{k+1\leq\ell\leq
j} S^\circ_\ell\}}
\Bigl(S^\circ_{k+1}- \min_{k+ 1\leq\ell\leq j} S^\circ_\ell\Bigr)
\stackrel{\mathrm{(d)}}{=} F_{j-i-1},
\end{eqnarray*}
where, for every $n\geq1$,
\[
F_n=\sum_{k=0}^{n-1} \ind_{\{ S^\circ_k\leq\min_{k+1\leq\ell\leq n}
S^\circ_\ell\}} \Bigl(S^\circ_{k+1}-\min_{k+1\leq\ell\leq n} S^\circ
_\ell\Bigr).
\]
Furthermore, a time reversal argument shows that $F_n$ has the same
distribution as~$G_n$, where
\[
G_n=\sum_{k=1}^{n} \ind_{\{ S^\circ_{k}\geq\max_{0\leq\ell\leq k-1}
S^\circ_\ell\}} \Bigl(\max_{0\leq\ell\leq k-1} S^\circ_\ell-S^\circ
_{k-1} \Bigr).
\]
\begin{lmm}
\label{thirdbound}
There exists a constant $K_5$ such that, for every integer $n\geq1$,
\[
\E[(G_n)^2]\leq K_5 n^{2/\alpha}.
\]
\end{lmm}

Combining Lemma \ref{thirdbound} with the preceding observations, we
see that the fourth moment of the second term in the right-hand side
of (\ref{decompos}) is bounded above by $K_6|j-i|^{2/\alpha}$ for
some constant $K_6$. We easily get the same bound for the third
term by using Lemmas \ref{estimatebridge} and
\ref{firstbound}. This completes the proof of Proposition
\ref{keybound}, but we still have to prove Lemmas \ref{firstbound},
\ref{secondbound} and \ref{thirdbound}.
\end{pf*}
\begin{pf*}{Proof of Lemma \ref{firstbound}} For every integer $k\geq0$, set
\[
V_k=\inf\{n\geq0\dvtx S^\circ_n=-k\}.
\]
Note that $V_k$ is the sum of $k$
independent copies of $V_1$. As a consequence of (\ref{Levyconv}),
$n^{-\alpha}V_n$ converges in distribution to the variable
$T_{c_0^{-1}}=\inf\{t\geq0\dvtx X_t<-c_0^{-1}\}$, which is stable with
index $1/\alpha$. By standard results concerning domains of attraction of
stable distributions (see, e.g., Section XVII.5 of \cite{feller2}),
there exists a constant $K>0$ such that
%
\begin{equation}
\label{equiv-hitting}
\mathbb{P}(V_1> n)\mathop{\sim}_{n\to\infty} K n^{-1/\alpha}.
\end{equation}
Consequently, there is a constant $K'>0$ such that, for every $n\geq1$,
\[
\mathbb{P}(V_1> n)\geq K'n^{-1/\alpha}.
\]
Then, for every $x\geq1$ and
$n\geq1$,
\begin{eqnarray*}
\mathbb{P}(|J_n|\geq xn^{1/\alpha})&\leq&\mathbb{P}\bigl(V_{[xn^{1/\alpha
}]}\leq n\bigr) \\
&\leq&
\mathbb{P}(V_1\leq n)^{[xn^{1/\alpha}]}\leq
(1-K'n^{-1/\alpha})^{[xn^{1/\alpha}]} \\
&\leq&
\exp(-K'x/2).
\end{eqnarray*}
It readily
follows that all moments of $n^{-1/\alpha}|J_n|$ are uniformly
bounded.
\end{pf*}
\begin{pf*}{Proof of Lemma \ref{secondbound}}
For all nonnegative integers $k\leq\ell$, we
set $J_{k,\ell}=\min_{k\leq
n\leq\ell} S^\circ_n$ so that $J_k=J_{0,k}$. We fix two nonnegative
integers $i< j$ and
first look for an expression of $\min_{i\leq\ell\leq j} H^\circ
_\ell$. To this end, we set
\[
g=\max\bigl\{r\in\{0,1,\ldots,i-1\}\dvtx S^\circ_r\leq J_{i,j}\bigr\}
\]
with the convention that $\max\varnothing=-\infty$. First, assume
that $g>-\infty$
and let
$k\in\{i,\ldots,j\}$. We then have
%
\begin{eqnarray}
\label{eq:3}\hspace*{30pt}
H^\circ_k&=&\#\bigl\{\ell\in\{0,\ldots,k-1\}\dvtx S^\circ_\ell=J_{\ell,k}\bigr\}
\nonumber\\[-8pt]\\[-8pt]
&=&\#\bigl\{\ell\in\{0,\ldots,g-1\}\dvtx S^\circ_\ell=J_{\ell,k}\bigr\}+
\#\bigl\{\ell\in\{g,\ldots,k-1\}\dvtx S^\circ_\ell=J_{\ell,k}\bigr\}.\nonumber
\end{eqnarray}
From the definition of $g$, it is easy to verify that
$J_{\ell,k}=J_{\ell,g}$ for every $\ell\in\{0,\ldots,g-1\}$. Thus,
the first term
in the right-hand side of (\ref{eq:3}) is equal to
$H^\circ_g$ and does not depend on $k$.
We then note that $S^\circ_g=J_{g,k}$, by the definition of $g$, so
the second term in the right-hand side of
(\ref{eq:3})
equals
\[
1+\#\bigl\{\ell\in\{g+1,\ldots,k-1\}\dvtx S^\circ_\ell=J_{\ell,k}\bigr\}.
\]
This expression attains its minimal value $1$ when $k$ equals
$\min\{\ell\geq i\dvtx S^\circ_\ell=J_{i,j}\}$. Thus, we have proved,
when $g>-\infty$,
that
\[
\min_{i\leq k\leq j} H^\circ_k=H^\circ_g + 1.
\]

When $g=-\infty$, by considering $k=\min\{\ell\geq i\dvtx S^\circ
_\ell= J_{i,j}\}$, we see that
\[
\min_{i\leq k\leq j} H^\circ_k=0.
\]

Using (\ref{eq:3}) and the preceding observations, we get that,
for every $k\in\{i,\ldots,j\}$,
%
\begin{equation}\label{eq:9}
H^\circ_k-\min_{i\leq\ell\leq j}H^\circ_\ell=
\#\bigl\{\ell\in\{0,1,\ldots,k-1\}\dvtx\ell> g\mbox{ and }S^\circ_\ell
=J_{\ell,k}\bigr\}.
\end{equation}
Specializing this formula to $k=i$, we have
%
\begin{equation}\label{eq:4}
H^\circ_i-\min_{i\leq\ell\leq j}H^\circ_\ell\leq
\#\bigl\{\ell\in\{g^+,\ldots,i-1\}\dvtx S^\circ_\ell=J_{\ell,i}\bigr\}.
\end{equation}
We now introduce the time-reversed walk $\wh{S}^{(i)}_\ell=S^\circ
_i-S^\circ_{i-\ell}$ for $0\leq
\ell\leq i$. Note that $(\wh{S}^{(i)}_\ell,0\leq\ell\leq i)$
has the same distribution as $(S^\circ_\ell,0\leq\ell\leq i)$. For
every integer $m\geq0$, set
\[
\wh{\rho}^{(i)}_m=\min\bigl\{k\in\{0,\ldots,i\}\dvtx\wh{S}^{(i)}_k\geq m\bigr\},
\]
where $\min\varnothing= +\infty$. For $k\in\{0,1,\ldots,i\}$, we
also set
\[
\wh{\Delta}^{(i)}(k)=\#\Bigl\{\ell\in\{1,\ldots,k\}\dvtx\wh{S}^{(i)}_\ell
=\max_{0\leq
n\leq\ell}\wh{S}^{(i)}_n\Bigr\},
\]
which is the number of (weak)
records of the time-reversed walk $\wh S^{(i)}$ before time~$k$.
Finally, let
$J^{(i)}_{j-i}=J_{i,j}-S^\circ_i$. With these definitions,
(\ref{eq:4}) can be rewritten in the form
%
\begin{equation}\label{eq:8}
H^\circ_i-\min_{i\leq\ell\leq j}H^\circ_\ell\leq
\wh{\Delta}^{(i)}\bigl(\wh{\rho}^{(i)}_{-J^{(i)}_{j-i}}\wedge i\bigr).
\end{equation}
Note that $J^{(i)}_{j-i}$ is independent of the time-reversed walk $\wh
S^{(i)}$ and that,
conditionally on
$\{-J^{(i)}_{j-i}=m\}$, the random variable $\wh{\Delta}^{(i)}(\wh
{\rho}^{(i)}_{-J^{(i)}_{j-i}}\wedge i)$
has the same distribution as $\Delta(\rho_m\wedge i)$, where, for
every integers $k,m\geq0$,
\[
\Delta(k)=\#\Bigl\{\ell\in\{1,\ldots,k\}\dvtx S^\circ_\ell=\max_{0\leq
n\leq\ell}S^\circ_n\Bigr\},\qquad
\rho_m=\inf\{k\geq0\dvtx S^\circ_k\geq m\}.
\]
We thus need to estimate the moments of
$\Delta(\rho_m)$. To this end, introduce the weak record times, which
are defined, by
induction, by $\tau_0=0$ and
\[
\tau_{n+1}=\inf\{k>\tau_n\dvtx S^\circ_k\geq S^\circ_{\tau_n}\} ,\qquad
n\geq0 .
\]
It
is well known (see, e.g., \cite{Trees}, Lemma 1.9) that the random
variables $S^\circ_{\tau_n}-S^\circ_{\tau_{n-1}}$, $n\geq1$, are
i.i.d. with
distribution
\[
\mathbb{P}(S^\circ_{\tau_1}=k)=\ov{\nu}(k),
\]
where $\ov\nu(k)=\nu([k,\infty))=\mu([k+1,\infty))$. From (\ref
{tail2}), we get that
there exists a positive constant $K'_1$ such that, for every $m\geq1$,
\[
\mathbb{P}(S^\circ_{\tau_1}\geq m)\geq K'_1m^{-\alpha+1} .
\]
Consequently, by
arguing as in the proof of Lemma \ref{firstbound}, we get, for every
real $y\geq1$,
\begin{eqnarray*}
\mathbb{P}\bigl(\Delta(\rho_m)> ym^{\alpha-1}\bigr)&\leq&\mathbb{P}(S^\circ_{\tau
_{[ym^{\alpha-1}]}} < m)\leq
P(S^\circ_{\tau_1}< m)^{[ym^{\alpha-1}]}\\&\leq&
\exp(-K'_1y/2) .
\end{eqnarray*}
Thus, the moments of
$\Delta(\rho_m)/m^{\alpha-1}$ are uniformly bounded. From
the remarks following (\ref{eq:8}), we get
\begin{eqnarray*}
\E\bigl[\bigl(\wh{\Delta}^{(i)}\bigl(\wh{\rho}^{(i)}_{-J^{(i)}_{j-i}}\wedge
i\bigr)\bigr)^2\bigr]&\leq&
K'_2\E\bigl[\bigl(-J^{(i)}_{j-i}\bigr)^{2(\alpha-1)}\bigr]=K'_2\E\bigl[(-J_{j-i})^{2(\alpha
-1)}\bigr]\\
&\leq&
K'_3|j-i|^{2(1-1/\alpha)} ,
\end{eqnarray*}
where we have used Lemma \ref{firstbound}
and Jensen's inequality in the
last bound.
By (\ref{eq:8}), this yields
%
\begin{equation}\label{eq:11}
\E\Bigl[\Bigl(H^\circ_i-\min_{i\leq\ell\leq j}H^\circ_\ell\Bigr)^2\Bigr]\leq
K'_3|j-i|^{2(1-1/\alpha)} .
\end{equation}

Next, let us take $k=j$ in (\ref{eq:9}). It follows
that
\[
H^\circ_j-\min_{i\leq\ell\leq j}H^\circ_\ell=
\#\bigl\{\ell\in\{i,\ldots,j-1\}\dvtx S^\circ_\ell=J_{\ell,j}\bigr\}.
\]
Using the same notation as above, we can rewrite the previous displayed quantity
as
\[
\#\Bigl\{\ell\in\{1,\ldots,j-i\}\dvtx\wh{S}^{(j)}_\ell=\max_{0\leq n\leq
\ell}\wh{S}^{(j)}_n\Bigr\}\stackrel{\mathrm{(d)}}{=} \Delta(j-i) .
\]
We claim that, for every integer $p\geq1$, the $p$th moment of
$\Delta(n)/n^{1-1/\alpha}$ is bounded independently of $n\geq1$.
Taking $p=2$, we then deduce, from the previous identity in
distribution, that
\[
\E\Bigl[\Bigl(H^\circ_j-\min_{i\leq\ell\leq j}H^\circ_\ell\Bigr)^2\Bigr]\leq
K'_4|i-j|^{2(1-1/\alpha)} .
\]
The statement of the lemma follows from the last
bound and from (\ref{eq:11}).

It thus remains to verify our claim. We
note that, for every real $y\geq1$ and every $n\geq1$,
\[
\mathbb{P}\bigl(\Delta(n)> yn^{1-1/\alpha}\bigr)\leq\mathbb{P}\bigl(\tau
_{[yn^{1-1/\alpha}]}< n\bigr).
\]
Since
$\tau_n=\sum_{k=1}^{n}(\tau_k-\tau_{k-1})$ and the random variables
$\tau_k-\tau_{k-1}$, $k\geq1$ are i.i.d., the same argument
as in the proof of Lemma \ref{firstbound} shows that
our claim will follow from the bound
%
\begin{equation}
\label{tailtau1}
\mathbb{P}(\tau_1\geq n)\geq K'_5n^{(1/\alpha)-1}
\end{equation}
for some positive constant $K'_5$. From formulas P5(b), page 181
and (3), page 187 of
\cite{spitzer}, IV.17, the generating
function of $\tau_1$ is given by the formula
%
\begin{equation}\label{eq:10}
1-\E[s^{\tau_1}]=\frac{1-s}{1-r_s} ,
\end{equation}
where,
for $0<s<1$, $r_s$ is the unique real solution in $(0,1)$ of equation
$r_s/s=\phi_\mu(r_s)$, with $\phi_\mu(s)=\sum_{k=0}^\infty s^k\mu(k)$.
From a standard Abelian theorem, the asymptotic formula (\ref{tail2})
implies that
$\phi_\mu(s)=s+ K_{(\mu)}(1-s)^{\alpha}+o((1-s)^\alpha)$ as $s\to1$,
with some positive constant $K_{(\mu)}$ depending on $\mu$.
From the equation $r_s/s=\phi_\mu(r_s),$ one then gets that the ratio
$K_{(\mu)}(1-r_s) ^\alpha/(1-s)$ tends to $1$ as $s\to1$.
From this and
(\ref{eq:10}), it follows that
\[
1-\E[s^{\tau_1}]= K_{(\mu)}^{1/\alpha}(1-s)^{1-1/\alpha
}+o\bigl((1-s)^{1-1/\alpha}\bigr)
\]
as $s\to1$.
The desired estimate (\ref{tailtau1}) then follows using
Karamata's Tauberian theorem for power series.
\end{pf*}
\begin{Remark*}
The previous proof may be compared with that of the analogous
statement in the continuous-time setting \cite{duqleg02}, Lemma 1.4.6.
\end{Remark*}
\begin{pf*}{Proof of Lemma \ref{thirdbound}}
To simplify notation, we set
\[
M_n=\max_{0\leq k\leq n} S^\circ_k
\]
for every $n\geq0$. We then have
%
\begin{equation}
\label{formulaG}
G_n=\sum_{k=0}^{n-1} \ind_{\{S^\circ_{k+1}\geq M_k\}} (M_k-S^\circ_k).
\end{equation}
By time reversal, $M_k-S^\circ_k$ has the same distribution as $-J_k$.
We start by deriving some information about the distribution of $J_k$.
From (\ref{tail2}), there exists a constant $K'_6$ such that, for every
$\ell\geq1$,
%
\begin{equation}
\label{3bound00}
\ov{\nu}(\ell)\leq K'_6 \ell^{-\alpha}.
\end{equation}
We use this to verify that, for every $k\geq1$ and $\ell\geq1$,
%
\begin{equation}
\label{3bound1}
\mathbb{P}(J_k>-\ell)\leq K'_7 \frac{\ell}{k^{1/\alpha}}
\end{equation}
with some constant $K'_7$.
Clearly, we may assume that $\ell<k^{1/\alpha}/10$. Recall the
notation $V_k$
introduced in the proof of Lemma \ref{firstbound}. As we already noted
in the
proof of this lemma,
$k^{-1}V_{[k^{1/\alpha}]}$ converges in distribution to a
stable variable with index $1/\alpha$ as $k\to\infty$. This implies
that there exists
a constant $c_*$ such that, for every $k\geq1$,
\[
\mathbb{P}\bigl(V_{[k^{1/\alpha}]}>k\bigr)\leq c_*<1.
\]
Let $U_1,U_2,\ldots$ be
independent random variables distributed as $V_\ell$. Then,
\begin{eqnarray*}
\mathbb{P}\bigl(V_{[k^{1/\alpha}]}>k\bigr)&\geq&\mathbb{P}\bigl(U_1+U_2+\cdots+U_{[\ell
^{-1}[k^{1/\alpha}]]}>k\bigr)\\ &\geq&1-\mathbb{P}(U_i\leq
k, \forall
i=1,\ldots,[\ell^{-1}[k^{1/\alpha}]] )\\ &=&1-\bigl(1-\mathbb{P}(V_\ell
>k)\bigr)^{[\ell^{-1}[k^{1/\alpha}]]}.
\end{eqnarray*}
Combining the last two displays, we get
\[
\bigl(1-\mathbb{P}(V_\ell>k)\bigr)^{[\ell^{-1}[k^{1/\alpha}]]}\geq1- c_*
\]
and, consequently,
\[
\mathbb{P}(V_\ell>k)\leq1- (1-c_*)^{1/[\ell^{-1}[k^{1/\alpha}]]}.
\]
The bound
(\ref{3bound1}) follows since $\mathbb{P}(J_k>-\ell)=\mathbb{P}(V_\ell
>k)$. Using
the bound
(\ref{3bound1}), we easily get that there exists a constant $K'_8$ such
that, for every $k\geq1$,
%
\begin{equation}
\label{3bound3}
\E[|J_k|^{1-\alpha}\wedge1]\leq K'_8 k^{(1/\alpha)-1}.
\end{equation}

Let us now bound $\E[(G_n)^2]$. From (\ref{formulaG}), we have
\begin{eqnarray*}
G_n &=& \sum_{k=0}^{n-1} \ov{\nu}(M_k-S^\circ_k) (M_k-S^\circ_k) +
\sum_{k=0}^{n-1} \bigl(\ind_{\{S^\circ_{k+1}\geq
M_k\}}-\ov{\nu}(M_k-S^\circ_k)\bigr)(M_k-S^\circ_k)\\
&=&\!:G'_n+G''_n.
\end{eqnarray*}
We first bound $\E[(G''_n)^2]$. Using the Markov
property for the random walk $S^\circ$ and, more precisely, the fact that
$\mathbb{P}(S^\circ_{k+1}\geq M_k\mid S^\circ_0,\ldots, S^\circ_k)=\ov
{\nu}(M_k-S^\circ_k)$, we
get
\begin{eqnarray*}
\E[(G''_n)^2]&=&\E\Biggl[\sum_{k=1}^{n-1} \bigl(\ind_{\{S^\circ_{k+1}\geq
M_k\}}-\ov{\nu}(M_k-S^\circ_k)\bigr)^2
(M_k-S^\circ_k)^2 \Biggr]\\ &=&\E\Biggl[\sum_{k=1}^{n-1}
(M_k-S^\circ_k)^2\ov{\nu}(M_k-S^\circ_k) \bigl(1-\ov{\nu}(M_k-S^\circ
_k)\bigr) \Biggr]\\ &\leq&\E\Biggl[\sum_{k=1}^{n-1}
(M_k-S^\circ_k)^2\ov{\nu}(M_k-S^\circ_k) \Biggr].
\end{eqnarray*}
Using the estimate (\ref{3bound00}), the fact that $M_k-S^\circ_k$
has the same
distribution as $|J_k|$ and then Lemma \ref{firstbound}
together with Jensen's inequality, we get
\[
\E[(G''_n)^2] \leq K'_6 \sum_{k=1}^{n-1}
\E[|J_k|^{2-\alpha}] \leq K'_6(K_2)^{(2-\alpha)/2} \sum_{k=1}^{n-1}
k^{2/\alpha-1} \leq K'_9 n^{2/\alpha}.
\]

We then turn to $E[(G'_n)^2]$. We have
\begin{eqnarray*}
\E[(G'_n)^2] &=& \E\Biggl[\sum_{k=0}^{n-1} \ov{\nu}(M_k-S^\circ_k)^2
(M_k-S^\circ_k)^2 \Biggr]\\
&&{} +2 \E\biggl[\sum_{0\leq k<j\leq n-1}
\ov{\nu}(M_k-S^\circ_k) (M_k-S^\circ_k) \ov{\nu}(M_j-S^\circ_j)
(M_j-S^\circ_j) \biggr].
\end{eqnarray*}
Since $\ov\nu(M_k-S^\circ_k)\leq1$, the first term in the
right-hand side is bounded above by
$K'_9n^{2/\alpha}$, as in the preceding calculation. Using (\ref
{3bound00}), the second term is
bounded above by
\[
2(K'_6)^2 \E\biggl[\sum_{0\leq k<j\leq n-1} \bigl((M_k-S^\circ_k)^{1-\alpha
}\wedge
1\bigr) \bigl((M_j-S^\circ_j)^{1-\alpha}\wedge1\bigr) \biggr].
\]
To bound this quantity,
we note that, for fixed $k$ and $j$ such that $0\leq k<j$, the
distribution of $M_j-S^\circ_j$, given the past of $S^\circ$ up to
time $k$,
dominates the (unconditional) distribution of $M_{j-k}-S^\circ_{j-k}$. Since
the function $x\to x^{1-\alpha}\wedge1$ is nonincreasing over $\R_+$,
it follows that the quantity in the last display is bounded above by
\begin{eqnarray*}
&&2(K'_6)^2 \sum_{0\leq k<j\leq n-1}\E[(M_k-S^\circ_k)^{1-\alpha
}\wedge1]
\E[(M_{j-k}-S^\circ_{j-k})^{1-\alpha}\wedge1]\\
&&\qquad \leq
2(K'_6)^2 \Biggl(\sum_{k=0}^{n-1} \E[(M_k-S^\circ_k)^{1-\alpha}\wedge
1] \Biggr)^2\\
&&\qquad = 2(K'_6)^2 \Biggl(\sum_{k=0}^{n-1}
\E[|J_k|^{1-\alpha}\wedge1] \Biggr)^2\\
&&\qquad \leq
2(K'_6)^2(K'_8)^2 \Biggl(1+\sum_{k=1}^{n-1}
k^{(1/\alpha)-1} \Biggr)^2\\
&&\qquad \leq K'_{10} n^{2/\alpha}.
\end{eqnarray*}
In the penultimate line of the
calculation, we have used the bound (\ref{3bound3}). We conclude that
$\E[(G'_n)^2]\leq(K'_9+K'_{10})n^{2/\alpha}$, which completes the proof
of Lemma \ref{thirdbound}.
\end{pf*}

\section{Contour processes and conditioned trees}
\label{contour-conditioned}

\subsection{Contour processes}

In view of our applications to random planar maps, it will be important
to reformulate Theorem \ref{functional} in terms of
contour processes associated with our
forest of mobiles. We consider the same general setting
as in the previous section. In particular, $u_0,u_1,\ldots$ are
the white vertices of the forest~$\F$, listed one tree after
another and in lexicographical order for every tree. Recall
that $H^\circ_n=\frac{1}{2} |u_n|$. We also denote by
$x_0,x_1,\ldots$ the sequence obtained by
concatenating the white contour sequences of
$\theta_1,\theta_2,\ldots.$ Notice that some of
the vertices $u_0,u_1,\ldots$ appear more than once in
the sequence $x_0,x_1,\ldots.$ More precisely, the number
of occurrences of a given white vertex of $\F$
is equal to $1$ plus the number of its black children. We set
$C^\circ_n=\frac{1}{2} |x_n|$ and denote by
$\Lambda_n$ the label of $x_n$.

In order to study the scaling limit of $(C^\circ_n)_{n\geq0}$, we define,
for every $n\geq0$,
\[
R_n=\inf\{j\geq0\dvtx x_j=u_n\}.
\]
Clearly,
\[
C^\circ_{R_n}= \tfrac{1}{2} |x_{R_n}|= \tfrac{1}{2} |u_n| = H^\circ_n.
\]
\begin{lmm}
\label{height-contour}
We have
\[
\lim_{n\to\infty} \frac{R_n}{n}= \frac{1}{\beta} \qquad\mbox{a.s.}
\]
\end{lmm}
\begin{pf}
For every $j=0,1,\ldots,$ let $B(j)$ denote the number of black
children of $u_j$. Notice that the random variables $B(0),B(1),\ldots$
are independent and distributed according to $\mu_0$. We first observe that
%
\begin{equation}
\label{upper-cont}
R_n\leq\sum_{j=0}^{n-1} \bigl(B(j)+1\bigr).
\end{equation}
This bound comes from the fact that any vertex that is visited
by the contour sequence $x_0,x_1,\ldots$ before the first visit of $u_n$
must be smaller than $u_n$ in lexicographical order. Hence, $R_n$
has to be smaller than the total number of visits by the contour
sequence of all vertices that are smaller than $u_n$ in
lexicographical order. The bound (\ref{upper-cont}) follows.

Since the mean of $\mu_0$ is $m_0=Z_qf_q(Z_q)= \frac{1}{\beta}-1$, the
law of large numbers gives
\[
\limsup_{n\to\infty} \frac{R_n}{n}\leq\frac{1}{\beta} \qquad\mbox{a.s.}
\]
We would like to derive the reverse inequality. To this end, note that
if a vertex $u_j$ with $j<n$ is not an ancestor of $u_n$, then all of
its visits
by the contour sequence will occur before the first visit of $u_n$. Thus,
\[
R_n\geq n + \sum_{j=0}^{n-1} B(j) \ind{\{u_j\mbox{ is not an
ancestor of }u_n\}}
\]
or, equivalently,
%
\begin{eqnarray}
\label{lower-count}
\sum_{j=0}^{n-1} \bigl(B(j)+1\bigr) - R_n
&\leq&\sum_{j=0}^{n-1} B(j) \ind{\{u_j\mbox{ is an ancestor of }u_n\}
}\nonumber\\[-8pt]\\[-8pt]
&\leq& H^\circ_n\times\sup_{0\leq j\leq n-1} B(j).\nonumber
\end{eqnarray}
A crude estimate gives, for every $\eps>0$,
\[
\lim_{n\to\infty} \frac{1}{n^\eps} \sup_{0\leq j\leq n-1} B(j)=0\qquad
\mbox{a.s.}
\]
On the other hand, by a special case of Lemma \ref{secondbound}, we
know that
$E[(H^\circ_n)^2]\leq K_3 n^{2(1-1/\alpha)}$. Using the
Markov inequality and then the Borel--Cantelli lem\-ma, we can find
$\eps>0$ such that
%
\begin{equation}
\label{height-contour-tech1}
\lim_{n\to\infty} \frac{1}{n^{1-\eps}} H^\circ_n=0 \qquad\mbox{a.s.}
\end{equation}
and we conclude that
\[
\lim_{n\to\infty} \frac{1}{n} H^\circ_n\times\sup_{0\leq j\leq
n-1} B(j)=0 \qquad\mbox{a.s.}
\]
The desired result then follows from (\ref{lower-count})
and the law of large numbers.
\end{pf}
\begin{Remark*}
Since the sequence $(R_n)_{n\geq0}$ is monotone
increasing, we also have,
for every $A>0$,
%
\begin{equation}
\label{height-contour-uni}
\lim_{n\to\infty} \frac{1}{n} \sup_{0\leq k \leq An} \biggl| R_k - \frac
{k}{\beta} \biggr| =0 \qquad\mbox{a.s.}
\end{equation}
\end{Remark*}

The next proposition is an analog of Theorem \ref{functional} for
contour processes.
\begin{prp}
\label{limit-contour}
We have
\[
\bigl(n^{-(1-1/\alpha)} C^{\circ}_{[nt]},n^{-1/2\alpha} \Lambda_{[nt]}
\bigr)_{t\geq0}
\mathop{\longrightarrow}_{n\to\infty}^{(d)} \bigl( c_0^{-1} H_{\beta t},
\sqrt{2c_0} D_{\beta t} \bigr)_{t\geq0},
\]
where the convergence holds in the sense of weak convergence of the
laws in
the Skorokhod space $\D(\R^2)$.
\end{prp}
\begin{pf}
Fix an integer $A>0$. The statement of the proposition will
be an immediate consequence of Theorem \ref{functional} once we have
verified that
%
\begin{equation}
\label{limit-cont1}
n^{-(1-1/\alpha)} \sup_{0\leq k\leq A n} \bigl|C^{\circ}_k - H^\circ
_{[\beta k]}\bigr|
\mathop{\longrightarrow}_{n\to\infty} 0 \qquad\mbox{in probability}
\end{equation}
and
%
\begin{equation}
\label{limit-cont2}
n^{-1/2\alpha} \sup_{0\leq k\leq A n} \bigl|\Lambda_k - L^\circ_{[\beta k]}\bigr|
\mathop{\longrightarrow}_{n\to\infty} 0 \qquad\mbox{in probability.}
\end{equation}

Let us start with the proof of (\ref{limit-cont1}).
It is elementary to check that for every integer $n\geq0$,
%
\begin{equation}
\label{limit-cont11}
{\sup_{R_n\leq j\leq R_{n+1}}} |C^{\circ}_j - C^{\circ}_{R_n}| \leq
|H^\circ_{n+1}- H^\circ_n| +1.
\end{equation}
Then, note that if
$k\in\{0,1,\ldots,An\}$ and $\ell$ is chosen so that $R_\ell\leq k
<R_{\ell+1}$, we have
\[
\bigl|C^\circ_k - H^\circ_{[\beta k]}\bigr| \leq
|C^{\circ}_k - C^{\circ}_{R_\ell}| + \bigl|H^\circ_\ell- H^\circ
_{[\beta k]}\bigr|
\]
since $C^{\circ}_{R_\ell}= H^\circ_\ell$.
By (\ref{limit-cont11}) and the fact that the limiting process $H$ in
(\ref{heightconv}) is continuous, we have
%
\begin{equation}
\label{limit-cont3}
{n^{-(1-1/\alpha)} \sup_{0\leq\ell\leq A n} \sup_{R_\ell\leq k <
R_{\ell+1}}} |C^{\circ}_k - C^\circ_{R_\ell}|
\mathop{\longrightarrow}_{n\to\infty} 0 \qquad\mbox{in probability.}
\end{equation}
On the other hand, for every fixed $\eps>0$, it follows from (\ref
{height-contour-uni}) that, with a probability close
to $1$ when $n$ is large, we have, for every $\ell=0,1,\ldots,An$,
\[
\ell-\eps n \leq\beta R_\ell\leq\beta R_{\ell+1} \leq\ell+ \eps n
\]
and thus
\begin{eqnarray*}
&&
n^{-(1-1/\alpha)} \sup_{0\leq\ell\leq A n} \sup_{R_\ell\leq k <
R_{\ell+1}} \bigl|H^\circ_\ell- H^\circ_{[\beta k]}\bigr|\\
&&\qquad \leq n^{-(1-1/\alpha)} \sup_{r,s\in[0,A+\eps], |r-s|\leq\eps}
\bigl|H^\circ_{[nr]}- H^\circ_{[ns]}\bigr|.
\end{eqnarray*}
The right-hand side will be small
in probability when $n$ is large, again by (\ref{heightconv}),
provided that $\eps$ has been chosen small enough. This completes the
proof of (\ref{limit-cont1}).

Let us now prove (\ref{limit-cont2}).
Notice that $L^\circ_n=\Lambda_{R_n}$ for every $n\geq0$. We can
therefore argue in a way similar
to the proof of (\ref{limit-cont1}), using Theorem \ref{functional}
in place of (\ref{heightconv}), provided that
we establish the analog of (\ref{limit-cont3}),
%
\begin{equation}
\label{limit-lab1}
{n^{-1/2\alpha} \sup_{0\leq\ell\leq A n} \sup_{R_\ell\leq k <
R_{\ell+1}}} |\Lambda_k - \Lambda_{R_\ell}|
\mathop{\longrightarrow}_{n\to\infty}0 \qquad\mbox{in probability.}
\end{equation}
So, let us verify that (\ref{limit-lab1}) holds. From the distribution
of labels, it is easy to
check that, for every fixed $n\geq0$, conditionally on the forest $\F$,
the sequence
\[
\bigl(\Lambda_{(R_n+j)\wedge R_{n+1}} - \Lambda_{R_n}\bigr)_{j\geq0}
\]
is a martingale (in fact, the increments of this sequence are both
independent and
centered, conditionally given $\F$). By Doob's inequality, there
are constants $K$ and $K'$ such that, for every $\ell\geq0$,
\[
\E\Bigl[ \sup_{R_\ell\leq k < R_{\ell+1}} (\Lambda_k - \Lambda_{R_\ell
})^4 \big| \F\Bigr]
\leq K \E[(\Lambda_{R_{\ell+1}} - \Lambda_{R_\ell})^4 | \F]
\]
and
\[
\E\Bigl[ \sup_{R_\ell\leq k < R_{\ell+1}} (\Lambda_k - \Lambda_{R_\ell
})^4 \Bigr]
\leq K \E[(\Lambda_{R_{\ell+1}} - \Lambda_{R_\ell})^4 ] \leq K',
\]
using\vspace*{1pt} Proposition \ref{keybound} with $i=\ell$ and $j=\ell+1$. Finally,
if $\eps>0$ is small enough so that $\frac{2}{\alpha} -4\eps-1 >0$,
we have
\[
\mathbb{P}\Bigl[ {\sup_{0\leq\ell\leq A n} \sup_{R_\ell\leq k < R_{\ell+1}}}
|\Lambda_k - \Lambda_{R_\ell}|
\geq n^{(1/2\alpha)-\eps} \Bigr] \leq(An+1) K'\bigl(n^{(1/2\alpha)-\eps}\bigr)^{-4},
\]
which tends to $0$ as $n\to\infty$.
This completes the proof of (\ref{limit-lab1}) and of the proposition.
\end{pf}

\subsection{Conditioning a mobile to have more than $n$ white vertices}

The definition of the continuous-time height process $(H_t)_{t\geq0}$
also makes sense under the excursion measure $\bN$, or under $\bN
(\cdot\mid\sigma= 1)$
(see Chapter 1 of \cite{duqleg02}).
Furthermore,
the law of the pair $(H_t,D_t)_{t\geq0}$ under $\bN(\cdot\mid\sigma
>1)$ coincides with
the law of $(H_{(g_{(1)}+t)\wedge d{(1)}},D_{(g_{(1)}+t)\wedge
d{(1)}})_{t\geq0}$ under $\mathbb{P}$,
where $(g_{(1)},d_{(1)})$ is the first excursion interval of $X-I$ with
length greater than $1$.
This follows from a minor extension of the arguments of Section \ref{sec43}.

For every integer $n\geq1$, we set $\wt\Q^{(n)}=\Q(\cdot\mid\#\mathcal{T}
^\circ\geq n)$.
\begin{theorem}
\label{heightcondi1}
The law of $\frac{1}{n}\#\mathcal{T}^\circ$ under $\wt\Q^{(n)}$ converges,
as $n\to\infty,$ to the law
of $\sigma$ under $\mathbf{N}(\cdot\mid\sigma>1)$. Moreover,
the law of the process
\[
\bigl(n^{-(1-1/\alpha)} H^\theta_{[nt]}, n^{-1/2\alpha} L^{\theta
}_{[nt]} \bigr)_{t\geq0}
\]
under $\wt\Q^{(n)}(d\theta)$
converges, as $n\to\infty$, to the law of the process
\[
\bigl(c_0^{-1} H_t, \sqrt{2c_0} D_t \bigr)_{t\geq0}
\]
under $\mathbf{N}(\cdot\mid\sigma>1)$.
Similarly, the law of the process
\[
\bigl(n^{-(1-1/\alpha)} C^{\theta}_{[nt]}, n^{-1/2\alpha} \Lambda
^{\theta}_{[nt]} \bigr)_{t\geq0}
\]
under $\wt\Q^{(n)}(d\theta)$
converges, as $n\to\infty,$ to the law of
\[
\bigl(c_0^{-1} H_{\beta t}, \sqrt{2c_0} D_{\beta t} \bigr)_{t\geq0}
\]
under $\mathbf{N}(\cdot\mid\sigma> 1)$.
\end{theorem}
\begin{pf}
Thanks to Theorem 1 and the Skorokhod representation theorem, we can construct,
for every integer $n\geq1$, a random labeled forest $\mathbf{F}^{(n)}$ having
the same distribution as $\mathbf{F}$, in such a way that
%
\begin{eqnarray}
\label{heightcondi11}
&&
\bigl(n^{-1/\alpha}S^{(n)}_{[nt]},n^{-(1-1/\alpha)} H^{(n)}_{[nt]},
n^{-1/2\alpha} L^{(n)}_{[nt]} \bigr)_{t\geq0}\nonumber\\[-8pt]\\[-8pt]
&&\qquad
\mathop{\longrightarrow}_{n\to\infty}^{\mathrm{a.s.}}
\bigl(c_0X_t, c_0^{-1} H_t, \sqrt{2c_0} D_t\bigr)_{t\geq0},\nonumber
\end{eqnarray}
where we have used the notation of the proof of Proposition \ref{finitedim}.
Let $\wt\theta^{(n)}$ be the first mobile in the forest
$\mathbf{F}^{(n)}$ with at least $n$ white vertices and note that
$\wt\theta^{(n)}$ is distributed according to $\wt\Q^{(n)}$. Let $[g_n,d_n]$
be the first excursion interval of $H^{(n)}$ away from $0$
with length greater than or equal to $n$. Then, writing $\wt H^{(n)}$
and $\wt L^{(n)}$ for the height process and the label process of $\wt
\theta^{(n)}$,
respectively, we have, for every $k\geq0$,
\[
\wt H^{(n)}_k = H^{(n)}_{(g_n+k)\wedge d_n} ,\qquad
\wt L^{(n)}_k = L^{(n)}_{(g_n+k)\wedge d_n} .
\]
This is\vspace*{1pt} the case because the interval $[g_n,d_n)$ corresponds exactly to
those integers $j$ such that the $(j+1)$st vertex of $\mathbf{F}^{(n)}$
(in lexicographical order) belongs
to $\wt\theta^{(n)}$.

One can then deduce from (\ref{heightcondi11}) that
%
\begin{equation}
\label{heightcondi12}
\frac{1}{n} g_n\mathop{\longrightarrow}_{n\to\infty}^{\mathrm{a.s.}}
g_{(1)} ,\qquad
\frac{1}{n} d_n\mathop{\longrightarrow}_{n\to\infty}^{\mathrm{a.s.}}
d_{(1)} .
\end{equation}
We omit
the details of the derivation of (\ref{heightcondi12}); see the proof of
Proposition 2.5.2 in \cite{duqleg02} or the proof of Corollary 1.13 in
\cite{Trees} for a very similar
argument.

The first assertion of the theorem readily follows from (\ref{heightcondi12})
since the number of white vertices of $\wt\theta^{(n)}$ is $d_n-g_n$
and the law of $d_{(1)}-g_{(1)}$
is precisely the law of $\sigma$ under $\mathbf{N}(\cdot\mid\sigma>1)$.

We then have
\begin{eqnarray*}
&&\bigl(n^{-(1-1/\alpha)} \wt H^{(n)}_{[nt]}, n^{-1/2\alpha} \wt
L^{(n)}_{[nt]}\bigr)\\
&&\qquad = \bigl(n^{-(1-1/\alpha)} H^{(n)}_{[n(({g_n}/{n}+t)\wedge
{d_n}/{n})]}, n^{-1/2\alpha} L^{(n)}
_{[n(({g_n}/{n}+t)\wedge{d_n}/{n})]}\bigr)
\end{eqnarray*}
and thus (\ref{heightcondi11}) and (\ref{heightcondi12}) give
\begin{eqnarray*}
&&\bigl(n^{-(1-1/\alpha)} \wt H^{(n)}_{[nt]}, n^{-1/2\alpha} \wt
L^{(n)}_{[nt]}\bigr)_{t\geq0}\\
&&\qquad \mathop{\longrightarrow}_{n\to\infty}^{\mathrm{a.s.}}
\bigl(c_0^{-1} H_{(g_{(1)}+t)\wedge d_{(1)}}, \sqrt
{2c_0}D_{(g_{(1)}+t)\wedge d_{(1)}}\bigr)_{t\geq0}.
\end{eqnarray*}
The first convergence stated in the theorem follows since we know that
the limiting process has the desired distribution.

Let us turn to the proof of the second convergence of the theorem.
From (\ref{limit-cont1}) and (\ref{limit-cont2}), we know that,
for every integer $A>0$,
\[
n^{-(1-1/\alpha)} \sup_{k\leq An} \bigl|C^{(n)}_k - H^{(n)}_{[\beta k]} \bigr|
\mathop{\longrightarrow}_{n\to\infty} 0 \qquad\mbox{in probability}
\]
and
\[
n^{-1/2\alpha} \sup_{k\leq An} \bigl| \Lambda^{(n)}_k - L^{(n)}_{[\beta
k]} \bigr|
\mathop{\longrightarrow}_{n\to\infty} 0 \qquad\mbox{in probability.}
\]
Write $\wt C^{(n)}$ and $\wt\Lambda^{(n)}$ for the contour process and
the contour label process, respectively, of $\wt\theta^{(n)}$. We
have for every $t\geq0$,
\[
\wt C^{(n)}_{[nt]} =
C^{(n)}_{(R_{g_n}+[nt])\wedge R_{d_n}}.
\]
Writing
\[
(R_{g_n}+[nt])\wedge R_{d_n}=n \biggl(\biggl(\frac{R_{g_n}}{n} + \frac
{[nt]}{n}\biggr)\wedge
\frac{R_{d_n}}{n} \biggr)
\]
and using Lemma \ref{height-contour} together with
(\ref{heightcondi12}), we get
\[
n^{-(1-1/\alpha)} \sup_{t\geq0}
\bigl| \wt C^{(n)}_{[nt]} - H^{(n)}_{[n((g_{(1)}+\beta t)\wedge d_{(1)})]} \bigr|
\mathop{\longrightarrow}_{n\to\infty} 0 \qquad\mbox{in probability.}
\]
Similarly, we have
\[
n^{-1/2\alpha} \sup_{t\geq0}
\bigl| \wt\Lambda^{(n)}_{[nt]} - L^{(n)}_{[n((g_{(1)}+\beta t)\wedge
d_{(1)})]} \bigr|
\mathop{\longrightarrow}_{n\to\infty} 0 \qquad\mbox{in probability.}
\]
The desired result now follows from (\ref{heightcondi11}).
\end{pf}

\subsection{Conditioning a mobile to have exactly $n$ white vertices}

We now set $\ov\Q{}^{(n)}=\Q(\cdot\mid\#\mathcal{T}^\circ= n)$. Note that
this makes sense
(the conditioning event has positive probability) for all sufficiently
large $n$. From now on, we consider
only such values of $n$. Our goal is to derive an analog of
Theorem \ref{heightcondi1} when $\wt\Q^{(n)}$ is replaced by $\ov\Q
{}^{(n)}$. The proof
is more delicate and will require a few preliminary lemmas.

Let $\theta=(\mathcal{T},(\ell(v))_{v\in\tc})$ be a mobile.
Recall that $w_0(\theta),\ldots,w_{\#\mathcal{T}^\circ-1}(\theta)$ are the
white vertices of
$\theta$ listed in lexicographical order. By convention, we put
$w_l(\theta)=\varnothing$
when $l\geq\#\mathcal{T}^\circ$.
For every $k\geq1$, we then define
another mobile $\theta^{[k]}=(\mathcal{T}_{[k]},(\ell_{[k]}(v))_{v\in
\mathcal{T}
_{[k]}^\circ})$ in the
following way. First,
$\mathcal{T}_{[k]}$ consists of the vertices $w_0(\theta),\ldots
,w_{k-1}(\theta)$,
together with all of the (black) children and all of the (white)
grandchildren of these
vertices in $\mathcal{T}$.
Then, $\ell_{[k]}(v)=\ell(v)$ for every $v\in\mathcal{T}_{[k]}^\circ$.
By convention, we also define $\theta^{[0]}$ as the trivial mobile
with just one vertex.

For every $k\geq0$, we let
$\g_k$ be the $\sigma$-field on $\Theta$ generated by the mapping
\mbox{$\theta\to\theta^{[k]}$}. It is
easily checked that the processes $H^\theta_k$ and $L^\theta_k$ are
adapted to the filtration $(\g_k)_{k\geq0}$.

Recall that, by definition of the Lukasiewicz path $S^\theta$, for
$j\in\{1,\ldots,\#\mathcal{T}^\circ\}$, $S^\theta_j-S^\theta_{j-1}+1$ is
the number of (white) grandchildren
of $w_{j-1}(\theta)$. It follows that, for every $k\geq0$, $S^\theta_k$
is $\g_k$-measurable. Furthermore, under the probability measure $\Q$,
the process $(S^\theta_k)_{k\geq0}$ is Markovian with respect to the
filtration
$(\g_k)_{k\geq0}$ and its transition kernels are those of the random
walk with jump distribution $\nu$ stopped at its first hitting time of $-1$.
The preceding properties can be derived by a minor modification of the
arguments found in Section 1 of \cite{Trees}. We leave the details to
the reader.

Recall our notation $(S_k)_{k\geq0}$ for a random walk with jump distribution
$\nu$. We assume that $S_0=j$ under the probability measure $\mathbb{P}_j$
for every $j\in\Z$. We set $V=\inf\{k\geq0\dvtx S_k=-1\}$.
\begin{lmm}
\label{absocont}
Let $k\in\{1,2,\ldots,n-1\}$. The Radon--Nikodym derivative of $\ov
\Q{}^{(n)}$ with respect to $\wt\Q^{(n)}$ on the $\sigma$-field $\g_k$
is equal to $\Gamma(k,n,S^\theta_k)$, where, for every integer $j\geq
0$,
\[
\Gamma(k,n,j)=\frac{\psi_{n-k}(j)/\psi_n(0)}{\varphi
_{n-k}(j)/\varphi_n(0)}
\]
and, for every integer $p\geq0$,
\begin{eqnarray*}
\psi_p(j)&=&\mathbb{P}_j(V=p),\\
\varphi_p(j)&=&\mathbb{P}_j(V\geq p).
\end{eqnarray*}
\end{lmm}
\begin{Remark*}
If $k\leq\#\tc$, then the number of white vertices of
$\theta^{[k]}$ is
$k+1+S^\theta_k$. If $\gamma$ has (strictly) more than $n$
white vertices, then $\ov\Q{}^{(n)}(\theta^{[k]}=\gamma)=0$. This is
consistent
with the fact that $\psi_{n-k}(j)=0$ if $j>n-k-1$.
\end{Remark*}
\begin{pf*}{Proof of Lemma \ref{absocont}}
Let $\gamma$ be a mobile with strictly more than $k$ white vertices
and such that $\gamma^{[k]}=\gamma$ (these are the necessary and
sufficient conditions for $\gamma$ to be of the form $\theta^{[k]}$
for some $\theta\in\Theta$ with at least $n$ white vertices). Then,
\[
\ov\Q{}^{(n)}\bigl(\theta^{[k]}=\gamma\bigr)= \frac{\Q(\{\theta^{[k]}=\gamma
\}\cap
\{\#\mathcal{T}^\circ= n\})}
{\Q(\#\mathcal{T}^\circ= n)}.
\]
On one hand,
\[
\Q(\#\mathcal{T}^\circ= n)=\mathbb{P}_0(V=n)=\psi_n(0).
\]
On the other hand, by the remarks preceding the statement of the lemma,
\begin{eqnarray*}
\Q\bigl(\bigl\{\theta^{[k]}=\gamma\bigr\}\cap\{\#\mathcal{T}^\circ= n\}\bigr)
&=& \Q\bigl(\bigl\{\theta^{[k]}=\gamma\bigr\}\cap\bigl\{\inf\{p\geq0\dvtx S^\theta_p=-1\}
=n\bigr\}\bigr)\\
&=&\Q\bigl(\ind_{\{\theta^{[k]}=\gamma\}} \mathbb{P}_{S^\theta_k}(V=n-k)\bigr)\\
&=&\Q\bigl(\ind_{\{\theta^{[k]}=\gamma\}} \psi_{n-k}(S^\theta_k)\bigr).
\end{eqnarray*}
We thus have
\[
\ov\Q{}^{(n)}\bigl(\theta^{[k]}=\gamma\bigr)= \Q\biggl( \ind_{\{\theta
^{[k]}=\gamma\}}
\frac{\psi_{n-k}(S^\theta_k)}{\psi_n(0)} \biggr).
\]
Similar arguments give
\[
\wt\Q^{(n)}\bigl(\theta^{[k]}=\gamma\bigr)= \Q\biggl( \ind_{\{\theta
^{[k]}=\gamma\}}
\frac{\varphi_{n-k}(S^\theta_k)}{\varphi_n(0)} \biggr).
\]
The desired result follows.
\end{pf*}
\begin{lmm}
\label{boundeddensity}
Let $a\in(0,1)$. There exist an integer $n_0$ and a constant $K$ such that,
for every $n\geq n_0$ and every $j\geq0$,
\[
\Gamma([an],n,j)\leq K.
\]
\end{lmm}
\begin{pf}
By Kemperman's formula (see, e.g., \cite{pitman06}, page 122), for
every $j\geq0$
and $n\geq1$,
%
\begin{equation}
\label{Kemp}
\mathbb{P}_j(V=n)=\frac{j+1}{n} \mathbb{P}_0(S_n=-j-1).
\end{equation}
On the other hand, Gnedenko's local limit theorem (see \cite{Ibragimov}, Theorem
4.2.1) shows that
%
\begin{equation}
\label{locallimit}
\lim_{n\to\infty} \sup_{k\in\Z}
\biggl| n^{1/\alpha} \mathbb{P}_0(S_n=k)-g\biggl(\frac{k}{n^{1/\alpha}}\biggr) \biggr| =0,
\end{equation}
where the function $g$ is continuous and (strictly) positive over $\R
$. Taking $k=-1$, we get that
there exist positive constants $K_1$ and $K_2$ such that, for $n$ large,
\[
\psi_n(0)=\frac{1}{n} \mathbb{P}_0(S_n=-1)\geq K_1 n^{-1-1/\alpha}
\]
and
\[
\varphi_n(0)=\sum_{m=n}^\infty\frac{1}{m} \mathbb{P}_0(S_m=-1) \leq K_2
n^{-1/\alpha}
\]
[the latter bound can also be derived from (\ref{equiv-hitting})].

So, in order to get the desired statement, we need to verify that the
quantity
\[
\frac{n\psi_{n-[an]}(j)}{\varphi_{n-[an]}(j)}
\]
is bounded when $n$ is large, uniformly in $j$.

First, consider the
case when $j\leq n^{1/\alpha}$. From (\ref{Kemp}) and (\ref
{locallimit}), we obtain that there
exist positive constants $K_3$ and $K_4$ such that, for $n$ large,
\[
\psi_{n-[an]}(j)=\frac{j+1}{n} \mathbb{P}_0\bigl(S_{n-[an]}=-j-1\bigr)\leq
K_3(j+1)n^{-1-1/\alpha}
\]
and
\begin{eqnarray*}
\varphi_{n-[an]}(j)&=& (j+1)\sum_{m=n-[an]}^\infty\frac{1}{m}
\mathbb{P}_0(S_{m}=-j-1)\\
&\geq& K_4(j+1)n^{-1/\alpha}.
\end{eqnarray*}
The desired bound follows.

Suppose, then, that $j\geq n^{1/\alpha}$. It easily follows from (\ref
{Levyconv})
that there exists a positive constant $K_5$ such that
\[
\varphi_{n-[an]}(j)\geq K_5 >0.
\]
On the other hand, we have already noted that the law of $V$ under
$\mathbb{P}
_0$ is in the domain
of attraction of a stable distribution with index $1/\alpha$. Another
application of
Gnedenko's local limit theorem shows that
\[
\lim_{k\to\infty} \sup_{n\geq1} \biggl| k^\alpha\mathbb{P}_k(V=n) - \wt
g\biggl(\frac{n}{k^\alpha}\biggr) \biggr| = 0 ,
\]
where the function $g$ is continuous and bounded over $(0,\infty)$.
Hence, there exists a constant $K_6$ such that, for all integers $n\geq
1$ and $k\geq n^{1/\alpha}$,
%
\begin{equation}
\label{dist-V}
n \mathbb{P}_k(V=n)\leq k^\alpha\mathbb{P}_k(V=n)\leq K_6.
\end{equation}
It immediately follows that
\[
n \psi_{n-[an]}(j)=\frac{n}{n-[an]} (n-[an]) \mathbb{P}_j(V=n-[an])\leq
\frac{K_6}{1-a},
\]
giving the desired bound when $j\geq n^{1/\alpha}$. This completes the proof.
\end{pf}
\begin{prp}
\label{condiduration1}
The law of the process
\[
\bigl(n^{-1/\alpha}S^\theta_{[nt]},n^{-(1-1/\alpha)} H^\theta_{[nt]}
\bigr)_{t\geq0}
\]
under $\ov\Q{}^{(n)}(d\theta)$
converges, as $n\to\infty,$ to the law of the process
\[
(c_0X_t,c_0^{-1} H_t )_{t\geq0}
\]
under $\mathbf{N}(\cdot\mid\sigma=1)$.
\end{prp}

This follows from Theorem 3.1 in \cite{duq03}. This theorem gives the
convergence in distribution of the rescaled height process
$(n^{-(1-1/\alpha)} H^\theta_{[nt]})_{t\geq0}$, under more general
assumptions.
A close look at the proof (see, in particular, formula (130) in \cite{duq03})
shows that the joint convergence stated in the proposition is indeed
a direct consequence of the arguments in \cite{duq03}.
\begin{lmm}
\label{condiduration2}
The finite-dimensional marginal distributions of the process
\[
\bigl(n^{-1/2\alpha}L^\theta_{[nt]}\bigr)_{0\leq t\leq1}
\]
under $\ov\Q{}^{(n)}(d\theta)$ converge, as $n\to\infty,$
to the finite-dimensional marginal distributions of the process
$(\sqrt{2c_0}D_t)_{0\leq t\leq1}$ under $\mathbf{N}(\cdot\mid\sigma=1)$.
Moreover, this convergence holds jointly with that
of Proposition \ref{condiduration1}.
\end{lmm}
\begin{pf}
This can be derived from the convergence of the rescaled process
$(n^{-1/\alpha}S^\theta_{[nt]})_{0\leq t\leq1}$ in Proposition
\ref{condiduration1}, in the same way as Proposition \ref{finitedim}
was derived from the convergence (\ref{Levyconv}).
The only delicate point is to verify that a suitable analog of Lemma
\ref{keyestimate} holds. To this end, we may argue as follows.
Suppose that we are interested in the finite-dimensional marginal
distribution at times $0\leq t_1<t_2<\cdots<t_p<1$. It then suffices to
prove that an analog of Lem\-ma~\ref{keyestimate} holds for the
vertices $w_0(\theta),w_1(\theta),\ldots, w_{[nt_p]-1}(\theta)$,
which are
the first $[nt_p]$ white vertices of $\theta$ in lexicographical order.
However, the desired property then involves an event that is measurable
with respect to the $\sigma$-field $\g_{[nt_p]}$ and so
we may use Lemmas \ref{absocont} and \ref{boundeddensity}
to see that it is enough to argue
under the probability measure
$\wt\Q^{(n)}$, rather than under $\ov\Q{}^{(n)}$. The same trick
that we used in the proof of Theorem \ref{heightcondi1}
then leads to the desired estimate. The remaining part
of the argument is straightforward and we leave the details
to the reader.
\end{pf}

Before stating and proving the main theorem of this section, we
need to establish an analog
of Lemma \ref{height-contour}. If $\theta$ is a mobile, then we still denote
(with a slight abuse of notation) by
$R_k=R_k(\theta)$ the time of the first visit of $w_k(\theta)$ by the
contour sequence of
$\theta$, for every $k\in\{0,1,\ldots,\#\mathcal{T}^\circ-1\}$.
\begin{lmm}
\label{height-contour-bis}
For every $\eps>0$,
%
\begin{equation}
\label{height-contour-uni2}
\lim_{n\to\infty} \ov\Q{}^{(n)} \biggl(\frac{1}{n} \sup_{0\leq k \leq
n-1} \biggl| R_k - \frac{k}{\beta} \biggr| >\eps\biggr)=0
\end{equation}
and
\[
\lim_{n\to\infty} \ov\Q{}^{(n)} \biggl(\biggl|\frac{1}{n}\#\mathcal{T}- \frac{1}{\beta
}\biggr|>\eps\biggr) = 0 .
\]
\end{lmm}
\begin{pf}
This follows by a minor modification of the proof of Lemma \ref
{height-contour}.
Starting from a forest $\mathbf{F}=(\theta_1,\theta_2,\ldots)$, as
previously, we note that
$\ov\Q{}^{(n)}(d\theta)$ is the distribution of $\theta_1$ under the
conditioned measure
$\mathbb{P}(\cdot\mid\#\mathcal{T}_1^\circ= n)$. Notice that $\mathbb
{P}(\#\mathcal{T}_1^\circ=
n)=\Q(\#\tc=n)=\psi_n(0)$
is of order $n^{-1-1/\alpha}$
when $n$ is large, by (\ref{Kemp}) and (\ref{locallimit}). Thus, we
can use standard large deviations
estimates for sums of independent random variables to verify that, for
every $\vep>0$,
%
\begin{equation}
\label{height-cont-tech9}
\lim_{n\to\infty} \mathbb{P}\Biggl(\frac{1}{n}\sup_{1\leq k\leq n} \Biggl|
\sum_{j=0}^{k-1} \bigl(B(j)+1\bigr)-\frac{k}{\beta} \Biggr|>\vep\Big| \#\mathcal{T}_1^\circ=
n \Biggr)=0.
\end{equation}
Similarly,
\[
\lim_{n\to\infty} \mathbb{P}\Bigl(\sup_{0\leq j\leq n-1} B(j) > n^\vep\big| \#
\mathcal{T}
_1^\circ= n \Bigr)=0.
\]
Furthermore, an analog of (\ref{height-contour-tech1}) follows from
Proposition \ref{condiduration1}, which implies that, for every $\eps>0$,
we have
\[
\mathbb{P}\Bigl(\sup_{0\leq k\leq n-1} H^\circ_k \geq n^{1-1/\alpha+\eps}
\big| \#\mathcal{T}_1^\circ= n \Bigr)\mathop{\longrightarrow}_{n\to\infty} 0 .
\]
The first assertion of the lemma follows from these remarks by the same
arguments as in
the proof of Lemma \ref{height-contour}. The second assertion is a consequence
of (\ref{height-cont-tech9}) since $\#\mathcal{T}_1= \sum_{j=0}^{n-1}
(B(j)+1)$, $\mathbb{P}$-a.s., on $\{\#\mathcal{T}_1^\circ= n\}$.
\end{pf}
\begin{theorem}
\label{condiduration0}
The law of the process
\[
\bigl(n^{-(1-1/\alpha)} H^\theta_{[nt]}, n^{-1/2\alpha} L^{\theta
}_{[nt]} \bigr)_{t\geq0}
\]
under $\ov\Q{}^{(n)}(d\theta)$
converges, as $n\to\infty,$ to the law of the process
\[
\bigl(c_0^{-1} H_t, \sqrt{2c_0} D_t \bigr)_{t\geq0}
\]
under $\mathbf{N}(\cdot\mid\sigma=1)$.
Similarly, the law of the process
\[
\bigl(n^{-(1-1/\alpha)} C^{\theta}_{[nt]}, n^{-1/2\alpha} \Lambda
^{\theta}_{[nt]} \bigr)_{t\geq0}
\]
under $\ov\Q{}^{(n)}(d\theta)$
converges, as $n\to\infty,$ to the law of
\[
\bigl(c_0^{-1} H_{\beta t}, \sqrt{2c_0} D_{\beta t} \bigr)_{t\geq0}
\]
under $\mathbf{N}(\cdot\mid\sigma= 1)$.
\end{theorem}
\begin{pf}
Fix a real $a\in(\frac{1}{2},1)$. Recall that a sequence of laws of
c\`adl\`ag processes
is $C$-tight if it is tight and any sequential limit is supported on
the space
of continuous functions. We first observe
that the sequence of the laws of the processes
%
\begin{equation}
\label{heightlabelduration}
\bigl(n^{-(1-1/\alpha)} H^\theta_{[nt]}, n^{-1/2\alpha} L^{\theta
}_{[nt]} \bigr)_{0\leq t\leq a}
\end{equation}
under $\ov\Q{}^{(n)}(d\theta)$ is $C$-tight. Indeed, by Lemmas \ref{absocont}
and \ref{boundeddensity}, the law under $\ov\Q{}^{(n)}$ of the process
in (\ref{heightlabelduration})
is absolutely continuous with respect to the law of the same
process under $\wt\Q^{(n)}$, with a Radon--Nikodym density that is bounded
uniformly in~$n$. The desired tightness then
follows from Theorem \ref{heightcondi1}.

Next, from Lemma \ref{height-contour-bis}
and the very same arguments as in the derivation of (\ref
{limit-cont1}) and
(\ref{limit-cont2}), we have, for every
$\eps>0$,
that
%
\begin{equation}
\label{condidur1}
\ov\Q{}^{(n)} \Bigl( n^{-(1-1/\alpha)} \sup_{0\leq k\leq{an}/{\beta}}
\bigl|C^\theta_k - H^\theta_{[\beta k]}\bigr| >\eps\Bigr)
\mathop{\la}_{n\to\infty} 0
\end{equation}
and
%
\begin{equation}
\label{condidur2}
\ov\Q{}^{(n)} \Bigl( n^{-1/2\alpha} \sup_{0\leq k\leq{an}/{\beta}}
\bigl|\Lambda^\theta_k - L^\theta_{[\beta k]}\bigr| >\eps\Bigr)
\mathop{\la}_{n\to\infty}0 .
\end{equation}
Note that we must restrict the supremum to $k\leq\frac{a}{\beta}n$ because
we need the $C$-tightness of the processes in (\ref{heightlabelduration}).

From (\ref{condidur1}) and (\ref{condidur2}), together with Lemma
\ref{condiduration2}, we obtain that the finite-dimensional marginal
distributions of the process
%
\begin{equation}
\label{condidur5}
\bigl(n^{-(1-1/\alpha)} C^{\theta}_{[nt]}, n^{-1/2\alpha} \Lambda
^{\theta}_{[nt]} \bigr)_{0\leq t\leq a/\beta}
\end{equation}
under $\ov\Q{}^{(n)}$ converge to those of $(c_0^{-1} H_{\beta t},
\sqrt{2c_0} D_{\beta t})_{0\leq t\leq a/\beta}$ under $\mathbf{N}(\cdot
\mid\sigma=1)$. Moreover, the sequence of the laws of the processes
in (\ref{condidur5}) is $C$-tight, by (\ref{condidur1}),
(\ref{condidur2}) and the tightness of the laws of the processes in
(\ref{heightlabelduration}).

This gives the second convergence stated in the theorem, but only over
the time interval $[0,a/\beta]$. To remove this restriction, we
may argue as follows. From Lemma \ref{height-contour-bis}, we have, for
every $\eps>0$,
\[
\ov\Q{}^{(n)} \biggl(\biggl|\frac{1}{n}\#\mathcal{T}- \frac{1}{\beta}\biggr|>\eps\biggr) \mathop
{\la}_{n\to\infty} 0.
\]
On the other hand, we know that $C^\theta_k=\Lambda^\theta_k=0$ for
every $k\geq\#\mathcal{T}- 1$. Furthermore, a simple argument shows
that the
processes
\[
(C^\theta_k,\Lambda^\theta_k)_{k\geq0} \quad\mbox{and}\quad
\bigl(C^\theta_{(\#\mathcal{T}- 1- k)^+},-\Lambda^\theta_{(\#\mathcal{T}- 1-
k)^+}\bigr)_{k\geq0}
\]
have the same distribution under $\ov\Q{}^{(n)}(d\theta)$. It is an easy
matter to combine these remarks in order to remove the
restriction $t\leq a/\beta$ in the convergence of the processes in
(\ref{condidur5}).

The first convergence of the theorem then follows from the second one,
using the identities $H^\theta_k= C^\theta_{R_k}$ and
$L^\theta_k=C^\theta_{R_k}$, together with Lemma \ref{height-contour-bis}.
\end{pf}

\section{Asymptotics for large planar maps}
\label{large-maps}

In this section, we apply the results of the preceding sections to
properties of planar maps distributed according to
$\mathbf{P}_q$ and conditioned to be large in some sense.
We recall our notation $v_*$ for the distinguished vertex
of a rooted and pointed bipartite planar map $\mathbf{m}$ and
$e_-$ for the origin of the root edge of $\mathbf{m}$. The radius of the
planar map $\mathbf{m}$ is defined by
\[
R(\mathbf{m}) = \max_{v\in V(\mathbf{m})} \dgr(v_*,v).
\]
The profile of distances in $\mathbf{m}$ is the point measure $\rho_\mathbf{m}$
on $\Z_+$ defined by
\[
\rho_\mathbf{m}(k)=\#\{v\in V(\mathbf{m})\dvtx \dgr(v_*,v)=k\},\qquad k\in\Z_+.
\]
Finally, we also set
$\Delta(\mathbf{m})= \dgr(e_-,v_*)$.

In the following theorem, we consider the distance process
$(D_t)_{t\geq0}$ under $\mathbf{N}(\cdot\mid\sigma=1)$ and
under $\mathbf{N}(\cdot\mid\sigma> 1)$. In both cases, we use the notation
\[
\ov D =\max_{t\geq0} D_t ,\qquad \underline D=\min_{t\geq0} D_t .
\]
\begin{theorem}
\label{asymptodistance}
Let $M_n$ be distributed according to $\mathbf{P}_q(\cdot\mid\# V(\mathbf{m})=n)$,
[resp., $\mathbf{P}_q(\cdot\mid\# V(\mathbf{m})\geq n)$]. Then:
\begin{longlist}
\item
${n^{-1/2\alpha} R(M_n)\mathop{\la}\limits_{n\to\infty}\limits^{(d)} \sqrt{2c_0}(
\ov D-
\underline D)}$;
\item if $\rho^{(n)}_{M_n}$ denotes the rescaled profile of
distances
in $M_n$ defined by
\[
\int\rho^{(n)}_{M_n}({d}x) \varphi(x)=n^{-1} \sum_{k\in\Z_+}
\rho_{M_n}(k) \varphi(n^{-1/2\alpha}k),
\]
then $\rho^{(n)}_{M_n}$ converges in distribution to the measure
$\rho^{(\infty)}$ defined by
\[
\int\rho^{(\infty)}({d}x) \varphi(x)
=\int_0^\sigma{d}t \varphi\bigl(\sqrt{2c_0}(D_t-\underline D)\bigr);
\]
\item
${n^{-1/2\alpha} \Delta(M_n) \mathop{\la}\limits_{n\to\infty}\limits^{(d)} \sqrt
{2c_0} \,\ov D.}$
\end{longlist}
In \textup{(i)}--\textup{(iii)}, the limiting distributions are to be
understood under the
probability measure $\mathbf{N}(\cdot\mid\sigma=1)$ [resp., $\mathbf{N}(\cdot\mid\sigma> 1)$].
\end{theorem}
\begin{pf}
Let $M_{n}$ be distributed according to
$\mathbf{P}_q(\cdot\mid\# V(\mathbf{m})=n)$
and let $\theta_n$ be the random mobile associated with $M_{n}$ by
the BDG bijection. By Proposition \ref{mobile-map}, $\theta_n$
is distributed according to $\ov\Q{}^{(n-1)}$. From Proposition
\ref{sec:boutt-di-franc},
\[
R(M_{n})=\ov\ell_n -\underline\ell_n+1 ,
\]
where $\ov\ell_n$ (resp., $\underline\ell_n$) denotes the maximal
(resp., minimal) label in $\theta_n$. It is now clear that
\[
\ov\ell_n -\underline\ell_n= \max_{k\geq0} \Lambda^{\theta_n}_k
- \min_{k\geq0} \Lambda^{\theta_n}_k
\]
and so (i) follows from the second assertion of Theorem \ref{condiduration0}.

Then, let $\varphi$ be a bounded continuous function on $\R_+$. We have
\begin{eqnarray*}
\int\rho^{(n)}_{M_n}({d}x)\varphi(x)&=&
n^{-1} \sum_{v\in V(M_n)} \varphi(n^{-1/2\alpha} \dgr(v_*,v))\\
&=&n^{-1}\sum_{i=0}^{n-2} \varphi\bigl(n^{-1/2\alpha}\bigl(\ell
_n(w_i)-\underline\ell_n+1\bigr)\bigr)\\
&&{} + n^{-1}\varphi(0),
\end{eqnarray*}
where $w_0=w_0(\theta_n), \ldots, w_{n-2}=w_{n-2}(\theta_n)$ denote the
white vertices of $\theta_n$ listed in lexicographical order and
$\ell_n(w_0),\ldots,\ell_n(w_{n-2})$ are their respective labels. Then,
\begin{eqnarray*}
&&n^{-1}\sum_{i=0}^{n-2} \varphi\bigl(n^{-1/2\alpha}\bigl(\ell
_n(w_i)-\underline\ell_n+1\bigr)\bigr) \\
&&\qquad =n^{-1}\sum_{i=0}^{n-2} \varphi\Bigl(n^{-1/2\alpha} \Bigl(L^{\theta_n}_i
-\min_{j=0,\ldots,n-2} L^{\theta_n}_j+1 \Bigr) \Bigr)\\
&&\qquad =\int_0^{1-n^{-1}} {d}t\, \varphi\Bigl(n^{-1/2\alpha} \Bigl(L^{\theta_n}_{[nt]}
-\min_{s\in[0,1]} L^{\theta_n}_{[ns]}+1 \Bigr) \Bigr).
\end{eqnarray*}
The convergence in (ii) is thus a consequence of the first assertion of Theorem
\ref{condiduration0}.

Finally, we have
\[
\Delta(M_n)=1-\underline\ell_n,
\]
except if $v_*=e_-$, in which case $\Delta(M_n)=0=-{\underline\ell}_n$.
Thus, the same argument as for (i) shows that
$n^{-1/2\alpha}\Delta(M_n)$
converges in distribution to $-\sqrt{2c_0} \underline D$, which has the
same law as $\sqrt{2c_0} \,\ov D,$ by symmetry.

The case where $M_n$ is distributed
according to $\mathbf{P}_q(\cdot\mid\# V(\mathbf{m})\geq n)$
is treated by similar arguments, using Theorem
\ref{heightcondi1} instead of Theorem \ref{condiduration0}.
\end{pf}

Recall from \cite{BBI} the notion of the Gromov--Hausdorff distance between
two compact metric spaces. The space
$\K$ of all isometry classes of compact
metric spaces, equipped with the Gromov--Hausdorff distance, is a
Polish space.
If $M$ is a random planar map, then the set $V(M)$
equipped with the metric $\dgr$ is a random variable with values in
$\K$.
\begin{theorem}
\label{GHconv}
For every $n\geq1$, let $M_n$ be distributed
according to $\mathbf{P}_q(\cdot\mid\break\# V(\mathbf{m})=n)$
[resp., $\mathbf{P}_q(\cdot\mid\# V(\mathbf{m})\geq n)$].
From every strictly increasing sequence of integers, one can extract a
subsequence along which
\[
(V(M_n),n^{-1/2\alpha}\dgr)
\mathop{\la}_{n\to\infty}^{(d)} (\mathbf{M}_\infty,\delta_\infty),
\]
where $(\mathbf{M}_\infty,\delta_\infty)$ is a random compact metric
space and the
convergence holds in distribution, in the Gromov--Hausdorff sense.
Furthermore, the Hausdorff dimension of $(\mathbf{M}_\infty,\delta
_\infty)$
is a.s. equal to $2\alpha$.
\end{theorem}
\begin{pf}
We consider only the case where $M_n$ is distributed
according to $\mathbf{P}_q(\cdot\mid\# V(\mathbf{m})=n)$. The first
assertion could be established
by using compactness criteria in the space $\K$ in order to derive the
tightness of the
distributions of the spaces $(V(M_n),n^{-1/2\alpha}\dgr)$. We will
use a different approach,
which is inspired by~\cite{legall06}, Section 3. This approach
provides additional
information about the limiting space $(\mathbf{M}_\infty,\delta_\infty
)$, which will be useful when
proving the second assertion of the theorem.

As in the previous proof,
let $\theta_n$ be the random mobile associated with $M_{n}$ by
the BDG bijection and write $v^n_0,v^n_1,\ldots,v^n_{r_n}$ for the
white contour sequence of $\theta_n$. Recall that the BDG bijection allows
us to identify the white vertices of $\theta_n$
with corresponding vertices of the map $M_n$. We can thus set, for
every $i,j\in\{0,1,\ldots,r_n\}$,
\[
d_n(i,j)=\dgr(v^n_i,v^n_j),
\]
where $\dgr$ refers to the graph distance in the map $M_n$. By convention,
we put $v^n_k=v^n_{r_n}=\varnothing$ for every $k\geq r_n$
so that the definition of $d_n(i,j)$ makes sense for all nonnegative integers
$i$ and $j$. We can use linear interpolation to extend the definition
of $d_n$ to real values of the parameters, by setting, for every
$s,t\geq0$,
\begin{eqnarray*}
d_n(s,t)&=&(s-[s])(t-[t])d_n([s]+1,[t]+1)\\
&&{} + (s-[s])([t]+1-t)d_n([s]+1,[t])\\
&&{} +([s]+1-s)(t-[t])d_n([s],[t]+1)\\
&&{} + ([s]+1-s)([t]+1-t)d_n([s],[t]).
\end{eqnarray*}

By \cite{legall06}, Lemma 3.1, we have, for all integers $i,j\geq0$,
%
\begin{equation}
\label{upperbounddist}
d_n(i,j)\leq d^0_n(i,j),
\end{equation}
where
\[
d_n^0(i,j)=\Lambda^{\theta_n}_i+ \Lambda^{\theta_n}_j
-2 \min_{i\wedge j\leq k\leq i\vee j} \Lambda^{\theta_n}_k +2.
\]
(To be precise, \cite{legall06} uses a slightly different version of
the BDG bijection, with nonnegative labels, but is straightforward to
verify that the argument of the proof of Lemma 3.1 in \cite{legall06}
goes through without change in our setting.)
In the same way as for $d_n$, we extend the definition of $d^0_n$
to real values of the parameters by linear interpolation. The bound
$d_n(s,t)\leq d^0_n(s,t)$ still holds for real $s$ and $t$.

Let $(H^{(1)}_t,D^{(1)}_t)_{t\geq0}$ be a random process which has the
distribution of $(H_t$,\break $D_t)_{t\geq0}$ under $\mathbf{N}(\cdot\mid\sigma
=1)$. From Theorem \ref{condiduration0},
%
\begin{equation}
\label{convdist0}
(n^{-1/2\alpha}d^0_n(ns,nt) )_{s,t\geq0}
\mathop{\la}_{n\to\infty}^{\mathrm{(d)}}
\bigl(\sqrt{2c_0} d_\infty^0(\beta s,\beta t) \bigr)_{s,t\geq0},
\end{equation}
where, for every $s,t\geq0$,
\[
d_\infty^0(s,t)= D^{(1)}_s + D^{(1)}_t -
2\min_{s\wedge t\leq r\leq s\vee t} D^{(1)}_r .
\]
In (\ref{convdist0}), the convergence holds, in the sense of weak convergence
of the laws in the space of continuous functions on $\R_+^2$.

We then observe that, for every $s,t,s',t'\geq0$,
%
\begin{equation}
\label{tightnessdist}
|d_n(s,t)-d_n(s',t')|\leq d_n(s,s') + d_n(t,t')\leq d^0_n(s,s')+ d^0_n(t,t').
\end{equation}
By the convergence (\ref{convdist0}),
we have, for every $\eta,\eps>0$,
\[
\limsup_{n\to\infty}
P \Bigl(\sup_{|s-s'|\leq\eta} n^{-1/2\alpha} d^0_n(ns,ns') \geq\eps\Bigr)
\leq P \biggl(\sup_{|s-s'|\leq\eta} d^0_\infty(\beta s,\beta s') \geq
\frac{\eps}{\sqrt{2c_0}} \biggr).
\]
If $\eps>0$ is fixed, then the right-hand side can be made arbitrarily small
by choosing $\eta>0$ to be small enough. Thanks to this remark and to
the bound
(\ref{tightnessdist}), one easily gets that the sequence of the laws
of the
processes
\[
(n^{-1/2\alpha}d_n(ns,nt) )_{s,t\geq0}
\]
is tight (see the proof of Proposition 3.2 in \cite{legall06}
for more details).

Also using Theorem \ref{condiduration0}, we obtain that, from any
strictly increasing sequence of positive integers, we
can extract a subsequence $(n_k)_{k\geq1}$ along which we have
the joint convergence
%
\begin{eqnarray}
\label{convdist}
&& \bigl(n^{-(1-1/\alpha)}C^{\theta_n}_{[nt]},n^{-1/2\alpha}\Lambda
^{\theta_n}_{[nt]},
n^{-1/2\alpha}d_n(ns,nt) \bigr)_{s,t\geq0}\nonumber\\[-8pt]\\[-8pt]
&&\qquad \mathop{\la}_{n\to\infty}^{\mathrm{(d)}} \bigl(c_0^{-1} H^{(1)}_{\beta t},\sqrt
{2 c_0}D^{(1)}_{\beta t},\sqrt{2 c_0}
d_\infty(\beta s, \beta t) \bigr)_{s,t\geq0},\nonumber
\end{eqnarray}
where $(d_\infty(s,t))_{s,t\geq0}$ is a continuous random process
indexed by $\R_+^2$ and
taking nonnegative values. From now on, we restrict our attention to
values of $n$ belonging to the sequence $(n_k)$.

By the Skorokhod representation theorem, we may, and will, assume that
the convergence (\ref{convdist}) holds almost surely. Note that the
bound $d_n\leq d^0_n$ immediately gives $d_\infty\leq d^0_\infty$.
From the convergence (\ref{convdist}), one also gets that
the function $(s,t)\to d_\infty(s,t)$ is symmetric and satisfies the triangle
inequality. Furthermore, the bound $d_\infty\leq d^0_\infty$
implies that $d_\infty(s,t)=0$ if $s\geq1$ and $t\geq1$. We define an
equivalence relation on $[0,1]$ by setting
\[
s\approx t \quad\mbox{if and only if}\quad d_\infty(s,t)=0.
\]
We let $\mathbf{M}_\infty$ be the quotient space $[0,1]/\approx$
and equip $\mathbf{M}_\infty$ with the metric $\delta_\infty=\sqrt
{2c_0} d_\infty$.
The continuity of $d_\infty$ ensures that the canonical
projection from $[0,1]$ (equipped with the usual metric) onto
$\mathbf{M}_\infty$ is continuous, so $\mathbf{M}_\infty$
is compact.

We claim that the convergence of the theorem holds almost surely [along
the sequence $(n_k)$] with this choice of the space
$(\mathbf{M}_\infty,\delta_\infty)$. To see this, define a correspondence
$\mathcal{C}_n$ between $(V({M}_n)\setminus\{v_*\},n^{-1/2\alpha
}\dgr)$ and $(\mathbf{M}_\infty,
\delta_\infty)$ by declaring that
a vertex $v$ of $V({M}_n)\setminus\{v_*\}$ is in correspondence with
$x\in\mathbf{M}_\infty$
if and only if there exists a representative $s$ of $x$ in $[0,1]$ such
that $v=v^n_{[ns/\beta]}$.
The desired convergence will follow if we can verify that the
distortion of $\mathcal{C}_n$ tends to $0$ as $n\to\infty$. To this end,
let $s,s'\in[0,1]$ and set $k=[ns/\beta]$ and $k'=[ns'/\beta]$.
If $v=v^n_k$ and $v'=v^n_{k'}$, and if $x$ and $x'$ are the respective
equivalence classes of $s$ and $s'$ in the quotient $[0,1]/\approx$,
then we have
\begin{eqnarray*}
&&\bigl|n^{-1/2\alpha} \dgr(v,v') -\sqrt{2c_0} d_\infty(x,x')\bigr|\\
&&\qquad =\bigl|n^{-1/2\alpha} d_n(k,k')-\sqrt{2c_0} d_\infty(s,s')\bigr|\\
&&\qquad =\biggl|n^{-1/2\alpha} d_n\biggl(\biggl[\frac{ns}{\beta}\biggr],\biggl[\frac{ns'}{\beta
}\biggr]\biggr)-\sqrt{2c_0} d_\infty(s,s')\biggr|\\
&&\qquad \leq\sup_{t,t'\geq0}
\bigl|n^{-1/2\alpha} d_n([nt],[nt'])-\sqrt{2c_0} d_\infty(\beta t, \beta t')\bigr|,
\end{eqnarray*}
which tends to $0$ as $n\to\infty$, by the (almost sure) convergence
(\ref{convdist}). This completes the
proof of the first assertion of the theorem.

Let us now turn to the Hausdorff dimension
of $(\mathbf{M}_\infty,\delta_\infty)$. From the bound $d_\infty\leq
d^0_\infty$ and the H\"older continuity properties of
the distance process, we get that for every $\eps\in(0,1/2\alpha),$
there is an
almost surely finite random constant $K_{(\eps)}$ such that, for every
$s,t\in[0,1]$,
\[
d_\infty(s,t)\leq K_{(\eps)} |t-s|^{(1/2\alpha)-\eps}.
\]
Hence, the projection mapping from $[0,1]$ onto $\mathbf{M}_\infty$ is
a.s. H\"older
continuous with exponent $(1/2\alpha)-\eps$. The almost sure bound
$\operatorname{dim}(\mathbf{M}_\infty,\delta_\infty)
\leq2\alpha$ immediately follows.

The proof of the lower bound $\operatorname{dim}(\mathbf{M}_\infty,\delta_\infty
)\geq2\alpha$ is more delicate.
We start with a useful lower bound on $d_\infty$.
\begin{lmm}
\label{lowerbounddist}
Almost surely, for every $0<s<t<1$ and $r\in(s,t)$ such that
$H^{(1)}_u>H^{(1)}_r$
for every $u\in[s,r)$, we have
\[
d_\infty(s,t)\geq D^{(1)}_s - D^{(1)}_r.
\]
Similarly, almost surely for every $0<t<s<1$ and $r\in(t,s)$ such that
$H^{(1)}_u>H^{(1)}_r$
for every $u\in(r,s]$, we have
\[
d_\infty(s,t)\geq D^{(1)}_s - D^{(1)}_r.
\]
\end{lmm}
\begin{pf}
We establish only the first assertion since the proof of the second one is
very similar. So, let $s,t,r$ be as in the first part of the lemma.
Let $(k_n)$ and $(k'_n)$ be two sequences of positive integers such that
$n^{-1}k_n\la\beta^{-1}s$ and $n^{-1}k'_n\la\beta^{-1}t$ as $n\to
\infty$ (both
sequences are indexed by the set of values of $n$ that we are considering).
Thanks to the convergence (\ref{convdist}) and our assumption
$H^{(1)}_u>H^{(1)}_r$
for every $u\in[s,r)$, we can find
another sequence $(m_n)$ of positive integers such that $n^{-1}m_n\la
\beta^{-1}r$
and, for $n$ large enough,
\[
C^{\theta_n}_j > C^{\theta_n}_{m_n} > \min_{i\in\{k_n,\ldots,k'_n\}
} C^{\theta_n}_i\qquad
\forall j\in\{k_n,\ldots,m_n-1\}.
\]
Recall our notation $v^n_0,v^n_1,\ldots$ for the white contour
sequence of $\theta_n$. The
preceding inequalities imply that $v^n_{m_n}$ is an ancestor of
$v^n_{k_n}$, but not an ancestor of
$v^n_{k'_n}$. Let $\ga_n=(\gamma_n(i),0\leq i \leq d_{\mathrm{gr}}
(v^n_{k_n},v^n_{k'_n}))$ be a geodesic
from $v^n_{k_n}$ to $v^n_{k'_n}$ in the planar map $M_n$ and let $i_n$
be the
largest integer $i\in\{0,1,\ldots, d_{\mathrm{gr}}(v^n_{k_n},v^n_{k'_n})\}
$ such that $\gamma_n(i)$
is a descendant of $v^n_{m_n}$. By the preceding remarks, we have
$0\leq i_n< d_{\mathrm{gr}}(v^n_{k_n},v^n_{k'_n})$. Furthermore, the contour
sequence of $\theta_n$
must visit $v^n_{m_n}$ between any time at which it visits the point
$\gamma_n(i_n)$
and any other time at which it visits $\gamma_n(i_n+1)$. Using the
construction of
edges in the BDG bijection, the existence of an
edge of $M_n$ between $\ga_n(i_n)$ and $\ga_n(i_n+1)$ implies that
\[
\ell_n(v^n_{m_n})\geq\ell_n(\ga_n(i_n)).
\]
It follows that
\begin{eqnarray*}
d_n(k_n,k'_n)&=&d_{\mathrm{gr}}(v^n_{k_n},v^n_{k'_n})\geq d_{\mathrm{gr}}(v^n_{k_n},\ga_n(i_n))\\
&\geq& d_{\mathrm{gr}}(v_*,v^n_{k_n}) - d_{\mathrm{gr}}(v_*,\ga_n(i_n))\\
&=&\ell_n(v^n_{k_n})-\ell_n(\ga_n(i_n))\\
&\geq&\ell_n(v^n_{k_n})- \ell_n(v^n_{m_n})\\
&=&\Lambda^{\theta_n}_{k_n} - \Lambda^{\theta_n}_{m_n}.
\end{eqnarray*}
The bound of the lemma follows by passing to the limit $n\to\infty$
using (\ref{convdist}).
\end{pf}

The next lemma will be used in combination with Lemma \ref{lowerbounddist}
to estimate the size of balls for the metric $\delta_\infty$. For
technical reasons,
we prove this lemma under the excursion measure $\bN$ and we will then
use a
scaling argument to get a similar result under $\bN(\cdot\mid\sigma=1)$.
For every $u>0$, $\lambda_u({d}s)$ denotes Lebesgue measure on $(0,u)$.
\begin{lmm}
\label{KeyHausdim}
For every $s\in(0,\sigma)$, set
\[
\mathcal{I}(s)=\{r\in[s,\sigma]\dvtx H_u>H_r\mbox{ for every }u\in
[s,r)\}
\]
and for every $\vep>0$,
set
\[
\tau^s_\vep=\inf\{t\in\mathcal{I}(s)\dvtx D_t\leq D_s-\vep\},
\]
where $\inf\varnothing= \infty$. Then,
for every $a\in(0,2\alpha)$,
\[
\lim_{\vep\da0} \vep^{-a} (\tau^s_\vep- s) =0 ,\qquad \lambda_\sigma
({d}s)\mbox{ a.e., }
\bN\mbox{ a.e.}
\]
\end{lmm}
\begin{pf}
For $s\in(0,\sigma)$ and $r\in[0,H_s)$, set
\[
\ga^s_r=\inf\{t\geq s \dvtx H_t < H_s -r\}.
\]
By convention, we put $\ga^s_r=\sigma$ if $r\geq H_s$. For our purposes,
it will be important to have information on the sample path behavior
of the function $r\la D_{\ga^s_r}$. This is the goal of the next
lemma, which relies heavily on
results from \cite{duqleg02}, to which we refer for
additional details. We first need to introduce some notation. For every
$s\in(0,\sigma)$, we define two positive finite measures on
$(0,\infty)$ by setting
\begin{eqnarray*}
\rho_s&=&\sum_{0\leq u\leq s} (I^u_s-X_{u-}) \ind_{\{X_{u-}<I^u_s\}}
\delta_{H_u} ,\\
\eta_s&=&\sum_{0\leq u\leq s} (X_u-I^u_s) \ind_{\{X_{u-}<I^u_s\}}
\delta_{H_u} .
\end{eqnarray*}
(It is not immediately obvious that $\eta_s$ is a finite measure; see
Chapter 3
of \cite{duqleg02}.) One can prove that, $\bN$ a.e., for every $s>0$,
the topological support of $\rho_s$
is $[0,H_s]$ and $\rho_s([0,H_s])=X_s$ (see Chapter 1 of \cite{duqleg02}).
Furthermore, the quantities $H_u$ corresponding to the values of $u$ that
give nonzero terms in the definition of $\rho_s$ are all distinct.

We denote by $\mathcal{N}({d}r\,{d}z\,{d}x)$ a Poisson point
measure on $[0,\infty)^3$
with intensity
\[
{d}r \,\pi({d}z) \ind_{[0,z]}(x) \,{d}x,
\]
where $\pi$ denotes the L\'evy measure of $X$. We can enumerate atoms
of $\mathcal N$
in a measurable way and write
\[
\mathcal{N}=\sum_{j\in J} \delta_{(r_j,z_j,x_j)}.
\]
\begin{lmm}
\label{momentrhoeta}
\textup{(i)} Let $\Phi$ be a nonnegative measurable function on $\R
_+\times M_f(\R_+)^2$. Then,
\[
\bN\biggl(\int_0^\sigma
{d}s\, \Phi(H_s,\rho_s,\eta_s) \biggr)
=\int_0^\infty{d}u\, \E\biggl[ \Phi\biggl(u, \sum_{0\leq r_j\leq u} x_j
\delta_{r_j},
\sum_{0\leq r_j\leq u} (z_j-x_j) \delta_{r_j} \biggr) \biggr].
\]

\textup{(ii)} Let $F$ be a nonnegative measurable function on $\D(\R)$. Then,
\[
\bN\biggl(\int_0^\sigma
{d}s\, F\bigl((D_s-D_{\ga^s_r})_{r\geq0}\bigr) \biggr)
=\int_0^\infty{d}u\, \E[ F((Z_{r\wedge u})_{r\geq0}) ],
\]
where $(Z_r)_{r\geq0}$ is a symmetric stable process with index
$2(\alpha-1)$.
\end{lmm}
\begin{pf}
Part (i) is a special case of Proposition 3.1.3 of \cite{duqleg02}.
Part (ii)
is essentially a consequence of (i) and our construction of the distance
process. Let us explain this in greater detail. We fix $s>0$, $r>0$
and argue on the
event $\{s<\sigma\}$. As in Section \ref{sec:continuous-object}, we assign a
Brownian bridge $b_u$ with length $\Delta X_u$ to each jump time $u$ of
$X$, in such a way that
\[
D_s = \sum_{u\leq s} b_u(I^u_s-X_{u-}) \ind_{\{X_{u-}<I^u_s\}}.
\]
We then also have, $\bN$ a.e.,
\[
D_{\ga^s_r} = \sum_{u\leq s} b_u(I^u_s-X_{u-}) \ind_{\{X_{u-}<I^u_s\}}
\ind_{\{H_u<H_s-r\}}.
\]
To see this, note that the identity
%
\begin{equation}
\label{momenttec1}
\ga^s_r=\inf\{t\geq s\dvtx X_t < X_s-\rho_s([H_s-r,H_s])\}
\end{equation}
is a consequence of formula (20) in \cite{duqleg02}. Moreover, by the
same formula,
$\rho_{\ga^s_r}$ is exactly the restriction of $\rho_s$ to the
interval $[0,H_s-r)$ (or the
zero measure if $r\geq H_s$). Hence, the values $u\leq\ga^s_r$ that give
a nonzero contribution to the sum defining $D_{\ga^s_r}$ are exactly
those $u\leq s$ such that
$X_{u-}<I^u_s$ and $H_u<H_s-r$, leading to the stated formula for
$D_{\ga^s_r}$.

It follows that
%
\begin{equation}
\label{momenttec2}
D_s - D_{\ga^s_r}
=\sum_{u\leq s} b_u(I^u_s-X_{u-}) \ind_{\{X_{u-}<I^u_s\}} \ind_{\{
H_s-r\leq H_u\leq H_s\}}
\end{equation}
and we can use part (i) to compute the Fourier transform of this quantity.
Note that, for every
jump time $u\leq s$ with the property $X_{u-}<I^u_s$, the duration of the
bridge $b_u$ is the sum of the masses assigned by $\rho_s$ and $\eta_s$,
respectively, to the point $H_u$.

Suppose that,
conditionally given $\mathcal{N}$, we are given a collection
$(b_j^{(z_j)})_{j\in J}$ of independent Brownian
bridges, with respective durations $(z_j)_{j\in J}$.
It then follows from (i), formula (\ref{momenttec2}) and the preceding
discussion that, for every $\lambda\in\R$,
\begin{eqnarray*}
&&\bN\biggl( \int_0^\sigma{d}s \exp\bigl(i\lambda(D_s- D_{\ga^s_r})\bigr) \biggr)\\
&&\qquad = \int_0^\infty{d}u\, \E\biggl[ \exp\biggl(i\lambda
\sum_{u-r\leq r_j\leq u} b_j^{(z_j)}(x_j) \biggr) \biggr]\\
&&\qquad = \int_0^\infty{d}u\,
\E\biggl[\exp\biggl(-\frac{\lambda^2}{2} \sum_{u-r\leq r_j\leq u} \frac
{x_j(z_j-x_j)}{z_j} \biggr) \biggr]\\
&&\qquad = \int_0^\infty{d}u\,
\E\biggl[\exp\biggl( - \int_{(u-r)_+}^u {d}v\int\pi({d}z)\\
&&\qquad\quad\hspace*{135.6pt}{}\times\int_0^z
{d}x \biggl(1-\exp\biggl(-\frac{\lambda^2}{2}
\frac{x(z-x)}{z}\biggr)\biggr) \biggr) \biggr]\\
&&\qquad =\int_0^\infty{d}u \exp\bigl(- K_\alpha(u\wedge r) |\lambda
|^{2(\alpha-1)}\bigr),
\end{eqnarray*}
by an easy calculation, using the fact that $\pi({d}z)=K'_\alpha
z^{-1-\alpha}\,{d}z$.

It follows that the formula of (ii) holds in the case where $F$ is of
the form
$F(\omega)=f(\omega(r))$ for a fixed $r>0$. A slight extension of the previous
calculation gives the case where $F$ depends only on a finite number of
coordinates.
This is enough to conclude since the process $(D_s-D_{\ga^s_r})_{r\geq0}$
has right-continuous paths.
\end{pf}

We now complete the proof of Lemma \ref{KeyHausdim}. We fix $a\in
(0,2\alpha)$.
We can then choose $b\in((2\alpha-2)^{-1},\infty)$, $b'\in
(0,(\alpha-1)^{-1})$
and $b''\in(0,\alpha)$ such that
\[
\frac{b'b''}{b} > a.
\]
By standard path properties of stable processes (see, e.g., \cite
{bertlev96}, Theorem VIII.6), we have
\[
\lim_{r\da0}r^{-b} \Bigl(\sup_{0\leq x\leq r} Z_x \Bigr)= \infty\qquad\mbox{a.s.}
\]
It then follows from Lemma \ref{momentrhoeta}(ii) that we also have
\[
\lim_{r\da0}r^{-b} \Bigl(\sup_{0\leq x\leq r} (D_s-D_{\ga^s_x}) \Bigr)=
\infty,\qquad
\lambda_\sigma({d}s)\mbox{ a.s., } \bN\mbox{ a.e.}
\]
Notice that $\ga^s_x\in\mathcal{I}(s)$ provided that $x$ is a
continuity point of the
mapping $r\to\ga^s_r$ and thus for all but countably many values of
$x$. Therefore,
the previous display also implies that
%
\begin{equation}
\label{Haustec2}
\tau^s_\vep\leq\ga^s_{\vep^{1/b}}
\end{equation}
for all sufficiently small $\vep>0$, $\lambda_\sigma({d}s)$
a.e., $\bN$ a.e.

The next step is to investigate the behavior of $\ga^s_x$
as $x\da0$. We first observe that
%
\begin{equation}
\label{Haustec3}
\lim_{x\da0} x^{-b'} \rho_s([H_s-x,H_s])=0 ,\qquad
\lambda_\sigma({d}s)\mbox{ a.s., } \bN\mbox{ a.e.}
\end{equation}
This is a consequence of Lemma \ref{momentrhoeta}(i): note that, for
every $u>0$, the process
\[
Y_x=\sum_{u-x\leq r_j\leq u} x_j,\qquad 0\leq x\leq u,
\]
is a stable subordinator with index $\alpha- 1$ and apply path
properties of
subordinators (see, e.g., \cite{bertlev96}, Theorem VIII.5).
Furthermore, by applying the
Markov property under $\bN$ and again using \cite{bertlev96}, Theorem
VIII.6, we get that
\[
\lim_{r\da0} r^{-1/b''} \sup_{0\leq x\leq r} (X_s - X_{s+x}) =
\infty,
\]
$\bN$ a.e. on $s<\sigma$, for every fixed $s>0$. It readily follows that
%
\begin{equation}
\label{Haustec4}
\inf\{x\geq0\dvtx X_{s+x}< X_s-r\}\leq r^{b''}
\end{equation}
for all sufficiently small $r>0$, $\lambda_\sigma({d}s)$ a.e.,
$\bN$ a.e.
Now, recall (\ref{momenttec1}) and use (\ref{Haustec3}) and (\ref{Haustec4})
to obtain
%
\begin{equation}
\label{Haustec5}
\ga^s_r\leq s + r^{b'b''}
\end{equation}
for all sufficiently small $r>0$, $\lambda_\sigma({d}s)$ a.e.,
$\bN$ a.e.
We get the statement of the lemma by combining
(\ref{Haustec2}) and (\ref{Haustec5}), recalling that $b'b''/b >a$.
\end{pf}

We now complete the proof of Theorem \ref{GHconv}. We again fix $a\in
(0,2\alpha)$. For every $s\in(0,1)$, we set
\[
\wt{\mathcal{I}}(s)=\bigl\{r\in[s,1]\dvtx H^{(1)}_u>H^{(1)}_r\mbox{ for
every }u\in[s,r)\bigr\}
\]
and for every $\vep>0$,
we set
\[
\wt\tau^s_\vep=\inf\bigl\{t\in\wt{\mathcal{I}}(s)\dvtx D^{(1)}_t\leq
D^{(1)}_s -\vep\bigr\}.
\]
From Lemma \ref{KeyHausdim} and an easy scaling argument, we get
\[
\lim_{\vep\da0} \vep^{-a} (\wt\tau^s_\vep- s) =0,\qquad \lambda_1({d}s)\mbox{ a.e., }
\mbox{ a.s.}
\]

However, if $\wt\tau^s_\vep\leq t < 1$, the first part of Lemma \ref
{lowerbounddist} implies that
$d_\infty(s,t)\geq\vep$. Thus,
\[
\int_s^1 {d}t\, \ind_{\{d_\infty(s,t)<\vep\}} \leq\wt\tau
^s_\vep-s
\]
and
\[
\lim_{\vep\da0} \vep^{-a}\int_s^1 {d}t\, \ind_{\{d_\infty
(s,t)<\vep\}}=0,\qquad
\lambda_1({d}s)\mbox{ a.e., a.s.}
\]
We can use a symmetric argument to handle the analogous integral where
$t$ varies between $0$ and $s$: use the second part of Lemma \ref
{lowerbounddist}
and note that the distribution of the pair
$(H^{(1)}_{t},D^{(1)}_t)_{0\leq t\leq1}$
is invariant under the change of parameter $t\to1-t$. We thus conclude that
\[
\lim_{\vep\da0} \vep^{-a}\int_0^1 {d}t\, \ind_{\{d_\infty
(s,t)<\vep\}}=0,\qquad
\lambda_1({d}s)\mbox{ a.e., a.s.}
\]
Finally, if $\kappa$ denotes the probability measure on $M_\infty$
which is the
image of Lebesgue measure on $(0,1)$ under the canonical projection,
then we see that
\[
\lim_{\vep\da0}
\frac{\kappa(B_\infty(x,\vep))}{\vep^a} =0,\qquad \kappa({d}x)\mbox
{ a.e., a.s.,}
\]
where $B_\infty(x,\vep)=\{y\in M_\infty\dvtx \delta_\infty(x,y)<\vep\}$.

The lower bound
$\operatorname{dim}(\mathbf{M}_\infty,\delta_\infty)\geq2\alpha$ now follows
from standard density theorems
for Hausdorff measures.
\end{pf}
\begin{Remark*}
As we already noted in Section \ref{sec1}, the results of this
section carry over to
Boltzmann distributions on nonpointed rooted planar maps. More precisely,
denote by $\wt W_{q}$ the Boltzmann distribution defined as in (\ref
{Bolt-weight}),
but now viewed as a measure on the set of all rooted planar maps.
Let $\wt M_n$ be a random rooted planar map distributed according
to the (suitably normalized) restriction of $\wt W_q$ to maps with $n$ vertices.
Then, Theorem \ref{asymptodistance} gives information about the distances
in $\wt M_n$ from a vertex chosen uniformly at random and both
assertions of
Theorem \ref{GHconv} remain valid if $M_n$ is replaced by $\wt M_n$.
\end{Remark*}

\section{Some motivation from
physics}\label{sec:some-motivation-from}

In this section, we describe a motivation for the models discussed in
this article that comes from the physics literature. In this
discussion, we rely on a number of nonrigorous predictions and our
only goal is to isolate some possible directions for future work. A
useful reference is Appendix B in the survey by Duplantier
\cite{duplantier04} and the references therein.

\begin{figure}[b]

\includegraphics{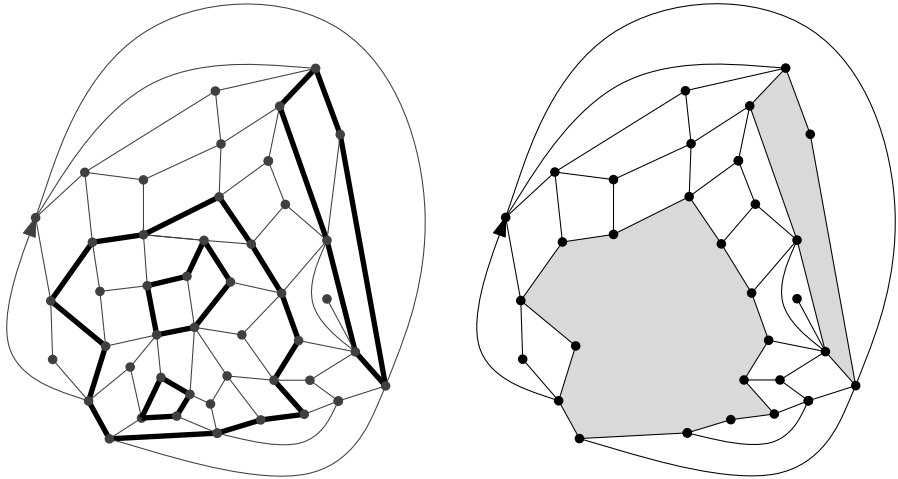}

\caption{An $O(N)$ configuration on a rooted quadrangulation, with $4$
cycles of total length $30$, and the external gasket associated with
this configuration, with shaded holes of degrees $6$ and $14$.}
\label{fig:onpic}
\end{figure}

As a starting point, we observe that models of random maps that are
very similar to
ours appear when studying annealed statistical physics
models on random maps. These models are similar to more familiar
models on regular lattices, such as percolation and Ising or Potts
models, but they are defined on a random map that is
chosen at the same time as the configuration of the model. To illustrate
this, we will first deal with the so-called $O(N)$ model on a random planar
quadrangulation. Let $\bq$ be a rooted quadrangulation. A \textit{loop
configuration}
on $\bq$ is a collection $\mathcal{L}=\{c_1,\ldots,c_k\}$, where
$c_1,\ldots,c_k$ are cycles, that is, paths on $\bq$ starting
and ending at the same point and never visiting the same vertex
twice. It is further required that the paths $c_i$ do not intersect.
We set
\[
\#\mathcal{L}=k \quad\mbox{and}\quad \operatorname{lg}(\mathcal{ L})=\sum
_{i=1}^k\operatorname{lg}(c_i) ,
\]
where $\lg(c_i)$ is the number
of edges in the path $c_i$; see Figure \ref{fig:onpic} for an
example.

%

Let $N\geq0$ be fixed. The annealed $O(N)$ measure is the
$\sigma$-finite measure over the set of all pairs $(\bq,\mathcal
{L})$, where
$\bq$ is a rooted quadrangulation and $\mathcal{L}$ is a loop
configuration on $\bq$, defined by
\[
W_{O(N)}(\bq,\mathcal{L})=e^{-\beta\#F(\bq)}x^{\operatorname{lg}(\mathcal{
L})}N^{\#\mathcal{L}} ,
\]
where $\beta$ and $x$ are positive
parameters. When the total mass $Z_{O(N)}(\beta,x)$ of $W_{O(N)}$ is
finite, we say
that the pair $(\beta,x)$ is \textit{admissible} and we can consider
the probability
measure $P_{O(N)}=Z_{O(N)}(\beta,x)^{-1}W_{O(N)}$.

Consider a configuration $(\bq,\mathcal{L})$. A cycle $c\in\mathcal{L}$
splits the sphere into two components. The one that contains the
face located to the left of the root edge of $\bq$ is called the
\textit{exterior} of $c$. The other component is called the \textit
{interior} of $c$. The \textit{external gasket} $\mathcal{E}(\bq,\mathcal{L})$ is the rooted planar
map obtained
from $\bq$ by deleting all the edges and vertices strictly contained
in the interior of some $c\in\mathcal{L}$; see Figure~\ref{fig:onpic}.

More precisely, ${\mathbf{m}}$ is defined as a rooted planar map with
two different
types of faces:
\begin{itemize}
\item faces that came from
the exterior of cycles of $\mathcal{L}$, which have degree $4$---we
denote by
$Q({\mathbf{m}})$ the set of all these faces;
\item faces of arbitrary even degree, called the \textit{holes} of
${\mathbf{m}}$, which came from the
deletion of the interior of a cycle of $\mathcal{L}$---we denote by
$H({\mathbf{m}})$
the set of all holes of ${\mathbf{m}}$ (note that certain holes may
have degree $4$).
\end{itemize}
Furthermore, the boundaries of the holes of ${\mathbf{m}}$ are
disjoint cycles. In particular, every edge of the boundary of a
hole is adjacent to a face of $Q({\mathbf{m}})$.

One can verify that the range of the external gasket mapping $(\bq
,\mathcal{L})\to\mathcal{E}(\bq,\mathcal{L})$
is the set of all rooted planar maps (with faces of two types)
satisfying the preceding conditions.
It is then an easy exercise to check that the push-forward of $W_{O(N)}$
under the external gasket mapping is
%
\begin{equation}
\label{model-O(N)}
W_{O(N)}\bigl(\{\mathcal{E}(\bq,\mathcal{L})={\mathbf{m}}\}\bigr)
=e^{-\beta\#Q({\mathbf{m}})}\prod_{f\in H({\mathbf{m}})}q_{\deg
f/2} ,
\end{equation}
where
\[
q_k= x^{2k}Z^\partial_{O(N),k}(\beta,x)
\]
and
$Z^\partial_{O(N),k}(\beta,x)$ is the partition function for the
$O(N)$-model with a boundary of length $2k$. This partition function is
defined in an
analogous way as $Z_{O(N)}(\beta,x)$, but configurations
$(\bq,\mathcal{L})$ now consist of rooted quadrangulations $\bq$
with a
boundary of length $2k,$ together with a collection $\mathcal{L}$
of disjoint cycles that do not intersect the boundary and such that
the boundary face lies on the left of the root edge. From formula
(\ref{model-O(N)}), we see that the
external gasket of a $P_{O(N)}$-distributed random map has a Boltzmann
distribution of a similar kind as those studied in the present work,
except that the maps that appear here have two distinct types of faces
and extra topological
constraints.

Ignoring these extra constraints, one can conjecture that for
appropriate values of $\beta$ and $x$, the scaling limits of these random
gasket configurations will be closely related to those depicted in Section
\ref{large-maps}, provided that the weights $q_k$ satisfy similar
asymptotics as in Section \ref{sec:choos-boltzm-weights}. At this
stage, some predictions from theoretical physics provide insight into
these questions. For fixed $\beta$ and $x$, we introduce the
generating function
\[
Z^\partial_{O(N)}(z)=\sum_{k\geq1}z^kZ^\partial_{O(N),k}(\beta,x)
.
\]
According to singularity analysis, for $a\in(3/2,2)\cup(2,5/2)$,
a behavior
\[
Z^\partial_{O(N)}(z)\mathop{\approx}_{z\uparrow z_c} (z_c-z)^{a-1}
,
\]
meaning that the singular part of $Z^\partial_{O(N)}$ near its
first positive singularity $z_c$ is of order $(z_c-z)^{a-1}$, leads to
asymptotics of the form $Z^\partial_{O(N),k}(\beta,x)\sim C k^{-a}$
for some finite $C>0$; see, for instance, \cite{FlSe09}, Corollary
VI.1. Of course, this requires additional hypotheses on
$Z^\partial_{O(N)}(z)$, which we ignore in this informal discussion.

We now summarize, and attempt to translate into a language more
familiar to
mathematicians, the discussion that can be found in \cite{duplantier04}, Appendix
B (see, in particular, equations B.48, B.64 and B.78,
and the
discussion at the end of Section~B.1.1 in \cite{duplantier04}). Assume
that $N\in(0,2)$ is
written in the form $N=2\cos(\pi\theta)$, where $\theta\in
(0,1/2)$. One conjectures that there exists a function $x_c(\beta)>0$
and a
critical value
$\beta_c>0$ such that:
\begin{itemize}
\item for fixed $\beta>\beta_c$ and $x=x_c(\beta)$,
\[
Z^\partial_{O(N)}(z)\mathop{\approx}_{z\uparrow z_c}
(z_c-z)^{1-\theta} ;
\]
\item for $\beta=\beta_c$ and $x=x_c(\beta_c)$,
\[
Z^\partial_{O(N)}(z)\mathop{\approx}_{z\uparrow z_c}
(z_c-z)^{1+\theta} .
\]
\end{itemize}
These two different behaviors, called the \textit{dense} and
the \textit{dilute} phase, respectively, hint at the asymptotics
\[
Z^\partial_{O(N),k}(\beta,x)\mathop{\sim}_{k\to\infty} C k^{-a}
,
\]
with $a=2-\theta$ and $a=2+\theta$, respectively. Recalling Section
\ref{sec:choos-boltzm-weights} and the preceding formula for $q_k$,
we see that the scaling limits of the distribution $W_{O(N)}$
in (\ref{model-O(N)}) should be related to the model studied in
the previous sections, with the particular value
$\alpha=a-1/2\in
\{3/2-\theta,3/2+\theta\}$. Note that the case $N=2$ appears as
a limiting critical situation where the dense and dilute phases
should coincide.

\begin{figure}

\includegraphics{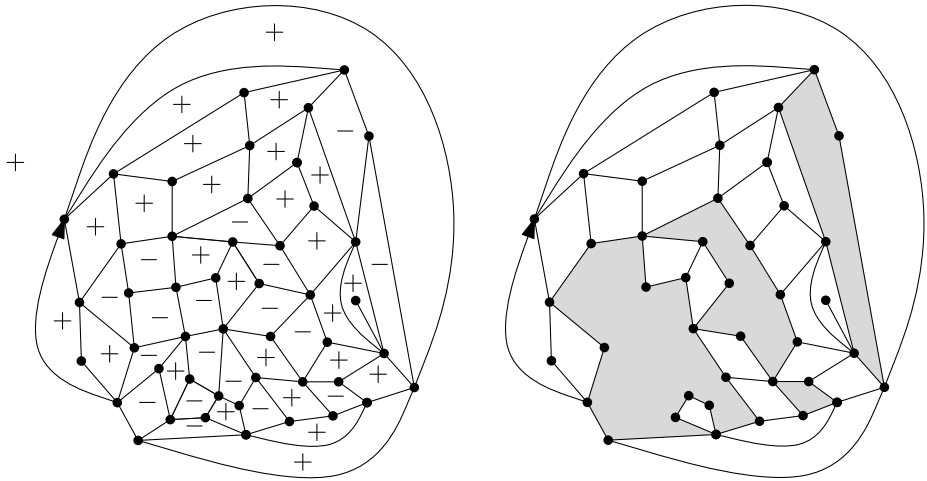}

\caption{An Ising (or percolation) configuration and the associated
exterior gasket.}
\label{fig:gasket}
\end{figure}

A similar description applies to other familiar statistical
physics models such as percolation or the Ising model on faces of a
random quadrangulation. In the latter setting, a configuration is a
pair $(\bq,\sigma),$
where $\bq$ is a rooted quadrangulation and
\[
\sigma=\bigl(\sigma_f,f\in
F(\bq)\bigr)\in\{-1,+1\}^{F(\bq)}.
\]
In the (annealed) Ising model, one
chooses the configuration with probability proportional to
\[
W_I(\bq,\sigma)=e^{-\beta\#F(\bq)} \exp\biggl(J\sum_{f\sim
f'}\sigma_f\sigma_{f'} \biggr) ,
\]
where $J$ is a real parameter and
the last sum is over all pairs of adjacent faces $f,f'$ in $\bq$. For
$J=0$, one gets the percolation model, where conditionally on the
quadrangulation $\bq$, all $\sigma\in\{-1,+1\}^{\#F(\bq)}$ are
equally likely to occur. One then defines the exterior gasket in a
way that should be clear from Figure \ref{fig:gasket}. This gasket
again has a
Boltzmann-type distribution when $(\bq,\sigma)$ is distributed
according to $W_I$. As previously, the relevant Boltzmann
weights correspond to partition functions for
the Ising model on a quadrangulation with a boundary. On the other
hand, the
topological constraints on the gaskets are now different: the boundaries
of holes need not be cycles and do not have
to be disjoint (however, an edge can be incident to at most one hole
and is
incident only once to this hole); see Figure \ref{fig:gasket}.

%

Kazakov \cite{kazakov06} identifies the value $J_c=\ln2$ as
critical. One conjectures that, respectively, for $J=J_c$ and $0\leq
J<J_c$ (and with the appropriate values of $\beta$), the Ising model
has the
same scaling limit as the dilute and dense phases of the $O(N=1)$
model, corresponding to $\theta=1/3$ and $\alpha\in
\{11/6,7/6\}$. This is confirmed (for $J=J_c$) by predictions for the
partition function of the Ising model with a boundary; see, for
example, Section 3.3 of \cite{COT96}.

Note that a discussion parallel to the present one appears in Sheffield
\cite{sheffield06}, Section 2.3, in
the case of regular hexagonal lattices, where it is conjectured that
the external gasket of
$O(N)$ models should converge to the so-called \textit{conformal loop
ensembles}, which are a conformally invariant family of random
curves related to the Schramm-Loewner evolutions. Such parallel
discussions might open some paths in
the mathematical understanding of the so-called KPZ formula, which
links scaling exponents for models on random and on regular
lattices. This approach would still be different from the
one developed recently by Duplantier and Sheffield \cite{DuSh08} as we
are focusing more
on the metric aspects of planar maps, rather than on the conformal
invariance properties that are intrinsic to \cite{DuSh08}.

At a rigorous level, it seems plausible that the topological
constraints that appear in
the random maps considered above can be handled using bijective methods,
in the spirit of Section \ref{sec:mobiles-bouttier-di}. Establishing
rigorous grounds for the conjectured behavior of
$Z^\partial_{O(N)}$ is another, probably much more challenging,
problem that would require a better understanding of the
combinatorial aspects of the $O(N)$ model on maps.


%
\printaddresses

\end{document}